\documentclass[a4paper, 12pt] {article}

\usepackage{fender}  

\begin{document}

    \title{RCF\,2 \\
           Evaluation and Consistency\footnote{
           Consideration of \emph{implicational} version $(\pidO)$ of
           \emph{Descent} axiom added} \\  
           $\eps\&\mathcal{C}*\piOR*\pidOR$}



\footnotetext{
Legend of LOGO: $\eps$ for Constructive evaluation, $\mathcal{C}$ 
for \emph{Self-Consistency} to be derived for suitable theories
$\piOR,\ \pidOR$ strengthening in a ``mild'' way the (categorical) 
Free-Variables Theory $\PRa$ of Primitive Recursion with predicate 
abstraction}

    \author{Michael Pfender\footnote{
              TU Berlin, Mathematik, pfender@math.tu-berlin.de}
           }    

    \date{July 2008\footnote{last revised \today}}                          

\maketitle

\textbf{Abstract:}
We construct here an \emph{iterative evaluation} of all (coded)
PR maps: progress of this iteration can be measured by 
\emph{descending complexity,} within Ordinal $O :\,= \N[\omega],$
of polynomials in one \emph{indeterminate,} called ``$\omega$''.
As (well) order on this Ordinal we choose the lexicographical one. 
Non-infinit descent of such iterations is added as a mild 
additional axiom schema $(\piO)$ to Theory $\PRa = \PR+(\abstr)$
of Primitive Recursion with \emph{predicate abstraction,} out of 
foregoing part RFC 1. This then gives (correct) \emph{on}-termination 
of iterative evaluation of \emph{argumented deduction trees} as well: 
for theories $\PRa$ and $\piOR = \PRa+(\piO).$ By means of this 
\emph{constructive} evaluation the \textbf{Main Theorem} is proved, 
on \emph{Termination-conditioned} \emph{(Inner) Soundness} for Theories 
$\piOR,$ $O$ extending $\N[\omega].$ As a consequence we get in fact 
\emph{Self-Consistency} for theories $\piOR,$ namely 
$\piOR$-\ul{derivabilit}y of $\piOR$'s own free-variable 
\emph{Consistency} formula 

\smallskip
  $\Con_{\piOR} = \Con_{\piOR} (k) \defeq 
      \neg\,\Pro_{\piOR} (k,\code{\false}): 
                           \N \to 2,\ k \in \N\ \free.$

\smallskip \noindent
Here PR predicate $\Pro_\T (k,u)$ says, for an arithmetical theory 
$\T:$ number $k \in \N$ is a $\T$-\emph{Proof} code \emph{proving} 
internally $\T$-\emph{formula} code $u,$ arithmetised \emph{Proof} 
in G\"odel's sense.

As to expect from classiccal setting, Self-Consistency of $\piOR$
gives (unconditioned) Objective Soundness. Eventually we show 
\emph{Termination-Conditioned} Soundness ``already'' for $\PRa.$
But it turns out that \emph{present} derivation of Self-Consistency,
and already that of \emph{Consistency formula} of $\PRa$ 
from this \emph{conditioned} Soundness ``needs'' schema $(\tilde{\pi})$ 
of \emph{non-infinit descent} in Ordinal $\N[\omega],$ which is
presumably not derived by $\PRa$ itself.

 
\newpage

\section{Summary}

\footnotetext{extended Poster Abstract  
  ``Arithmetical Consistency via Constructive evaluation'',
  Conference celebrating Kurt G\"odel's 100th birthday, 
  Vienna april 28, 29, 2006}

G\"odel's first Incompleteness Theorem for \emph{Principia Mathematica}
and \emph{``verwandte Systeme'',} on which in particular is based the 
second one, on non-prova\emph{bility} of $\PM$'s own 
\emph{Consistency formula} $\Con_{\PM},$ exhibits a (closed) $\PM$ 
formula $\ph$ with property that 
  $$\PM \derives\ [\,\ph 
     \iff \neg\,(\exists\,k \in \N)\,\Pro_{\PM} (k,\code{\ph})\,],
                                            \ \text{in words:}$$
Theory $\PM$ \ul{derives} $\ph$ to be equivalent to its
``own'' \emph{coded, arithmetised non-Provability.}  

Since this equivalence needs already for its \emph{statement} ``full'' 
formal, \emph{``not testable''} quantification, the 
\emph{Consistency Provability} issue is not settled for Free-Variables 
Primitive Recursive Arithmetic and its strengthenings -- Theories $\T$ 
which express (formalised, ``internal'') Consistency as free-variable 
formula 
  $$\Con_\T = \Con_\T(k) = \neg\,\Pro_\T(k,\code{\false}): \N \to 2:$$
``No $k \in \N$ is a \emph{Proof} code \emph{proving} $\code{\false}.$''

This is the point of depart for investigation of ``suitable'' 
strengthenings $\piOR = \PRa+(\piO)$ of categorical Theory $\PRa$ 
of Primitive Recursion, enriched with 
\emph{predicate abstraction Objects} 
$\set{A\,|\,\chi} = \set{a \in A\,|\,\chi(a)}:$ Plausibel axiom schema
$(\piO),$ more presisely: its contraposition $\tilde{\pi}_O,$ 
states ``weak'' impossibility of infinite descending chains in
any \emph{Ordinal} $O$ extending polynomial semiring $\N[\omega],$ 
with its canonical, \emph{lexicographical} order.

\smallskip
\textbf{Central \emph{Non-Infinite Descent} Schema,\ \emph{Descent} Schema}
for short:

We need an  \textbf{axiom-schema} for expressing 
-- in \emph{free variables} -- \textbf{Finite descent} 
\textbf{(endo-driven) chains,} \emph{descending} in 
\emph{complexity value} out of Ordinal 
$O \succeq \N[\omega],$ a schema called $(\piO),$
which gives the ``name'' to 
\emph{Descent}\footnote{notion added 2 JAN 2009} Theory 
$\piOR = \PRa+(\piO):$
This theory is a \emph{pure strengthening} of $\PRa,$ 
it has the same \emph{language.}   

Easier to interprete logically is $(\piO)$'s equivalent, 
\emph{Free-Variables contraposition,} on ``absurdity'' of 
\emph{infinite descending chains,} namely:
\inference{ (\tilde{\pi}_O) }
{ $c = c(a): A \to O$ PR \emph{(complexity),} \\
& $p = p(a): A \to A$ PR \emph{(predecessor endo),} \\
& $\PRa\ \derives\ c(a) > 0_O \implies c\,p(a) < c(a)$ \emph{(descent),} \\
& $\PRa\ \derives\ c(a) \doteq 0_O \implies p(a) \doteq a$
                                       \emph{(stationarity at zero)} \\
& $\psi = \psi(a): A \to 2$ \emph{absurdity test predicate,} \\
& $\PRa \derives\ \psi(a) \implies c\,p^n(a) > 0_O,$ \\
& \qquad 
    with quantifier decoration: \\
& $\PRa \derives\ \forall\,a\,[\,\psi(a) 
               \implies \forall\,n\,c\,p^n(a) > 0_O\,]$ \\
& \qquad
    the \emph{latter statement:} ``infinit descent'', is felt absurd, \\
& \qquad    
    and ``therefore'' so ``must be'', by \emph{axiom,} \\
& \qquad 
    \emph{condition} $\psi$ \emph{implying} this ``absurdity'':
}  
{ $\piOR \derives\ \psi(a) \doteq \false: A \to 2,$ intuitively: \\
& $\piOR \derives\ \forall\,a\ \neg\,\psi(a).$
}           

\smallskip
$[\,$The first four lines of the \emph{antecedent} constitute $(p,c)$ 
as (the data of) a $\mrCCIO:$ of a \emph{Complexity Controlled Iteration,} 
with (stepwise) descending order values in \emph{Ordinal} $O.$ 
Central \textbf{example:} \emph{General Recursive,} \NAME{Acckermann} 
type \emph{PR-code evaluation} $\eps$ will be \emph{resolved} into 
such a $\mrCCIO,$ $O :\,= \N[\omega] \subset \N.]$ 


\smallskip
My \textbf{Thesis} then is that these theories $\piOR,$ weaker 
than $\PM,$ \textbf{set theories} and even Peano Arithmetic $\PA$
(when given its \emph{quantified} form), \ul{derive} their own 
internal (Free-Variable) \emph{Consistency formula} 
$\Con_{\piOR} (k): \N \to 2,$ see above. 

\smallskip
\textbf{Notions and Arguments for Self-Consistency of} $\piOR:$ 
In order to obtain \emph{constructive} Theories -- candidates for
\emph{self-Consistency} -- we introduce first, into \emph{fundamental}
Theory $\PR$ of (categorical) \emph{Free-Variables} Primitive 
Recursion, \emph{predicate abstraction} of PR maps 
$\chi = \chi(a): A \to 2$ ($A$ a finite power of NNO $\N$),
into \emph{defined Objects} $\set{A\,|\,\chi},$  
and then \emph{strengthen} Theory $\PRa$ obtained this way,
by a free-variables, \emph{(inferential)} schema $(\piO)$ of 
\emph{``on''-terminating descent,} into Theorie(s) 
$\piOR,$ \emph{on}-terminating descent of 
\emph{Complexity Controlled Iterations} ($\mrCCIO$'s, see above),
with (descending) complexity values in \emph{Ordinal} 
$O \succeq \N[\omega].$  

Strengthened Theory $\piOR = \PRa+(\piO),$ with its 
\emph{language} equal to that of $\PRa,$ is asserted to derive 
the (Free-Variable) formula $\Con_{\piOR} (k)$ which expresses
internally: within $\piOR$ itself, \emph{Consistency} of Theory 
$\piOR,$ see above. 

\textbf{Proof} is by $\mrCCI_{\N[\omega]}$ (descent) property of a 
suitable, \emph{atomic} PR evaluation \emph{step} $e$ applied to 
\emph{PR-map-code/argument} pairs 
$(u,x) \in \mrPRa \times \X.$ 

\smallskip
$[\,$Here $\X \subset \N$ denotes the \emph{Universal Object} 
of all (codes of) \emph{singletons} and (nested) \emph{pairs} 
of natural numbers, enriched by a shymbol $\bott$ equally
coded in $\N,$ to designate \emph{undefined values,} of
\emph{defined partially defined} PR maps. Objects $A$ of 
$\PRa,\ \piOR$ admit a \emph{natural embedding} $A \sqsub \X$ 
into this this universal Object$.]$

\smallskip
Iteration $\eps,$ of step $e,$ is in fact \emph{controlled} by a 
\emph{syntactic complexity} $c_{\mrPR} (u) \in \N[\omega],$ 
descending with each application of $e$ as long as minimum complexity 
$0 = c_{\mrPR} (\code{\id})$ is not ``yet'' reached.

\emph{Strengthening} of $\PRa$ by schema $(\piO)$ 
-- \cf its free-variables contraposition $(\tilde{\pi}_O)$ above -- 
into Theory $\piOR = \PR+(\piO),$ is ``just'' to allow for 
a so to say \emph{sound,} canonical evaluation ``algorithm'' 
for $\piOR:$

On one hand it is proved straight forward that evaluation $\eps$ 
above has the expected recursive properties of an \emph{evaluation,} 
this within (categorical, Free-Variables) Theory $\muR$ of 
$\mu$-Recursion. 

On the other hand, $\piOR$ has the same \textbf{Language} as $\PRa,$ 
so that this $\eps$ is a natural candidate for likewise 
-- \emph{sound} -- evaluation of internal version of theory $\piOR,$ 
and for being \emph{totally defined} in a suitable \emph{Free-Variables} 
sense, technically: to \emph{on-terminate,} this just by its property 
to be a \emph{Complexity Controlled Iteration,} with order values in 
$\N[\omega].$

In fact, by schema $(\piO)$ itself ($O$ extending $\N[\omega]$), 
$\eps$ \emph{preserves} the \textbf{extra} equation instances 
inserted by internalisation of $(\piO).$

\smallskip
\textbf{Dangerous bound:} is there a good reason that this evaluation 
is not a \emph{self-evaluation} for Theory $\piOR?$

Answer: $\eps$ is -- by definition -- \emph{not} PR: 
If you take the \emph{diagonal} 
  $$\diag(n) \defeq \eps(\mr{enum}_{\mrPR} (n),\cantor_\X(n)): 
                                                          \N \to \N,$$
$\mr{enum}_{\mrPR}$ an internal PR \emph{count} of all PR map codes, and
$\cantor_\X: \N \overset{\iso} {\lto} \X$ ``the'' Cantor's 
\emph{count} of $\X \subset \N,$  
then you get \NAME{Ackermann's} 
original diagonal function\footnote{
    for a two-parameter, simple genuine \NAME{Ackermann} function
    \cf Eilenberg/Elgot 1970}
which grows faster than any PR function: but $\piOR$ has only 
$\mrPR$ maps as its \emph{maps,} it is a (pure) \emph{strengthening} 
of $\PRa.$

On the other hand, $\eps$ is \emph{intuitively} total, 
since, intuitively, complexity $c\ e^m(u,x)$ 
``must'' reach $0$ in \emph{finitely many} $e$-steps. The latter 
intuition can be, in free variables (!), expressed \emph{formally} 
by $\piOR$'s \textbf{schema} $(\tilde{\pi}_O):$ Free-Variables 
contraposition of $(\piO).$ Schema $(\tilde{\pi}_O)$ says 
that a condition which implies \emph{infinite descent} of such 
a chain (on all $x$), must be \emph{false} (on all $x$), ``absurd''.    

\smallskip
\textbf{Complexity Controlled Iteration} $\eps$ of $e$ extends 
canonically into a Complexity Controlled evaluation $\eps_d,$ 
of \textbf{argumented deduction trees,} $\eps_d$ again defined 
by $\mrCCI_{\N[\omega]}:$ this time by iteration of a 
\emph{tree evaluation step} $e_d$ suitably extending basic 
evaluation step $e$ to argumented deduction trees. 
 
Deduction-tree evaluation starts on trees of form $\dtree_k/x,$ 
obtained as follows from $k$ and $x:$ Call $\dtree_k$ the (first) 
\emph{deduction tree} which (internally) \emph{proves} $\,k$\,th 
internal equation $u\,\checkeq_k\,v$ of theory $\piOR,$ enumeration of 
\emph{proved} equations being (lexicographically) by code of (first) 
\emph{Proof.} 
This argument-free deduction tree $\dtree_k$ then is provided 
-- node-wise top down from given $x \in \X$ -- 
with its \emph{spread down} arguments in 
$\Xb \defeq \X\,\dot\cup\,\set{\Box} 
                        = \X\,\dot\cup\,\set{\an{}} \subset \N;$   
(empty list $\Box = \an{}$ refers to a not yet known argument, not 
``yet'' at a given time of stepwise \emph{evaluation} $e_d$.) 

\emph{Spreading down} arguments this way eventually converts 
argument-free $k$\,th deduction tree $\dtree_k$ into 
(partially non-dummy) \emph{argumented deduction tree} $\dtree_k/x.$

\emph{Iteration} $\eps_d,$ of tree evaluation step $e_d,$ again is 
\emph{Complexity Controlled descending} in Ordinal $\N[\omega],$ 
when controlled by deduction tree \emph{complexity} $c_d.$ 
This complexity is defined essentially as the (polynomial) \emph{sum} 
of all (syntactical) complexities $c_{\mrPR} (u)$ of \emph{map codes} 
appearing in the deduction tree. 

So, as it does to \emph{basic} evaluation $\eps,$ schema 
$\tilde{\pi}_{\N[\omega]}$ applies to complexity controlled evaluation 
$\eps_d$ of argumented deduction-trees as well, and gives

\smallskip
\textbf{Deduction-Tree Evaluation non-infinit Descent:} Infinit strict
descent of endo map $e_d$ -- with respect to complexity $c_d$ 
-- is \emph{absurd.} 

\smallskip
This deduction-tree evaluation $\eps_d$ externalises, 
\emph{as far as terminating,} $k$\,th internal equation 
$u\,\checkeq_k\,v$ of theory $\piOR$ into \emph{complete evaluation} 
$\eps(u,x) \doteq \eps(v,x):$
  
\emph{Termination-Conditioned Inner Soundness,} our \textbf{Main Theorem.} 

\smallskip  
For a given PR predicate $\chi = \chi(x): \X \to 2,$ the 
\textbf{Main Theorem} reads:

\medskip
Theory $\piOR$ \ul{derives}: \textbf{If} for $k \in \N$ and for 
$x \in \X \sminus \set{\bott}$ given, $\Pro_{\piOR} (k,\code{\chi})$ 
``holds'', and \textbf{if} \emph{argumented $\piOR$ deduction tree} 
$\dtree_k/x$ admits \emph{complete evaluation} by $m$ (``say'') 
deduction-tree evaluation-steps $e_d,$ 

\textbf{Then} the pair $(k,x)$ is a \textbf{Soundness}-Instance, 
\ie \textbf{then} $k$\,th given (internal) $\piOR$-\emph{Provability}
$\Pro_{\piOR} (k,\code{\chi})$ \emph{implies} $\chi(x),$ 
for the given argument $x \in \X \sminus \set{\bott}.$ All this within 
Theory $\piOR$ itself.  

\medskip
\textbf{Corollary: Self-Consistency Derivability for Theory} $\piOR:$
\begin{align*}
\piOR \derives\ 
& \Con_{\piOR}, 
    \ \text{\ie ``necessarily'' in \emph{Free-Variables} form:} \\
\piOR \derives\ 
& \neg\,\Pro_{\piOR} (k,\code{\false}): \N \to 2,
                                 \ \text{\ie equationally:} \\  
\piOR \derives\ 
& \neg\,[\,\code{\false} \checkeq_k \code{\true}\,]: 
                                   \N \to 2,\ k \in \N \ \free:
\end{align*}
\emph{Theory $\piOR$ \ul{derives} that no $k \in \N$ is the 
internal $\piOR$-\emph{Proof} for $\code{\false}.$}
 
\smallskip
\textbf{Proof} of this \textbf{Corollary} to 
\emph{Termination-Conditioned Soundness:} 

By the last assertion of the 
\textbf{Theorem,} with $\chi = \chi(x) :\,= \false_\X (x): \X \to 2,$ 
and $x :\,= \an{0} \in \X,$ we get: 
  
\smallskip
\emph{Evaluation-effective internal inconsistency} of $\piOR,$ 
\ie availability of an \emph{evaluation-terminating} internal 
\emph{deduction tree} of $\code{\false},$ \emph{implies} $\false:$ 
  $$\piOR \derives\ 
      \code{\false} \checkeq_k \code{\true}\ 
        \land\ c_d\ e_d^m(\dtree_k/\an{0}) \doteq 0 
                             \implies \false_\X(\an{0}).$$
Contraposition to this, still with $k,m \in \N$ free:
  $$\piOR \derives\
      \true_\X(\an{0}) \implies 
        \neg\,[\,\code{\false} \checkeq_k \code{\true}\,]
                           \,\lor\,c_d\ e_d^m(\dtree_k/\an{0}) > 0,$$
\ie by Free-Variables (Boolean) tautology: 
  $$\piOR \derives\
      \code{\false} \checkeq_k \code{\true} 
        \implies c_d\ e_d^m(\dtree_k/\an{0}) > 0: \N^2 \to 2.$$
This $\piOR$ \ul{derivative} invites to apply schema 
$(\tilde{\pi}_{\N[\omega]})$ of $\piOR:$  
 
\smallskip
\emph{``infinite endo-driven descent with order values 
in $\N[\omega]$ is absurd.''}
 
\smallskip
We apply this schema to deduction tree evaluation $\eps_d$ given
by \emph{step} $e_d$ and complexity $c_d$ which descends 
-- this is \emph{Argumented-Tree Evaluation Descent} -- 
with each application of $e_d,$ as long as complexity $0$ is not 
(``yet'') reached. We combine this with choice of 
``overall'' \emph{absurdity condition}  
  $$\psi = \psi(k) :\,= [\,\code{\false} \checkeq_k \code{\true}\,]: 
                                      \N \to 2,\ k \in \N\ \free\ (!)$$
and get, by schema $(\tilde{\pi}_{\N[\omega]}),$ overall negation of this 
(overall) ``absurd'' predicate $\psi,$ namely
  $$\piOR \derives\ \neg\,[\,\code{\false} \checkeq_k \code{\true}\,]: 
                                     \N \to 2,\ k \in \N\ \text{free.}$$
This is $\piOR$-\ul{derivation} of the \emph{free-variable} 
\emph{Consistency Formula} of $\piOR$ itself.


\smallskip
From this \emph{Self-Consistency} of Theorie(s) $\piOR,$ which is
equivalent to \emph{injectivity} of (special) internal 
\emph{numeralisation} $\nu_2: 2 \into [\one,2]_{\piOR}\,,$ we get 
immediately injectivity of \emph{all} these numeralisations 
  $\nu_A = \nu_A(a): A \into [\one,A] = \cds{\one,A}/\checkeq\,,$
and from this, with \emph{naturality} of this family, ``full'' 
objective \textbf{Soundness} of Theory $\piOR$ which reads:

\smallskip
\emph{Formalised} $\piOR$-\emph{Provability} of (code of) 
PR predicate $\chi: \X \to 2$ \emph{implies} -- within Theory 
$\piOR$ -- \emph{``validity''} $\chi(x)$ of $\chi$ at ``each'' 
of $\chi$'s arguments $x \in \X.$


But for derivation of \emph{Self-Consistency} from Termination-conditioned 
Soundness, a suitable \textbf{strengthening} of $\PRa,$ here by schema 
$(\tilde{\pi}) = (\tilde{\pi}_{\N[\omega]}),$ stating \emph{absurdity} 
of infinite descent in Ordinal $\N[\omega],$ seems to be necessary:
my guess is that Theories $\PRA$ as well as $\PR$ and hence $\PRa,$ 
are \emph{not strong enough} to derive their own (internal) Consistency. 
On the other hand, we know from G\"odel's work that Principia Mathematica
``und verwandte Systeme'' are \emph{too strong} for being self-consistent. 
This is true for any (formally) \emph{quantified} Arithmetical Theory 
$\bfQ,$ in particular for the (classical, quantified) version $\PA$ of 
Peano Arithmetic: Such theory $\bfQ$ has all ingredients for G\"odel's 
Proof of his two \emph{Incompleteness Theorems.}

\smallskip
In section 7 We discuss\footnote{insertion ? JAN 2009} 
a formally stronger, \emph{implicational,} ``local'' variant $(\pidO)$ 
of \ul{inferential} \emph{Descent} axiom $(\piO),$ with respect 
to \emph{Self-Consistency} and (Objective) \emph{Soundness:} 
In particular, \emph{Self-Consistency} \textbf{Proof} becomes 
technically easier for corresponding theory $\pidOR.$

\smallskip
The final section 8\footnote{inserted 2 JAN 2009} 
gives a \textbf{proof} of (Objective) \ul{Consistenc}y
for Theorie(s) $\pidOR$ (hence $\piOR$) relative to basic Theory 
$\PRa$ of Primitive Recursion and hence relative to fundamental Theory
$\PR$ of Primitive Recursion ``itself''.  

For \textbf{proof} of this (relative) Consistency, we use a schema, 
$(\rho_O),$ of recursive \emph{reduction} for predicate validity,
reduction along a Complexity Controlled Iteration ($\mrCCIO$), 
admitted by Theory $\PRa$ (and its strengthenings.)






\section{Iterative Evaluation of PR Map Codes}

Object- and map terms of all our theories are coded straight ahead, 
in particular since formally we have no (individual) \emph{variables} 
on the Object Language level: We code all our terms just as prime power 
products ``over'' the \LaTeX source codes describing these terms,
this externally in \ul{naive numbers}, out of $\ul{\N}$ as well as 
into the NNO $\N$ of the (categorical) arithmetical theory itself.


\smallskip
\textbf{Equality Enumeration:} As ``any'' theories, \emph{fundamental}
Theory $\PR$ of Primitive Recursion as well as \emph{basic} Theory
$\PRa = \PR+(\abstr),$ definitional enrichement of $\PR$ by the schema 
of \emph{predicate abstraction:} 
$\bfan{\chi: A \to 2} \bs{\mapsto} \set{A\,|\,\chi},$ a ``virtual'',
\emph{abstracted} Object in $\PRa,$ admit an (external) primitive 
recursive \ul{enumeration} of their respective \textbf{theorems,} 
ordered by length (more precisely: by lexicographical order) 
of the first \textbf{proofs} of these (equational) Theorems, here:
\begin{align*}
& =^{\PR}(\ul{k}): \ul{\N} \bs{\to} \PR \bs{\times} \PR 
             \bs{\subset} \ul{\N} \bs{\times} \ul{\N}\ \text{and} \\
& =^{\PRa}(\ul{k}): \ul{\N} \bs{\to} \PRa \bs{\times} \PRa 
                      \bs{\subset} \ul{\N} \bs{\times} \ul{\N}
\end{align*}
respectively.

\smallskip
By the PR Representation Theorem 5.3 of \NAME{Rom\`an} 1989, these 
\ul{enumerations} give rise to their internal versions
\begin{align*}
& \checkeq^{\PR}_k: \N \to \mrPR \times \mrPR \subset \N^2\ \text{and} \\
& \checkeq^{\PRa}_k: \N \to \mrPRa \times \mrPRa \subset \N^2,
\end{align*}
with internalisation \emph{(representation)} property
\begin{align*}
& \PR \derives\ \checkeq_{\num(\ul{k})}\ =\ \num(=^{\PR}_{\ul{k}}):
                  \one \to \mrPR \times \mrPR \subset \N^2\ \text{and} \\
& \PR \derives\ \checkeq_{\num(\ul{k})}\ =\ \num(=^{\PRa}_{\ul{k}}):
                  \one \to \mrPRa \times \mrPRa \subset \N^2.
\end{align*}
Here (external) \ul{numeralisation} is given externally PR as 
\begin{align*}
& \num(\ul{n}) = s^{\ul{n}}: 
                    \one \ovs{0} \N \ovs{s} \ldots \ovs{s} \N, \\
& \num(\ul{m},\ul{n}) = (\num(\ul{m}),\num(\ul{n})): \one \to \N \times \N,
                 \ \ul{m},\ul{n}\ \text{(``meta'') \ul{free} in}\ \ul{\N},
\end{align*} 
$\mrPR = \set{\N\,|\,\mrPR}$ is the predicative, PR decidable
subset of $\N$ ``of all $\PR$ codes'' (a $\PRa$-Object), 
\emph{internalisation} of $\PR \bs{\subset} \ul{\N}$ of all 
$\PR$-\ul{terms} on Object Language level. Analogeous meaning for
\emph{internalisation} $\mrPRa \subset \N$ of $\PRa \bs{\subset} \ul{\N}.$


\smallskip
For discussion of ``constructive'' evaluation, we need representation
of all $\PRa$ maps within one $\PR$ endo map \ul{monoid}, namely
within $\PR(\Xbott,\Xbott),$ where $\X \subset \N,$ 
$\X = \set{\N\,|\,\X: \N \to 2}$ is the (predicative)
\emph{Universal Object} of $\N$-\emph{singletons} 
$\set{\an{n}\,|\,n \in \N},$ possibly nested $\N$-\emph{pairs}
$\set{\an{a;b}\,|\,a,b \in \X},$ and 
  $$\Xbott \defeq \X\,\dot\cup\,\set{\bott} 
              = \X(a)\ \dot\lor\ a \doteq \bott: \N \to 2$$
is $X$ augmented by symbol (code) $\bott: \one \to \N,$ 
$\bott$ taking care of \emph{defined undefined} 
arguments of \emph{defined partial maps.}\footnote{
                        \cf Ch.\ 1, final section $\X$} 

Here we view (formally) 
$\X = \X(a),\ \Xbott = \Xbott(a): \N \to \N$ as $\PR$-\emph{predicates,} 
not ``yet'' as \emph{abstracted} Objects $\X = \set{\N\,|\,\X},$ 
$\Xbott = \set{\N\,|\,\Xbott},$ of Theory $\PRa = \PR+(\abstr).$

We allow us to write ``$a \in \X$'' instead of
$\X(a) \doteq \true: \N \to \N,$ and ``$a \in \Xbott$'' for
$\Xbott(a) \doteq \true,$ and similarly for other predicates. 

This way we introduce -- \`a la \NAME{Reiter} -- ``Object'' $2$
just as target for predicates $\chi: A \to 2,$ meaning
$\chi: A \to \N$ to be a \emph{predicate} in the exact sense
that $\chi: \A \to \N$ satisfies
  $$\chi \circ \sign \bydefeq \chi \circ \neg \circ \neg = \chi:
      \N \ovs{\chi} \N \ovs{\sign} \N,
                            \ \text{``still''}\ A\ \emph{fundamental.}$$
We \textbf{define,} within endo map \ul{set} $\PR(\N,\N)$ a 
subTheory $\PRX$ externally PR as follows, by mimikry of schema $(\abstr)$
for the special case of predicate $\X = \X(a): \N \to \N,$ but
\emph{without} introduction of a coarser notion of equality, as in
case of schema of abstraction constituting Theory 
$\PRa = \PR+(\abstr).$

So Theory $\PRX \bs{\subset} \PR(\N,\N)$ comes in, by external PR
\ul{enumeration} of its Object and map \ul{terms} as follows:

\emph{Objects} of $\PRX$ are \emph{predicates} 
$\chi: \X \to 2,$ \ie $\PR$-predicates $\chi: \N \to 2$ such that 
\begin{align*}
& \PR \derives\ \chi(a) \implies \X(a): \N \to 2,
                                   \ \text{\ie such that} \\
& \PR \derives \chi(a) \implies \Xbott(a)\,\land\,a \neq \bott: 
                                                       \N \to 2.
\end{align*}
$\PRX$-maps in $\PRX(\chi,\psi)$ are $\PR$-maps $f: \N \to \N$ 
such that 
  $$\neg\,\X(a) \implies f(a) \doteq \bott,\ \text{and}\ 
      \chi(a) \implies \psi \circ f(a): \N \to 2,$$
observe the ``truncated'' parallelism to \textbf{definition} of 
$\PRa$-maps \\
$f: \set{A\,|\,\chi} \to \set{B\,|\psi}.$

\smallskip
Then ``assignment'' $\bfI: \PR \overset{\bs{\sqsub}} {\To} \PRX$ 
is \textbf{defined} as follows externally PR:
\begin{align*}
& \bfI\,\one = \dot{\one} \defeq \set{\an{0}}: 
                        \N \supset \Xbott \supset \X \to 2, \\
& \bfI\,\N = \dot{\N} 
    \defeq \an{\N} \defeq \set{\an{n}\,|\,n \in \N}: 
                        \N \supset \Xbott \supset \X \to 2, \\
& \quad
    \text{and further \ul{recursivel}y:} \\
& \bfI\,(A \times B) \defeq \an{A \times B} 
    \defeq \set{\an{a;b}\,|\,(a,b) \in (A \times B)}: 
                                    \N \supset \X \to 2,
\end{align*}
Functorial \textbf{definition} of $\bfI$ on $\PR$ maps:
  $$\PR(A,B) \owns f \overset{\bfI} 
               {\bs{\mapsto}} \bfI\,f = \dot{f} \in \PRX$$
then is ``canonical'', by external PR on the structure of
$\PR$-map $f: A \to B,$ in particualar by mapping all ``arguments'' 
in $\N \sminus \dot{A} = \N \sminus \bfI\,A$ into 
$\bott \in \Xbott:$ one \emph{waste basket} outside all Objects of 
$\PRX.$\footnote{for the details see Ch.\ 1, final section $\X.$} 

\smallskip
Interesting now is that we can extend embedding $\bfI$ above
into an \ul{embeddin}g $\bfI: \PRa \To \PRX,$ by the following

\smallskip
\textbf{Definition:} For a (general) $\PRa$ Object, of form 
$\set{A\,|\,\chi},$ define
\begin{align*}
& \bfI\,\set{A\,|\,\chi} \defeq \set{\dot{A}\,|\,\dot\chi}
                           \bydefeq \set{\bfI\,A\,|\,\bfI\,\chi} \\
& \bydefeq \set{a \in \bfI\,A\,|\,\bfI\,\chi(a) \doteq \an{\true}}:
                                              \N \supset \Xbott \to 2.
\end{align*}
We replace here ``don't-worry arguments'' in the complement $\neg\,\chi$
of $\PRa$-Object $\set{A\,|\,\chi}$ by \emph{cutting} them \emph{out}
in the definition of \emph{replacing} $\PRX$-Object 
$\bfI\,\set{A\,|\,\chi} = \set{\dot{A}\,|\,\dot\chi}.$ 
``Coarser'' notion $=^{\PRa}$ (coarser then $=^{\PR}$) is then replaced by
original notion of equality, $=^{\PR}$ itself, notion of
map-equality of \emph{roof} 
$\PRX\,\text{``$\bs{\subset}$''}\,\PR(\N,\N):$
This formal ``sameness'' of PR equality was the goal of the
considerations above: The new version $\PRaX$ replacing $\PRa$
isomorphically, is a \textbf{sub}Theory of $\PR$ with 
\emph{notion of equality} -- objectively as well as (then) 
\emph{internally} -- inherited from \emph{fundamental} Theory $\PR.$
 
\medskip
\textbf{Universal Embedding Theorem:}\footnote{from 
                                         Ch.\ 1, final section $\X$} 
\begin{enumerate} [(i)]
\item
$\bfI: \PR \To \PRX \bs{\subset} \PR(\N,\N)$ above is an \ul{embeddin}g 
which preserves composition.

\item
(Enumerative) Restriction 
$\bfI: \PR \overset{\bs{\iso}} {\To} \PR^\X \defeq\bfI\,[\,\PR\,]$ 
of this embedding to its (\ul{enumerated}) Image defines
an \ul{isomor}p\ul{hism} of categories. It is \textbf{defined} above as
  $$\bfan{f: A \to B} \overset{\bfI} {\To} 
                           \bfan{\dot{f}: \dot{A} \to \dot{B}},$$
by the ``natural'' (primitive) \ul{recursion} on the structure
of $f$ as a \ul{ma}p in fundamental Theory $\PR$ of (Cartesian)
Primitive Recursion.  

\item $\PR$ embedding $\bfI$ ``canonically'' extends into an 
\ul{embeddin}g (!)
  $$\bfI: \PRa \To \PR(\N,\N)$$ 
\ul{of} Theory $\PRa = \PR+(\abstr)$ -- Theory $\PR$ with \emph{abstraction}
of \emph{predicates into} (``new'', ``virtual'') \emph{Objects} 
$\set{A\,|\,\chi: A \to 2}$ -- \ul{to} the \ul{Set} of $\PR$ endomaps of
$\N,$ of which -- by the way -- $\PRa(\Xbott,\Xbott)$ is (formally) a 
\ul{Sub}Q\ul{uotient}.
 
\smallskip
$[\,$Equality $\ =^{\PRa}$ of (distinguished) $\PR$ endo maps when 
viewed as \\
$\PRa$ endo maps on $\Xbott = \set{\N\,|\,\Xbott: \N \to 2},$ 
is \ul{embedded} to \\
$\PRX$- ($\PR$-)equality by 
  $\bfI: \PRa \To \PRX\,\text{``$\bs{\subset}$''}\,\PR(\N,\N).]$

\item \textbf{Main} assertion: Embedding $\bfI$ above \textbf{defines} 
an \ul{isomor}p\ul{hism} of \ul{cate}g\ul{ories} 
  $$\bfI: \PRa \overset{\bs{\iso}} {\To} \PRaX$$
onto a ``naturally choosen'' (\ul{emumerated}) category $\PRaX$ 
of $\PR$ predicates on \emph{Universal Object} ($\PR$-predicate) 
$\Xbott: \N \to \N,$ with canonical maps in between (see above), and 
whith composition inherited from that of $\PR(\N,\N).$ This isomorphism 
is defined (naturally) by
\begin{align*}
& \bfI\,(\,f: \set{A\,|\,\chi} \to \set{B\,|\,\psi}\,)
    = \bfan{\dot{f}: \dot\chi \to \dot\psi}, \\
& \dot\chi: \N \supset \Xbott \supset \X \supset \dot{A} \to 2, \\
& \dot\psi: \N \supset \Xbott \supset \X \supset \dot{B} \to 2, \\
& \dot{f} \bydefeq \bfI_{\PR} (f): 
            \N \supset \dot{A} \to \dot{B} \subset \N\ \text{above.}
\end{align*}
By this isomorphism of categories, $\PRaX$ inherits from category
$\PRa$ all of its (categorically described) structure: the 
isomorphism transports Cartesian PR structure, equality predicates
on all Objects, schema of predicate abstraction, equalisers,
and -- trivially -- the whole algebraic, logic and order structure
on NNO $\N$ and truth Object $2.$ 
    
\smallskip
We have furthermore:

\item
For each fundamental Object $A,$ embedded Object 
$\dot{A} = \bfI\,A \subset \Xbott$ comes with a \emph{retraction} 
  $\mr{retr}^\X_A: \Xbott \to \dot{A}\,\dot\cup\,\set{\bott},$ 
\textbf{defined} by 
$\mr{retr}^\X_A (a) \defeq a$ for $a \in \dot{A},$ 
$\mr{retr}^\X_A (a) \defeq \bott$ \emph{otherwise.}

\smallskip
This family of retractions clearly extends to a retraction family 
  $$\mr{retr}^\X_{\set{A\,|\,\chi}}: 
      \Xbott \to \set{\dot{A}\,|\,\dot\chi} \,\dot\cup\,\set{\bott}
                 = \bfI\,\set{A\,|\,\chi}\,\dot\cup\,\set{\bott}$$
for all $\PRa$-Objects $\set{A\,|\,\chi}:$ This is what 
$\bott \in \Xbott$ is good for. 

\item 
For each Object $\set{A\,|\,\chi}$ of $\PRa,$ in particular
for each \emph{fundamental} Object $A \identic \set{A\,|\,\true_A},$ 
$\PRa$ comes with the characteristic (predicative) \emph{subset} 
$\dot\chi: \bfI\,\set{A\,|\,\chi}: \Xbott \to 2$ of $\Xbott$ 
\textbf{defined} PR above, isomorphic to $\set{A\,|\,\chi}$ within 
$\PRa$ (!) via ``canonical'' $\PRa$-isomorphism 
  $$\mr{iso}^\X_{\set{A\,|\,\chi}}: 
      \set{A\,|\,\chi} \ovs{\iso} \bfI\,\set{A\,|\,\chi} 
                              = \set{\dot{A}\,|\,\dot\chi},$$
the $\PRa$-isomorphism \textbf{defined} PR on the ``structure'' of 
$\set{A\,|\,\chi},$ as restriction of $\mr{iso}^\X_A: A \to \bfI\,A$ 
for \emph{fundamental} Object $A,$ in turn (externally/internally) 
PR defined by
\begin{align*}
& \mr{iso^\X_{\one}} (0) \defeq \an{0}: 
                          \one \to \bfI\,\one \subset \Xbott, \\
& \mr{iso^\X_\N} (0) \defeq \an{0}: 
    \one \to \ [\ \bfI\,\one \subset\ ]\ \bfI\,\N \subset \Xbott, \\
& \quad
    \text{further externally PR:} \\
& \mr{iso}^\X_{(A \times B)} (a,b) 
    \defeq \an{\mr{iso}^\X_A(a);\mr{iso}^\X_B(b)}:
              A \times B \ovs{\iso} \bfI\,(A \times B) \subset \Xbott.
\end{align*}
We name the \emph{inverse isomorphism}
  $\mr{jso}^\X_{\set{A\,|\,\chi}}: 
      \bfI\,\set{A\,|\,\chi} \ovs{\iso} \set{A\,|\,\chi}.$

\item
$\ulfamily\ \mr{iso}^\X_{\set{A\,|\,\chi}}: 
   \set{A\,|\,\chi} \ovs{\iso} \bfI\,\set{A\,|\,\chi} 
                                   \subset \Xbott \subset \N$ 
above, $\set{A\,|\,\chi}$ Object of $\PRa,$ 
is \emph{natural,} in the sense of the following commuting 
$\PRa$-\textsc{diagram} for a $\PRa$-map
  $f: \set{A\,|\,\chi} \to \set{B\,|\,\psi}:$

\begin{minipage} {\textwidth}
\xymatrix{
& \set{A\,|\,\chi}
  \ar[r]^f
  \ar[d]_{\mr{iso}^\X_{\set{A\,|\,\chi}}}^{\iso}
  \ar @{} [dr]|{=} 
  & \set{B\,|\,\psi}
    \ar[d]_{\iso}^{\mr{iso}^\X_{\set{B\,|\,\psi}}} 
\\
\set{\dot{A}\,|\,\dot{\chi}}
\ar @{=} [r]
& \bfI\,\set{A\,|\,\chi}
  \ar[d]^{\subset}
  \ar[r]^{\bfI\,f}
  & \bfI\,\set{B\,|\,\psi}
    \ar[r]^(0.4){\subset}
    & \bfI\,\set{B\,|\,\psi}\,\dot\cup\,\set{\bott}
      \ar[d]^{\subset}
\\
& \Xbott 
  \ar[rr]^{\dot{f} \bydefeq \bfI_{\PR}\,f}
  \ar[d]^{\subset}
  \ar @{} [drr]|{=} 
  & & \Xbott
      \ar[d]^{\subset}
\\
& \N 
  \ar[rr]^{\dot{f}}
  & & \N
}
\medskip
\quad $\PRa$ Embedding \textsc{diagram} for $\bfI\,f = \bfI_{\PRa}\,f$

\smallskip
$\bs{\in}
    \ \PRaX(I\,\set{A\,|\,\chi},\set{B\,|\,\psi})
        = \PRX(I\,\set{A\,|\,\chi},\set{B\,|\,\psi}).$
\end{minipage}

\medskip
In particular

\item
$$
\bfI\,f(a) \bydefeq
\begin{cases} 
\mr{iso}^\X_B \circ f \circ \mr{jso}^\X_A(a): 
            \dot{A} \ovs{\iso} A \ovs{f} B \ovs{\iso} \dot{B} \\
\quad \myif\ \dot\chi (a) \doteq \an{\true}_A,\ \ie 
                       \ \myif\ \chi(\mr{jso}^\X_A(a)), \\
\bott \in \dot{B}\,\cup\,\set{\bott} \subset \Xbott\ \emph{otherwise,} \\
\quad \ie \ \myif\ \neg\,\chi(\mr{jso}^\X_A(a)).
\end{cases}
$$

\end{enumerate}

\medskip
By PR \emph{internalisation} we get from the above the following

\medskip
\textbf{Internal Embedding Theorem:}
With \emph{Internalisitions} $\mrPR: \N \to 2$ of
$\PR \bs{\subset} \ul{\N},$\ \ $\mrPRa: \N \to 2$
of $\PRa \bs{\subset} \ul{\N},$\ \ 
$\mrPRaX \subset \mrPRX \subset \cds{\N,\N}_{\PR}: \N \to 2,$ 
and the corresponding internalised notions of equality 
  $$\ \checkeq^{\PR}_k,\ \checkeq^{\PRa}_k,
     \ \checkeq^{\PRaX} \subset \checkeq^{\PRX}: \N \to \N \times \N$$
we get $\PRa$ \emph{injections}
\begin{align*}
& I = I(u): \mrPR \xto{\iso} I\,[\,\mrPR\,] 
                \subset \mrPRX/\checkeq^{\PRX}\ = \\
& \qquad\qquad\qquad
    = \mrPRX/\checkeq^{\PR}
         \subset [\N,\N] \defeq \cds{\N,\N}_{\PR}/\checkeq^{\PR}, \\
& \quad
    \text{as well as an extension of this}\ I\ \text{into} \\
& I = I(u): \mrPRa \xto{\iso} \mrPRaX = I\,[\,\mrPRa\,]
        \subset \mrPRX/\checkeq^{\PRX} \subset [\N,\N] 
                             = \cds{\N,\N}_{\PR}/\checkeq^{\PR}.
\end{align*}
Both injections are \emph{internal (Cartesian PR) functors,}
isomorphic onto their (enumerated) images 
$\mrPR^\X = I\,[\,\mrPR\,]$ and 
$\mrPRaX = I\,[\,\mrPRaX\,] \subset \N$ respectively.

(\emph{Enumerated}) \emph{injectivity} of $I$ is meant injectivity as 
a $\PRa$ map, more precisely: as a map in Theory $\PRaQ = \PRa+(\Quot):$ 
Theory $\PRa$ definitionally (and conservatively) enriched 
with \emph{Quotients} by (enumerated) equivalence \emph{relations} 
(\cf \NAME{Reiter} 1980), such as in particular the different
internal notions $\checkeq_k: \N \to \N^2$ above. The ``mother'' of 
all these is here 
$\checkeq\ =\ \checkeq^{\PR}_k: \N \to \mrPR \times \mrPR \subset \N^2.$

The second \emph{injectivity} -- corresponding to theories $\PRa,$
$\PRaX,$ and $\PRX$ reads, in terms of $\PR$ and $\PRa$ alone:
\begin{align*}
& I(u)\,\checkeq^{\PR}_k\,I(v) 
      \implies u\,\checkeq^{\PRa}_{j(k)}\,v:
                      \N \times \cds{A,B}^2 \to 2, \\
& k \in \N\ \free,\ u,v \in \cds{A,B}^2\ \free,
    \ j = j(k): \N \to \N\ \text{available in $\PR,$} \\
& A,B\ \text{in}\ \PRa\ \text{(meta) \ul{free};}
\end{align*}
analogeous meaning for the former internal (parallel: \emph{objective})
injectivity properties\ \textbf{\qed}


\smallskip
$[\,$As mentioned above, \emph{Coding} $\mrPR = \mrPR/\checkeq^{\PR}$ 
of Theory $\PR = \PR/=^{\PR}$ restricts to coding
  $\mrPRX = \mrPRX/\checkeq 
     = \mrPRX/\checkeq^{\PR} \subset 
                    \cds{\N,\N}_{\PR}/\checkeq^{\PR}:$
coding of Object and map terms of $\PRX$ as well as internalising 
its inherited (enumerated) notion of equality.$]$  

\smallskip
We now have all formal ingredients for \textbf{stating}
\emph{Recursive Characterisation} of 
(wanted) -- double recursive -- \emph{evaluation algorithms}
\begin{align*}
& \eps^{\PR} = \eps^{\PR} (u,a): 
    \mrPR \times \Xbott \iso \mrPR^\X \times \Xbott \parto \Xbott, \\
& \quad \text{and its extension} \\
& \eps = \eps^{\PRaX} (u,a): \mrPRaX \times \Xbott \parto \Xbott.
\end{align*}
These evaluations are to become formally \emph{partial} $\PRa$-maps,
\ie maps of Theory $\hatPRa,$ see Ch.\ 1.

(Formal) \emph{partiality} will be here \emph{not} of PR decidable
nature, in contrast to that of \emph{defined partial} -- $\PRa$ -- 
maps, of form $f: \set{A\,|\,\chi} \to \set{B\,|\,\psi}$ 
discussed above.

\medskip
\textbf{Double Recursive Characterisation of Evaluation Algorithms}
  $$\eps^{\PR}: \mrPR \times \Xbott \parto \Xbott\quad\text{and}\quad
      \eps = \eps(u,a): \mrPRaX \times \Xbott \parto \Xbott$$
to \emph{evaluate} all \emph{map codes} in $\mrPR \bs{\iso} \mrPR^\X$ 
on all \emph{arguments of} -- free variable on -- 
Universal Object $\Xbott.$


\smallskip
The (wanted) \textbf{characterisation} is the following: 
\begin{itemize}
\item
Exceptional case of $x = \bott \in \Xbott$ -- \emph{undefined argument case:} \\
$\eps(u,\bott) \doteq \bott: \mrPRa \to \mrPRa \times \Xbott \to \Xbott:$ 
Once a value is \emph{defined undefined,} it remains so under evaluation
of any map code.

\item
case of basic map constants $\bas: A \to B,$ namely $\bas$ one of 
$0: \one \to \N,$ $s: \N \to \N,$ $\id_A: A \to A,$ 
$\Delta_A: A \to A \times A,$ $\Theta_{A,B}: A \times B \to B \times A,$ 
\ $\ell_{A,B}: A \times B \to A,$ and $r_{A,B}: A \times B \to B,$ first
$A,B$ fundamental Objects, in $\PR:$ 
\begin{align*}
& \eps^{\PR} (\code{\bas},a) = \bas(a): 
                       \Xbott \sqsup A \to B \sqsub \Xbott, \\
& \text{\ie (formally) in terms of theory}\ \PR^\X \bs{\iso} \PR: \\
& \eps^{\PR^\X} (\code{\bfI\,\bas},a) 
                     = \bfI\,\bas(a) = \bfI_{\PR}\,\bas(a): \\
& \Xbott \supset \dot{A} \to \dot{B} \subset \Xbott.
\end{align*}
Extension $\eps = \eps^{\PRaX}$ to the case of all -- \emph{basic} --
Objects of $\PRaX \bs{\supset} \PR^\X \bs{\iso} \PR:$ 
\begin{align*}
& \eps(\code{\bfI\,\bas},a) = \bfI\,\bas(a): 
               \Xbott \supset \bfI\,A \to \bfI\,B \subset \Xbott
                                          \ \text{(``again''),} \\
& \bydefeq
  \begin{cases}
  \mr{iso}^\X_B \circ \bas \circ \mr{jso}^\X_A(a): \\ 
  \quad
    \bfI\,A \ovs{\mr{jso}} A \ovs{\bas} B \ovs{\iso} \bfI\,B
                                        \ \myif\ a \in \bfI\,A, \\
  \bott\ \text{otherwise, \ie if}\ a \in \Xbott \sminus \bfI\,A
  \end{cases} \\
& : \Xbott \supset \bfI\,A \to \bfI\,B \subset \Xbott,
\end{align*}
this time $A$ and $B$ (suitable, basic) Objects, of $\PRa.$

\medskip
\begin{minipage} {\textwidth}
\textbf{Example:} 
\begin{align*}
& \eps(\code{\bfI\,\ell_{\set{\N\,|\,\mr{even}},\N \times \N}},x) \\
& =
  \begin{cases}
  \an{x_1} \in \an{\N} = \bfI\,\N\ \myif\ x 
    = \an{x_1;\an{x_{2\,1};x_{2\,2}}} \in \an{\N \times \N^2}\ \land\ 2|x_1, \\
  \bott\ \text{otherwise} 
\end{cases} \\
& : \Xbott \supset \an{\set{\N\,|\,\mr{even}} \times \N^2} 
         \to \an{\set{\N\,|\,\mr{even}}} \subset \an{\N} \subset \Xbott.
\end{align*}
\end{minipage}

\medskip
The \emph{compound} cases are the following ones:

\item \textbf{case} of evaluation of internally \emph{composed} 
\begin{align*}
& \an{v \odot u} \bydefeq \an{v \code{\circ} u}, 
                                         \ \text{for} \\
& u \in \cds{A,B}_{\PRaX},\ v \in \cds{B,C}_{\PRaX}
                       \ \text{``$\subset$''}\ \cds{\N,\N}_{\PR}:
\end{align*} 
\textbf{Characterisation} in this composition case is (is wanted): 
\begin{align*} 
& \eps(\an{v \odot u},a) = \eps(v,\eps(u,a)) 
          = \eps \parcirc (v,\eps \parcirc (u,a)): & (\odot) \\
& \cds{B,C} \times \cds{A,B} \times \Xbott \parto \Xbott,\ \text{in particular} \\
& \eps(\an{v \odot u},a) \doteq \bott \iff a \in \Xbott \sminus A,
                 \ \emph{defined undefined.}
\end{align*}

$[\,$Formally we cannot ``yet'' guarantee that $\eps$ be 
\emph{enumeratively terminating} at ``all'' \emph{regular} arguments,
``termination'' in a sense still to be \textbf{defined.}$]$

\smallskip
\textbf{Remark:} ``Definition'' in this -- central -- composition 
case is recursively \emph{legitimate,} by structural recursion on 
$\depth \an{v \odot u}$ \emph{down to} 
$\depth(u)$ and $\depth(v),\ u,v \in \mrPRaX,$ PR \textbf{definition}
of $\depth(u)$ for (general) 
  $$u = \an{\dot{\chi},\ring{u},\dot{\psi}} 
      \in \cds{\bfI\,\set{A\,|\,\chi},\bfI\,\set{B\,|\,\psi}}_{\PRaX}
        \subset \cds{\Xbott \setminus \set{\bott},\set{B\,|\,\psi}}_{\PRX}$$
see below.

\item cylindrified $\an{A \times v},$ $v \in \cds{B,B'}_{\PRaX}:$
\begin{align*}
& \eps(\an{A \times v},x) =
  \begin{cases}
  \an{x_1;\eps(v;x_2)} \in \an{A \times B'} \subset \Xbott & (\code{\times}) \\
  \quad 
    \myif\ x = \an{a;b} \in \an{A \times B} \subset \Xbott, \\
  \bott\ \text{otherwise}
  \end{cases} \\
& : \Xbott \supset \an{A \times B} \to \an{A \times B'} \subset \Xbott:
\end{align*}
\emph{evaluation in the cylindrified component.}

\item internally \emph{iterated} $u^{\code{\S}},$ for $u \in \cds{A,A}:$ 
\begin{align*}
& \eps(u^{\code{\S}}, \an{a;0}) = a,  & (\mr{iteration}\ \anchor) \\
& \eps(u^{\code{\S}}, \an{a;s\,n})  
  = \eps(u,\eps(u^{\code{\S}},\an{a;n})) \\
& = \eps \parcirc (u,\eps \parcirc (u^{\code{\S}},\an{a;n})):
                                      & (\mr{iteration}\ \step) \\
& (\mrPRaX \times \N) \times \Xbott \supset
    (\cds{A,A} \times \N) \times A \parto A \subset \Xbott,
\end{align*}
``$\supset$'' meaning ``again'': $\eps(u^{\code{\S}},x) \doteq \bott$ 
in all other cases. This case distinction is always here PR.

\item \emph{abstracted} map code $u$, of form 
\begin{align*}
& u = \an{\dot{\chi},\ring{u},\dot{\psi}} 
    \in \cds{\bfI\,\set{A\,|\,\chi},\bfI\,\set{B\,|\,\psi}}_{\PRaX}: \\
& \eps(u,a) =
  \begin{cases}
  \eps^{\PR} (\ring{u},a) \in \set{\dot{B}\,|\,\dot{\psi}} 
                                      = \bfI\,\set{B\,|\,\psi} \\
  \quad\myif\ \chi(a) \doteq \true \\
  \bott\ \emph{otherwise}
    \ \ie\ \myif\ a \in \Xbott \sminus \bfI\,\set{A\,|\,\chi}
  \end{cases} \\
& : \mrPRaX \times \Xbott \supset \cds{\set{\dot{A}\,|\,\dot{\chi}},
                                     \set{\dot{B}\,|\,\dot{\psi}}} 
       \parto \set{\dot{B}\,|\,\dot{\psi}} \subset \X
                   \subset \Xbott.
\end{align*}

\end{itemize}

\smallskip
\textbf{Remark:} If we restrict (wanted) evaluation $\eps$
to \emph{fundamental} map codes, out of 
  $$\mrPR\ [\ \sqsubset \mrPRaX\ ] \sqsubset \mrPRX 
                                       \subset \cds{\N,\N}_{\PR},$$
-- omit last case above and the ``\,$\bfI$\,'' in description of $\eps$ 
above throughout --  \\
we get, by $\PRa$ implications in cases above for 
\emph{basic map constants,} \emph{composition,} \emph{cylindrification,} 
as well as of \emph{iteration} 
characterisation of (wanted) \emph{fundamental} evaluation
\begin{align*}
& \eps^{\PR} = \eps^{\PR} (u,a): 
    \mrPR \times \Xbott \sqsupset \cds{A,B}_{\PR} \times A 
                                         \parto B \sqsubset \Xbott, \\
& \quad A,B \sqsubset \Xbott\ \text{fundamental,}\ \text{restriction of} \\
& \eps = \eps(u,a) = \eps^{\PRaX} (u,a): 
                  \mrPRaX \times \Xbott \parto \Xbott\ \text{above,}
\end{align*}
both to be characterised (within Theorie(s) $\piOR$ to come),
as formally \emph{partial} $\PRa$ maps -- out of Theory $\hatPRa$ --, 
but \emph{on-terminating} in $\piOR,$ and to be \textbf{defined} 
below as \emph{Complexity Controlled Iterations} ``$\mrCCIO$'s'' with
\emph{complexity values} in Ordinal $\N[\omega].$

Considering this restricted, \emph{fundamental} evaluation 
$\eps^{\PR}: \mrPR \times \Xbott \parto \Xbott$ will be 
helpfull, in particular since the Objects of $\PRa$ are
nothing else then \emph{fundamental} predicates
$\chi: A \to 2,$ still more formal: \emph{fundamental maps} 
$\chi: A \to \N$ such that
$\neg \circ \neg \circ \chi =^{\PR} \chi: A \to \N \to \N \to \N.$

\smallskip
\textbf{Recursive Legitimacy} for ``\textbf{definition}'' above of
evaluation $\eps$ is obvious for all cases above, except
for second subcase of case of \emph{iterated,} since in the other
cases recursive reference is made (only) to map terms of lesser 
$\depth.$ 

\smallskip
Here $\depth(u): \mrPRaX \to \N$ is \textbf{defined} PR as follows:
\begin{align*}
& \depth(\code{\id_A}) \defeq 0\ \text{for}\ A\ \text{fundamental,} \\
& \quad 
    \text{as well as for}\ A = \set{A'\,|\,\chi}
                               \ \text{basic, in}\ \PRa. \\
& \depth(\code{\bas'}) \defeq 1\ \text{for}\ \bas': A \to B \\ 
& \quad
    \text{one of the other basic map constants, in}\ \PRa;
                                               \ \text{further PR:} \\
& \depth(\an{v \odot u}) \defeq \depth(u)+\depth(v)+1: \\ 
& \cds{B,C}_{\PRaX} \times \cds{A,B}_{\PRaX} \to \N^2 \to \N.
\end{align*}
We then get automatically
\begin{align*}
& \depth_{\PRaX} \an{\code{\set{\dot{A}\,|\,\dot{\chi}}},u,
                              \code{\set{\dot{B}\,|\,\dot{\psi}}}} \\
& = \depth_{\PR^\X} \an{\code{\dot{A}},\code{\dot{B}}}
       = \depth_{\PR} (u): \cds{A,B}_{\PR} \subset \mrPR \to \N: \\
& \text{forget about (depth of) Domain and Codomain.}
\end{align*}

  
  
Using this $\depth = \depth(u): \mrPRaX \to \N,$ (wanted) characterisation 
above of $\eps^{\PR}$ and $\eps = \eps^{\PRaX}$ is recursively
\emph{legitimate} for all cases except -- a priori -- the iteration case,
since in those cases it recurs to its ``definition'' for map terms
with (strictly) lesser $\depth.$ 

In case of an iterated, reference is made to a term with 
\emph{equal} $\depth,$ but with decreased \emph{iteration counter:} 
from
  $$\iter(u^{\code{\S}}, \an{a;s\,n}) \defeq s\,n
      \quad\text{down to}\quad\iter(u^{\code{\S}}, \an{a;n}) \defeq n.$$
This shows \emph{double recursive,} (intuitive) \emph{legitimacy}
of our \textbf{``definition'',} more precisely: (double recursive) 
\textbf{description} of formally partial evaluation \\ 
$\eps: \mrPRaX \times \Xbott \parto \Xbott.$ A possible such 
(formally partial) map is \emph{characterised} by the above
\emph{general recursive} equation system. This system
constitutes a \emph{definition} by a (nested) \emph{double recursion} 
\`a la \NAME{Ackermann,} and hence in particular it constitutes
a \textbf{definition} in classical recursion theory.

\medskip
We now attempt to \textbf{resolve} basic evaluation $\eps,$ to be
\textbf{characterised} by \\
the above \emph{double recursion,} into a \textbf{definition} as 
an \emph{iteration} of a suitable \\ 
evaluation \emph{step} 
  $$e = e(u,x): \mrPRaX \times \Xbott \to \mrPRaX \times \Xbott,$$
first of a step 
$e = e^{\PR} (u,x): \mrPR \times \Xbott \to \mrPR \times \Xbott.$
 
In fact resolution into a \emph{Complexity Controlled} Iteration, 
$\mrCCI,$ which is to give, upon reaching complexity $0,$ evaluation 
\emph{result} $\eps(u,x) \in \Xbott$ in its right component.

\smallskip
For discussion of \emph{termination} of this (content driven)
iteration, we consider

\smallskip
\textbf{Complexity Controlled Iterations} in general: Such a $\mrCCIO$
is given -- in Theory $\PRa$ by data a (``predecessor'')
\emph{step} $p: A \to A$ coming with a \emph{complexity} $c: A \to O,$ 
such that $\PRa \derives\ \DeSta\,[\,p\,|\,c\,]\,(a): A \to 2,$
\text{where}
\begin{align*}
& \DeSta\,[\,p\,|\,c\,]\,(a) \defeq 
    [\,c(a) > 0 \implies p\ c(a) < c(a)\,] \\ 
& \qquad \emph{(strict \ul{De}scent above complexity zero)} \\
& \qquad
    \land\ [\,c(a) \doteq 0 \implies p(a) \doteq_A a\,] \\ 
& \qquad \emph{(\ul{Sta}tionarity at complexity zero)}.
\end{align*}
$O$ is an \emph{Ordinal,} here a suitable extension 
$O \succeq \N[\omega]$ of the semiring of polynomials in one 
indeterminate, with lexicographical order. \emph{Suitable} 
in the sense that we are convinced that it does not allow for 
infinitely descending chains.

\smallskip
\textbf{Examples of such ``Ordinals'',} besides $\N[\omega]:$ 
\begin{itemize}
\item 
$[\,\N$ itself as well as $\N \times \N,\ \N^{\ul{m}}$ with
hierarchical order are Ordinals \emph{below} $\N[\omega],$ but
we will need for our complexity values Ordinals 
$O \succeq \N[\omega] \iso \N^+\,]:$ 

\item
$O = \N^+ \identic \N[\xi] \identic \N[\omega]:$ $\N^+$ is the
set of non-empty strings, ordered lexicographically, and
to be interpreted here as \emph{coefficient strings} of
(the semiring of) polynomials over $\N$ in one indeterminate.
The order choosen on $\N[\omega]$ is in fact the lexicographical 
one on its coefficient strings in $\N^+.$

\item
$O$ the semiring $O = \N[\xi_1,\ldots,\xi_{\ul{m}}]$ in 
$\ul{m}$ indeterminates, the \emph{later} indeterminates
having \emph{higher priority} with respect to $O$'s order.

\item
$O$ the semiring $\N[\,\vec{\xi}\,] = \union_m \N[\xi_1] \ldots [\xi_m]$ 
in several variables (in arbitrary finitely many ones).
Order ``extrapolated'' from foregoing example.

\item
$O$ the \emph{ultimate (?) (countable)} Ordinal $\E$ given 
by arbitrarily \emph{balanced bracketing} of strings of 
natural numbers:
 
All of the above examples can be given the form of such 
sets of balanced-bracketed strings, but not containing 
\emph{singletons of singletons,} of form $\an{\an{\ldots}}.$
 
Admitting these \emph{pairs of double,triple,$\ldots$} brackets
leads to interpretation of $\E$ as the semi-algebra of strings 
of polynomials in (finitely many) indeterminates out of (countable) 
\emph{families of families of $\ldots$ families} of (candidates for)
indeterminates: indeterminates out of \emph{later families} then get 
\emph{higher priority} with respect to the order of $\E.$

\end{itemize}

Abbreviating predicate $\DeSta\,[\,p\,|\,c\,]\,(a): A \to 2$ 
given, ``positive'' \textbf{axiom} schema
$(\piO),$ of all $\mrCCIO$'s to \emph{on-terminate} 
-- whose equivalent \emph{contraposition} is schema 
$(\tilde{\pi}_O)$ of \emph{non-infinit descent} of the 
$\mrCCIO$'s --, reads:
\inference{ (\piO) }
{ $c: A \to O,\ p: A \to A$ $\PRa$ maps \\
& $\PRa \derives\ \DeSta\,[\,p\,|\,c\,]\,(a): A \to 2$ 
                                          \text{(see above);}  \\
& furthermore: for $\chi: A \to 2$ \emph{``test''} predicate, in $\PRa:$ \\ 
& ``test on reaching $0_O$'' by chain $p^n(a):$ \\
& $\PRa \derives\ \TerC\,[p,c,\chi] = \TerC\,[p,c,\chi]\,(a,n): 
                                                    A \times \N \to 2,$ \\ 
& \qquad\qquad
  $\defeq [\,c\ p^n(a) \doteq 0 \implies \chi(a)\,]: A \times \N \to 2$ \\ 
& \qquad
    (\emph{\ul{Ter}mination \ul{C}omparison} condition), \\ 
& \qquad
    with \emph{quantifier decoration:} \\
& $\PRa \derives\ (\forall a)\,[\,(\exists n)\,
                                 c\ p^n(a) \doteq 0_O \implies \chi(a)\,]$
}
{ $\piOR \derives \chi: A \to 2,$ \quad
    \ie $\chi =^{\piOR} \true_A: A \to 2.$ }

It is important to note in context of \emph{evaluation} -- that 
``emerging'' Theory $\piOR$ has same \emph{language} as basic 
$\mrPR$ Theory $\PRa.$ It just adds equations \emph{forced} by 
the additional schema.
\emph{Axis case} is $O :\,= \N[\omega],$ $(\pi) \defeq (\pi_{\N[\omega]}),$ 
$\pi\bfR \defeq \PRa+(\pi).$ Theory $\pi_\N\bfR$ would be just 
Theory $\PRa.$ 

\medskip
\textbf{Characterisation Theorem} for $\mrCCIO$'s:
Let \emph{complexity} $c = c(a): A \to O$ and \emph{predecessor} 
$p = p(a): A \to A$ be given, as in the antecedent of $(\piO)$ above.
Then (formally partial) $\hatPRa$ map 
  $$f(a) = p^\S \parcirc (a,\mu\,[\,c\,|\,p\,] \parcirc\ a): 
                                   A \parto A \times \N \to A$$
is nothing else then the $\hatPRa$ map ($\while$ loop)
$f = \wh[\,c > 0_O\,|\,p\,]: A \parto A,$ and we ``name'' it
$\wh_O[\,c\,|\,p\,]: A \parto A.$

Written with free variable, and \emph{dynamically:} 
  $$\wh_O[\,c\,|\,p\,]\,(a) 
       \pareq \wh[\,c(a) > 0_O\,|\,a :\,= p(a)\,]: A \parto A.$$
By $\while$ loop \textbf{Characterisation} in RFC1, 
this complexity controlled iteration ($\mrCCIO$) is characterised by
$$
  \wh_O = \wh_O[\,c\,|\,p\,] \parcirc a = 
  \begin{cases}
  a\ \, \text{if}\ c(a) \doteq 0_O \\
  \wh \parcirc p(a) \ \, \text{if}\ c(a) > 0_O
  \end{cases}
  : A \parto A.
$$
The standard $\hatPRa$ form of this $\mrCCIO$ reads:
\begin{align*}
& \wh_O = \wh_O[\,c\,|\,p\,] 
    = \bfan{(d_{\wh_O},\widehat{\wh}_O): D_{\wh_O} \to A \times A}: 
                                          A \parto A,  \ \text{with} \\
& D_{\wh_O} = \set{(a,n)\,|\,p^n(a) \doteq 0_O} \\
& d_{\wh_O} = d_{\wh_O} (a,n) = \ell(a,n) = a: D_{\wh_O} \to A,\ \text{and} \\
& \widehat{\wh}_O (a,n) 
    = p^\S(a,\min\set{m \leq n\,|\,p^m(a) \doteq 0_O}) 
                                       = p^n(a): D_{\wh_O} \to A,
\end{align*}
the latter because of \emph{stationarity} of $p: A \to A$ at 
\emph{zero-complexity.}

\smallskip
\textbf{Comment:} In terms of these $\while$ loops, equivalently: 
\emph{formally partial} PR maps, schema $(\piOR)$ says map theoretically: 
\emph{Defined-arguments} enumeration of the $\mrCCIO$'s \emph{have} 
image \emph{predicates,} and these predicative images equal \emph{true,} 
on the common \emph{Domain,} $A,$ of the given step and complexity. 
By \textbf{definition,} this means that these enumerations are \emph{onto,}
become so by axiom; and by this, all $\mrCCI_O$'s \emph{on-terminate.} 
In our context -- use \emph{equality definability} -- this is equivalent 
with \emph{epi} property of the defined-arguments enumerations of 
the $\mrCCIO$'s -- but \emph{not} with these enumerations
to be \emph{retractions.}

\smallskip
\textbf{Dangerous bound:}\footnote{added 2 Nov 2008} 
For complexity $c: A \to O$ above, descending with ``each'' step 
$p: A \to A,$ we have 
\begin{align*}
& \widehat{\wh}_O\,[\,c\,|\,p\,] \parcirc (\id_A,\mu_O) \pareq \wh_O:
                                     A \parto D_{\wh_O} \to A,
                                              \ \text{where} \\ 
& \mu_O = \mu_O[\,c\,|\,p\,]\,(a) \defeq \mu\set{n\,|\,c\,p^n \doteq_O 0}:
                                                           A \parto \N.
\end{align*}
But this $\mu_O = \mu_O[\,c\,|\,p\,]: A \parto \N$ cannot in general 
be a ($\hatPRa$) \emph{section} to 
$d_{\wh_O[\,c\,|\,p\,]}: D_{\wh_O[\,c\,|\,p\,]} \to A,$ since otherwise -- by 
\textbf{Section Lemma} in Ch.\ 1 -- $\hatPRa$ map 
$\mu_O: A \parto D_{\wh_O[\,c\,|\,p\,]}$ would become a PR (!) 
\emph{section} to defined-arguments (PR) enumeration 
$d_{\wh_O[\,c\,|\,p]},$ and hence $\wh_O[\,c\,|\,p]: A \to A$ 
would become PR itself. But at least for evaluation $\eps,$ which 
\emph{is} of $\mrCCI_O$ form, this is excluded by \NAME{Ackermann}'s
result that diagonalisation of $\eps$ -- ``evaluate $n$-th (unary) map
at argument $n$'' -- grows faster than any PR map. 

\smallskip
$[\,$Here we use the \NAME{Church} type result of Ch.\ 1, that any 
$\mu$-recursive map has a representation as a \emph{partial} 
$\mrPRa$ map, \ie that it can be viewed as a map within Theory $\hatPRa,$ 
as well as \emph{Objectivity} of evaluation $\eps$ which will be 
\textbf{proved} below.$]$



\smallskip
With motivation above, we now \textbf{define} $\PRa$ maps 
  $$e = e^{\PR} (u,a): \mrPRaX \times \Xbott \to \mrPRaX \times \Xbott$$ 
\emph{evaluation step,} and $c = c_{\PRaX}: \mrPRaX \to \N[\omega]$ 
\emph{(evaluation) complexity,} to give \textbf{evaluation} in fact
as a formally \emph{partial map} 
  $$\eps = \eps^{\PRaX} (u,a): \mrPRaX \times \Xbott \parto \Xbott,
                                 \ \text{within theory}\ \hatPRa,$$ 
$e$ and $c$ maps within Theory $\PRa.$ 

Partial \emph{evaluation} map $\eps$ then will be \textbf{defined} 
by iteration of PR \emph{evaluation step} 
$e: \mrPRaX \times \Xbott \to \mrPRaX \times \Xbott,$ descending in 
\emph{complexity}
  $$c = c(u,x) = c_{\eps} (u,x) \defeq c_{\PRaX} (u): 
            \mrPRaX \times \Xbott \onto \mrPRaX \to \N[\omega].$$
The (endo) \emph{evaluation step}
  $$e = e(u,x) = (e_{\map} (u,x),e_{\argg} (u,x)): 
                    \mrPRaX \times \Xbott \to \mrPRaX \times \Xbott$$
is \textbf{defined} below as a $\PRa$ map. Here left component
\begin{align*}
& e_{\map} (u,x): \mrPRaX \times \Xbott \to \mrPRaX
                       \ \text{designates the by-one-step} \\
& \quad
    \text{\emph{evaluated, reduced} map code, and right component} \\
& e_{\argg} (u,x): \mrPRaX \times \Xbott \to \Xbott\ \text{is to designate} \\ 
& \quad
    \text{the by-one-step (``in part'') \emph{evaluated argument.}}
\end{align*}
So here is the \textbf{definition} of evaluation step 
$e = (e_{\map},e_{\argg}),$ endo map of \\
$\mrPRaX \times \Xbott,$ by $\PRa$ \textbf{case distinction,} 
\cf (wanted) \textbf{characterisation} of \\
$\eps$ above:
\begin{itemize}
\item 
case of \textbf{basic} maps, of form $\bas: A \to B$ in $\PRaX(A,B):$
\begin{align*}
& e(\code{\mr{d\dot{a}s}},a) \defeq (\code{\id_{\dot{B}}},\mr{d\dot{a}s}(a)): 
      \Xbott \supset \dot{A} 
          \ovs{\mr{d\dot{a}s}} \dot{B} \ovs{\subset} \Xbott, \\
& \dot{A} \bydefeq \bfI\,A,\ A = \set{A'\,|\,\chi}\ \text{in}\ \PRa,
   \ \text{analogeously for}\ \dot{B.}
\end{align*}
``finished''. 

\smallskip
\textbf{Recall:} $\bas: A \to B$ is one out of the basic 
\emph{map constants} 
  $$\id_A,\ 0: \one \to \N,\ s: \N \to \N,\ !_A,\ \Theta_{A,B},\ \Delta_A,
   \ \ell_{A,B},\ r_{A,B},$$
$A,B$ Objects of $\PRa,$ in particular: $A,B\ \PR$-Objects.

\item 
\textbf{composition} cases: ``for'' (free variable) 
                 $v \in \cds{A,B},\ \cds{A,B} = \cds{A,B}_{\PRaX}:$ 
\begin{align*}
& e(\an{v \odot \code{\id_A}},a) \defeq (v,a)  & (\odot\ \text{anchoring}) \\
& \quad
    \in \cds{A,B} \times A \subset \mrPRaX \times \X 
                               \subset \mrPRaX \times \Xbott. 
\end{align*} 
For $((u,v),a) \in \cds{B,C} 
                 \times (\cds{A,B} \sminus \set{\code{\id_A}}) \times A 
                                       \subset (\mrPRaX)^2 \times \Xbott:$    
\begin{align*}
& e(\an{v \odot u},a) 
      \defeq (\an{v \odot e_{\map} (u,x)}, e_{\argg} (u,x)) \\ 
& \qquad
    \in \cds{\Dom(e_{\map} (u,x)),C} \times \Xbott 
                                   \subset \PRaX \times \Xbott,
\end{align*}
where $\Dom(e_{\map} (u,x)),$ Object of $\mrPRaX,$ is ``known'' 
-- \textbf{defined} PR on $\depth,$ in particular -- ``anchoring'' --
for $e_{\map} (u,x) = \dotbas$ above, $\Dom$ of form $\dot{A}$ 
in $\PRaX$ ($A$ \text{in} $\PRa$) is known, ``\etc'' PR.

So \textbf{definition} of $e$ in this composition case in toto, is PR 
on $\depth(\an{v \odot u}),$ \emph{``down to''} 
$\depth \an{v \odot e_{\map} (u,x)}.$

\item \textbf{cylindrified} cases:
\begin{itemize}
  \item 
  ``trivial'', \emph{termination} (sub)case:
    $$e(\an{\code{\id_A} \code{\times} \code{\id_B}},\an{a;b})
                       \defeq (\code{\id_{(A \times B)}},\an{a;b})$$
  ``finished'', and 

  \item 
  genuine cylindrified case: 
  for $v \in \cds{B,B'} \sminus \set{\code{\id_B}}:$
  \begin{align*}  
  & e(\an{\code{\id_A} \code{\times} v},\an{a;b}) \\
  & \defeq (\an{\code{\id_A} \code{\times} e_{\map} (v,b)},
                                            \an{a;e_{\argg} (v,b)}):
  \end{align*}
  apply evaluation (step) to right component $v$ and its argument $b.$
\end{itemize}

\item 
\textbf{iteration} case 
\begin{align*}
& u^{\code{\S}} \in \cds{\an{A \times \N},A},
                           \ \an{a;n} \in \an{A \times \N}\ (\free): \\
& e(u^{\code{\S}},\an{a;n}) \defeq (u^{[n]},a),
                               \ \text{where, by PR definition} \\  
& u^{[0]} \defeq \code{\id_A} \in \mrPRaX,\ \text{and}\
            u^{[s\,n]} \defeq \an{u^{[n]} \odot u} \in \mrPRaX \\
& \text{is \emph{code expansion} ``at run time''.}
\end{align*}
$[\,$This latter case of \textbf{definition} by \emph{code expansion,} 
is not very \emph{``effective'',} but logically simple.$]$

\end{itemize}

\smallskip
\textbf{Definition} of \emph{evaluation complexity,} to descend
with each application of \emph{evaluation (endo) step,} first
of $\PR$ map codes $u \in \mrPR:$ 


$c(u) = c_{\PRaX(u)}: \mrPRaX \to \N[\omega],$ is \textbf{defined} as a 
$\PRa$-map as follows: 
\begin{align*} 
& c\,\code{\id_A} \defeq 0 \cdot \omega^0 
                   = \mr{min}_{\N[\omega]},\ A\ \PRaX-Object, \\
& c\,\code{\bas'} 
        \defeq 1 \cdot \omega^0: \one \to \N[\omega], \\ 
& \qquad 
    \text{for $\bas'$ one of the other basic map constants of}\ \PRaX; \\ 
& \qquad
    \text{for}\ (u,v) \in \cds{B,C} \times \cds{A,B}  
             = \cds{B,C}_{\PRaX} \times \cds{A,B}_{\PRaX}: \\ 
& c\,\an{v \odot u} \defeq c(u)+c(v) +1 \cdot \omega^0 \in \N[\omega] \\ 
& \qquad 
    (\text{internal composition}\ \odot\ ); \\ 
& c\,\an{A \times v} = c\,\an{\dot{A} \code{\times} v}  
      \defeq c(v) +1 \cdot \omega^0: \mrPRaX \to \N[\omega] \\
& \qquad 
    (\text{internal cylindrification;}) \\
& \qquad
    \text{for}\ u \in \cds{A,A}_{\PRaX}: \\
& c\,(u^{\code{\S}}) 
      \defeq \omega^1 \cdot (c(u) + 1)
            = (c(u)+1) \cdot \omega^1: \\
& \mrPRaX \supset \cds{A,A} \to \N[\omega]
              \quad (\text{internal iteration}), \\
\qquad 
& \text{where}\  
   \omega = \omega^1 \identic 0;1\,,
       \ \omega^2 \identic 0;0;1\,,\ \omega^3 \identic 0;0;0;1\,
             \ \text{\etc in}\ \N[\omega], \\
& \N[\omega] \identic \N^+ = \N^* \sminus \set{\bot} \identic \N_{>0},  
                                                 \ \text{Ch.\ 1.}
\end{align*}
\textbf{Motivation} for above \textbf{definition} -- in particular
for this latter iteration case -- will become clear
with the corresponding case in \textbf{proof} of 
\textbf{Descent Lemma} below for \emph{basic evaluation}
  $$\eps = \eps(u,v) \defeq \wh[\,c_{\eps}\,|\,e\,]: 
      \PRaX \times \Xbott \parto \PRaX \times \Xbott \overset{r} 
                                                       {\onto} \Xbott.$$
\textbf{Remark:} As pointed out already above, restriction of
a $\PR^\X$ map code $u \in \cds{\dot{A},\dot{B}}$ to 
$u' \in \cds{\set{\dot{A}\,|\,\dot\chi},\set{\dot{B}\,|\,\dot\psi}}$
has no effect to complexity: If $u$ restricts this way, then
  $$c(u') = c^{\PRaX} (u') 
      = c^{\PR^\X} (u) = c^{\PR} (u) = c^{\PRaX} (u).$$




\smallskip
\textbf{Example:} Complexity of \emph{addition,} with 
 $+\ \bydefeq s^{\S}: \N \times \N \to \N,$ identified with
$\dot{+}: \an{\bfI\,\N \times \bfI\,\N} \to \bfI\,\N$ within $\PRaX:$
\begin{align*}
& c\,\code{+} = c\,\code{s^\S} = c\,(\code{s}^{\code{\S}}) \\
& = \omega^1 \cdot (c\,\code{s}+1) = 2 \cdot \omega \in \N[\omega] 
                            \ [\ \identic\ 0;2\ \in \N^+\ ].
\end{align*}
Evaluation \emph{step} and \emph{complexity} above are the right ones
to give 

\medskip
\textbf{Descent Lemma} for formally \emph{partially defined} 
and ``nevertheless'' \emph{on-terminating} evaluation map
  $$\eps = \eps(u,a) \bydefeq \wh[\,c_{\eps}\,\,|\,e\,]: 
       \mrPRaX \times \Xbott \parto \PRa \times \Xbott \ovs{r} \Xbott,$$
\ie for step \\
  $e = e(u,a) = (e_{\map},e_{\argg}): 
                     \mrPRaX \times \Xbott \to \mrPRaX \times \Xbott,$
and complexity 
  $$c_{\eps} = c_{\eps} (u,a) \defeq c(u): \mrPRaX \to \N[\omega]$$
we have Descent \emph{above} $0 \in \N[\omega],$ and Stationarity
\emph{at} complexity $0:$
\begin{align*}
\PRa \derives\ 
& c_{\eps} (u,a) > 0 \implies c_{\eps}\,e(u,a) < c_{\eps}\,c(u,a): \\ 
& \mrPRaX \times \Xbott \to \N[\omega] \times \N[\omega] 
                   \xto{ < \times < } \to 2^2 \ovs{\implies} 2,\ \ie \\
\PRa \derives\ 
& c(u) > 0 \implies c\,e_{\map} (u,a) < c(u): 
                      \mrPRaX \times \Xbott \to 2, & (\Desc) \\
& \text{as well as} \\
\PRa \derives\ 
& c(u) \doteq 0\ [\ \iff\ u\ \text{of form}\ u = \id_A\ ] \\
& \implies c_{\eps}\,e(u,a) \doteq 0\ [\ \land\ e(u,a) \doteq (u,a)\ ],
                                                         & (\Sta) 
\end{align*}
this with respect to the canonical, ``lexicographic'', and -- intuitively -- 
\emph{finite-descent} order of the polynomial semiring $\N[\omega].$ 

\textbf{Proof:} The only non-trivial case $(v,b) \in \mrPRaX \times \Xbott$
for the descent condition $c\ e(v,b) < c(v,b)$ 
is the iteration case 
  $$(u^{\code{\S}},\an{a;n}) 
      \in \cds{\an{A \times \N},A} \times A \subset \mrPRaX \times \Xbott.$$
In this ``acute'' iteration case we have in fact by induction on $n,$ 
\begin{align*} 
& c\,(u^{[n]}) = n \cdot c(u) + (n \dotminus 1), 
                                \ \text{since -- recursion:} \\
& c\,(u^{n+1}) = c\,\an{u \odot u^{[n]}} = c\,(u^{[n]})+c(u)+1
                                           = (n+1) \cdot c(u) + n, \\
& \qquad
    \text{whence} \\
& c_{\eps}\,e(u^{\code{\S}},\an{a;n}) = c\,(u^{[n]}) 
                                 \ (\text{\textbf{definition} of}\ e) \\ 
& = n \cdot c(u) + (n \dotminus 1) < \omega \cdot (c(u)+1), \\
& \qquad
    \text{since}\ \omega > m,\ m \in \N.
\end{align*}
$[\,$``$+1$'' in $c(u^{\code{\S}}) \defeq \omega \cdot c(u)+1$ is to 
account for the (trivial) case $\code{\id}^{\code{\S}}.]$

\smallskip
Stationarity at complexity $0 \in \N[\omega]$ is obvious \,\textbf{\qed} 
 
\medskip
This \emph{Basic Descent Lemma} makes plausible 
\textbf{global termination} of the ($\mu$-recursive) version of 
evaluation $\eps = \eps(u,x): \mrPRaX \times \Xbott \to \Xbott,$
in a suitable framework, here: it \textbf{proves} that this
\emph{basic} (formally) \emph{partial} evaluation map out of $\hatPRa:$ 
  $$\eps = \eps(u,x): 
       \mrPRaX \times \Xbott \parto \mrPRaX \times \Xbott \onto \Xbott$$
\emph{on-terminates} within Theory $\piOR = \PRa+(\piOR),$ for
Ordinal $O \succeq \N[\omega].$ This means that evaluation $\eps$
has an \emph{onto, epi} \emph{defined arguments} enumeration
\begin{align*}
& d_{\eps} = d_{\eps} (n,(u,x)) \defeq (u,x): \\
& D_{\eps} = 
    \set{(m,(u,x))\,|\,c\,\ell\ e^n(u,a) \doteq 0} 
                                               \to \mrPRaX \times \Xbott
\end{align*}
within $\piR \defeq \pi_{\N[\omega]}\bfR,$ and a fortiori in $\piOR,$ 
Ordinal $O \succeq \N[\omega],$ such choice of $O$ taken always here. 

\smallskip
\textbf{Remark:} Even if intuitively \emph{terminating,} and derivably
\emph{on-}terminating, partial map $\eps$ does not give (by
\emph{isomorphic translation}), a \emph{self-evaluation} of Theory 
  $$\piR = \PRa+(\pi) = \piR+(\pi_{\N[\omega]}),$$  
``\textbf{Dangerous bound}'' in \textbf{Summary} above.
Nothing is said (above) on evaluation of Theory 
$\hatpiOR = \widehat{\piOR}.$

\smallskip
In present context, we need an ``explicit''

\smallskip
Free-Variable Termination \textbf{Condition,} in particular
for our \emph{basic} evaluation $\eps,$ and later for its extension, 
$\eps_d,$ into an evaluation for \emph{argumented deduction trees.}

For a $\while$ loop in general, of form
\begin{align*}
& \wh[\,\chi\,|\,f\,]\,(a): A \parto A
           \ (\text{read:}\ \while\ \chi(a)\ \ul{do}\ a :\,= f(a)), \\
& \textbf{define}\ \
    [\,m\ \deff\ \wh[\,\chi\,|\,f\,]\,(a)\,] 
      \defeq [\,\neg\,\chi\ f^m(a)\,]: \N \times A \to 2:
\end{align*}
$m$ ``defines'' argument $a$ for $\while$ loop $\wh[\,\chi\,|\,f\,],$
to \emph{terminate} on this \emph{defined argument} after at most 
$m$ \emph{steps.}
 
\smallskip
This gives in addition:
\begin{align*}
& [\,m\ \deff\ \wh[\,\chi\,|\,f\,]\,(a)\,] 
     \implies \wh(a) \doteq_A \widehat{\wh} (a,m): \N \times A \to 2; \\ 
& [\,\wh(a) \doteq_A \widehat{\wh} (a,m)\,] 
     \bydefeq f^\S(a,\min\set{n \leq m\,|\,\neg\,\chi\ f^n(a)}):
                                                  \N \times A \to 2.
\end{align*}
Things become more elegant for $\mrCCIO$'s, because of 
\emph{stationarity} of $\mrCCI$'s at complexity $0 = 0_O \in O:$
\begin{align*}
\PRa \derives\ 
& [\,m\ \deff\ \wh_O[\,c\,|\,p\,]\,(a)\,] 
    = [\,c\ p^m(a) \doteq 0_O\ \land\ \wh_O(a) \doteq_A p^m(a)\,]: \\ 
& A \times \N \to 2,\quad\text{in particular:} \\
\PRa \derives\ 
& [\,m\ \deff\ \eps(u,x)\,]  
    = [\,c\,\ell\, e^m(u,x) \doteq 0\ 
        \land\ \eps(u,x) \doteq r\,e^m(u,x)\,]: \\ 
& \N \times (\mrPRaX \times \Xbott) \to 2.
\end{align*}
We will use this \emph{given} termination counter 
``$m\ \deff\ \ldots$'' only as a \emph{(termination)} condition (!), 
in \emph{implications} of form $m\ \deff\ \wh_O(a) \implies \chi(a),$ 
\quad $\chi = \chi(a)$ a \emph{termination conditioned} predicate. 
And we will make assertions on formally \emph{partial} maps such as 
evaluation $\eps$ and \emph{argumented deduction-tree evaluation} 
$\eps_d$ below, mainly in this termination-conditioned, 
``total'' form. 

So the main stream of our story takes place in theory 
$\PRa:$ we go back usually to the $\PRa$-building 
blocks of formally partial maps occurring, in particular to those 
of \emph{basic} \emph{evaluation} $\eps$ as well as those of 
\emph{tree evaluation} $\eps_d$ to come.

\smallskip
\textbf{Iteration Domination} above, applied to the \emph{Double}
\emph{Recursive} equations for $\eps,$ makes out of these
the following 

\begin{center}
\textbf{Dominated Characterisation Theorem} for \textbf{evaluation} 

\smallskip
$\eps = \eps(u,a): \mrPRaX \times \Xbott \parto \Xbott,$ \\

\smallskip
and hence equally for its \emph{isomorphic translation} \\

\smallskip
$\eps = \eps(u,a): \mrPRa \times \X \parto \X:$
\end{center}    
\begin{align*}
\PRa \derives\  
& [\,\eps(\code{\dotbas},a) \doteq \dotbas(a)\ \text{resp}
      \ \eps(\code{\bas},a) \doteq \bas(a)\,]\ \land\ : \\ 
& [\,m\ \deff\ \eps(v \odot u,a)\,] \implies
     \eps(\an{v \odot u},a) \doteq \eps(v,\eps(u,a)) \\
& \land\ [\,m\ \deff\ \eps(v,b)\,] \implies 
     \eps(\an{\code{\id} \code{\times} v}, \an{a;b}) 
                     \doteq \an{a;\eps(v,b)} \\
& \land\ \eps(u^{\code{\S}},\an{a;0}) 
           \doteq e^1(u^{\code{\S}},\an{a;0}) \doteq a \\ 
& \land\ [\,m\ \deff\ \eps(u^{\code{\S}},\an{a;s\ n})\,] \implies: \\
& \qquad
    m\ \text{defines all}\ \eps\ \text{instances below, and}: \\
& \eps(u^{\code{\S}},\an{a;s\ n}) 
    \doteq \eps(u^{\code{\S}},\an{\eps(u,a);n})
      \doteq \eps(u,\eps(u^{\code{\S}},\an{a;n})): \\
& \quad \N \times (\mrPRaX)^2 \times \X^2 \times \N \to 2, \\
& m \in \N\ \free,\ u,v \in \mrPRaX \subset\ \N \free
    \ \text{\resp} u,v \in \mrPRa \subset \N\ \free, \\ 
& a,b \in \X \subset \N,\ n \in \N\ \free.
\end{align*}

\textbf{Proof} of this \textbf{Theorem} by Primitive Recursion
(Peano Induction) on \\ 
$m \in \N\ \free,$ via case distinction on codes $w,$  
  $$w \in \mrPRaX \subset \cds{\X,\X}_{\PRX} 
                      \subset \cds{\N,\N}_{\PR} \subset \N,$$ 
and arguments $z \in \X$  appearing in the  different cases of 
the asserted conjunction, as follows, case $w$ one of the basic
map constants being trivial: 

All of the following -- \textbf{induction step} -- is situated in $\PRa,$ 
read: \\
$\PRa \derives\ \etc:$

\begin{itemize}

\item case $(w,z) = (\an{v \odot u},a)$ 
of an (internally) \emph{composed,} subcase $u = \code{\id}:$ obvious. 

\smallskip
Non-trivial subcase $(w,z) = (\an{v \odot u},a),$ $u \neq \code{\id}:$ 
\begin{align*}
& m+1\ \deff\ \eps(w,a) :\,= \eps(\an{v \odot u},a) \implies: \\ 
& \eps(w,a)
    \bydefeq e^\S((\an{v \odot e_{\map} (u,x)},e_{\argg} (u,a)),m) \\
& \qquad\qquad \text{by iterative definition of $\eps$ in this case} \\ 
& \doteq \eps(v,\eps(e_{\map} (u,a),e_{\argg} (u,a))) \\
& \qquad\qquad 
    \text{by induction hypothesis,} 
                                      \ \text{namely:} \\
& \qquad\qquad 
    m\ \deff\ \mu[\,c\,|\,e\,]\,(\an{v \odot e_{\map} (u,a)},
                               e_{\argg} (u,a)),\ [\ \ie\ \mu \leq m\ ] \\ 
& \implies: \\
& m+1\ \deff\ \eps(v,\eps(e_{\map} (u,a),e_{\argg} (u,a)))
                                     \doteq \eps(v,\eps(u,a)):
\end{align*} 
Same way back, by the same induction hypothesis, on $m,$ map code $v$
unchanged, ``passive'', in both directions of reasoning.

\item case $(w,z) = (\an{\code{\id} \code{\times} v},\an{a;b})$ 
of an (internally) \emph{cylindrified:} Obvious by definition of $\eps$
on a cylindrified map code.



\item case $(w,z) = (u^{\code\S},\an{a;0})$ \\
$\in \cds{\an{A \times \N},A} \times \an{A \times \N} 
                                  \subset \mrPRaX \times \X$ \\
of a null-fold (internally) iterated: again obvious.

\item case $(w,z) = (u^{\code\S},\an{a;n+1})$ \\
$\in \cds{\an{A \times \N},A} \times \an{A \times \N} 
                       \subset \mrPRaX \times \X$ \\   
of a genuine (internally) iterated: for $a \in \dot{A},$ $n \in \N$ free:
    %
    %
\begin{align*}
& (w,z) \doteq (u^{\code\S},\an{a;n+1}) \implies: \\
& m+1\ \deff\ \eps(w,z) \implies \\
& \eps(w,z) \doteq 
             \eps(e_{\map} (u^{\code\S},\an{a;n+1}),
                    e_{\argg} (u^{\code\S},\an{a;n+1})) \\
& \doteq \eps(u^{[n+1]},a) \doteq \eps(\an{u^{[n]} \odot u},a) 
                                  \doteq \eps(u^{[n]},\eps(u,a)) \\
& \qquad 
    \text{the latter by induction hypothesis on}\ m, \\ 
& \qquad 
    \text{case of internal composed} \\  
& \doteq \eps(u^{\code{\S}},\an{\eps(u,a);n}): 
\end{align*} 
same way back -- using \emph{bottom up characterisation} 
of the \emph{iterated} -- with $\eps(u,a)$ in place of $a,$ 
and $n$ in place of $n+1.$

\end{itemize}
This shows the (remaining) predicative---truncated---\emph{iteration} 
equations ``anchor'' and ``step'',
for an (internally) iterated $u^{\code{\S}},$ and so \textbf{proves} 
fullfillment of the above \textbf{Double Recursive} system of 
\textbf{truncated equations} for $\eps: \mrPRaX \times \X \parto \X,$
as well ``then'' for \emph{isomorphic translation}
$\eps: \mrPRa \times \X \parto \X,$ in terms of its defining 
components, within basic theory 
$\PRa \bs{\sqsubset} \hatPRa$ ``itself''\ \,\textbf{\qed}

\medskip
\textbf{Characterisation Corollary:} Evaluations -- $\hatPRa$-maps -- 
\begin{align*}
& \eps = \eps(u,a): \mrPRaX \times \X 
    \supset \cds{\bfI\,A,\bfI\,B}_{\PRaX} \times \bfI\,A 
                                     \parto \bfI\,B \subset \X \\
& \qquad
    \text{as well as -- \emph{back-translation} -- } \\
& \eps = \eps(u,a): \mrPRa \times \X 
             \sqsup \cds{A,B}_{\PRa} \times A \parto B \sqsub \X,
\end{align*}
now (both) \textbf{defined} as \emph{Complexity Controlled iterations} 
-- $\mrCCI$'s -- with complexity values in Ordinal $O :\,= \N[\omega],$ 
\emph{on-}terminate in Theorie(s) $\piOR$ ($O \succeq \N[\omega]$), by 
\textbf{definition} of these theory strengthenings of $\PRa,\ \hatPRa,$ 
and satisfy there the \textbf{characteristic} Double-Recursive 
equations stated for $\eps$ at begin of \textbf{section.}

\medskip
\textbf{Evaluation Objectivity:} We ``rediscover'' here the 
logic \emph{join} between the \emph{Object Language} level and 
the external PR Metamathematical level, join by 
\ul{externalisation} via evaluation $\eps$ above. The corresponding, 
very plausible Theorem says that evaluation $\eps$ \emph{mirrors} 
``concrete'' \emph{codes,} $\code{f}$ of maps 
$f: A \to B$ of Theories $\PR$ (via $\PR^\X = \bfI\,[\PR]$), $\PRaX$ 
as well as $\PRa,$ the latter via $\PRaX \bs{\iso} \PRa,$ back into 
these maps themselves.

\smallskip
\textbf{Objectivity Theorem:} Evaluation $\eps$ is \emph{objective,} 
\ie: for each \emph{single,} (meta \ul{free})
$f: \Xbott \sqsupset A \to B \sqsubset \Xbott$ in Theory $\PRa$ itself,
we have, with ``isomorphic translation'' of evaluation from $\PRaX:$
\begin{align*}
& \PRa \derives\ \eps(\code{f},a) = f(a): 
    \X \sqsupset A \to B \sqsubset \X,\ \text{symbolically:} \\
& \PRa \derives\ \eps(\code{f},\,\_\,\,) = f: A \to B,
\end{align*}
a fortiori: $\piOR \derives\ \eps(\code{f},a) = f(a): 
                                   \X \sqsupset A \to B \sqsubset \X.$ 

\smallskip
\textbf{Remark:} For such $f$ \ul{fixed},  
  $$\eps(\code{f},a) = \eps \parcirc (\code{f},a): 
                         A \to \cds{A,B} \times A \parto B$$ 
is in fact a $\PRa$ map
$\eps(\code{f},\,\_\,\,) = \eps(\code{f},a): A \to B,$
although in the \textbf{Proof} of the \textbf{Theorem} intermediate
steps are formally $\hatPRa$ equations ``$\pareq$'': But 
$\PRa \bs{\sqsubset} \hatPRa$ is a diagonal monoidal $\mrPR$ 
\emph{Embedding.}

\smallskip
\textbf{Proof} of \emph{Evaluation Objectivity} by \textbf{first:}
External structural \ul{recursion} on the nesting depth 
$\ul{de}p\ul{th}\,[\,f\,]$  
(``\ul{bracket de}p\ul{th}'') of $\PRa$-map $f: A \to B$ in question,
seen as external \ul{code}: $f \bs{\in} \ul{\N},$ and second:
in case of an \ul{iterated}, $g^\S = g^\S(a,n): A \times \N \to A,$ 
by $\PRa$-\emph{recursion} on \emph{iteration count} $n \in \N.$ 
This uses (dominated) Double Recursive Characterisation of evaluation 
$\eps$\ \,\textbf{\qed}

\medskip
\textbf{Finally} here: as forshadowed above, \emph{evaluations} ``split'' 
into (externally) \ul{indexed} Objective \ul{evaluation families}
  $$[\ \eps_{A,B} = \eps_{A,B} (u,a): 
           \cds{A,B} \times A \parto B\ ]_{A,B\ \text{Objects}},$$
with all of the above characteristic properties ``split''.


\smallskip
\textbf{Central} for all what follows is 
\textbf{(Inner) Soundness Problem} for \emph{evaluation}
  $$\eps = \eps(u,a): \PRaX \times \Xbott \parto \Xbott,\ \text{namely:}$$
Is there a ``suitable'' \emph{Condition}
  $\Gamma = \Gamma(k,(u,v)): \N \times (\mrPRaX)^2 \to 2,$ 
under which Theory $\PRa$ exports internal equality $u\,\checkeq_k\,v$
into Objective, predicative equality
  $\eps(u,a) \doteq \eps(v,a)?$ Formally: such that
\begin{align*}
\PRa \derives\ 
& \Gamma(k,(u,v)) 
     \implies [\,u \checkeq_k\,v \implies \eps(u,a) \doteq \eps(v,a)\,]: \\
& \N \times (\mrPRaX)^2 \times \X \parto \X \times \X \ovs{\doteq} 2\,?
\end{align*}
Such (``suitably conditioned'') \emph{evaluation Soundness} is strongly 
expected, and \ul{derivable} \emph{without condition} in classical 
Recursion Theory (and \textbf{set} theory) -- the latter two in the 
r\^ole of frame theory $\PRa$ above:

\smallskip\noindent
The formal \textbf{problem} here lies in \emph{termination.}


\section{Deduction Trees and Their Top Down 
                                                           \emph{Argumentation}}

As a first step for ``solution'' of the 
\textbf{(Conditioned) Soundness Problem} for evaluation
$\eps: \mrPRaX \times \X \parto \X,$ 
we fix in present \textbf{section} \emph{internal,} ``formalised'' 
\emph{Proofs} $\Proof_\T$ of map Theorie(s) $\T :\,= \piOR$ as 
(internal) \emph{deduction trees} $\dtree_k$ with nodes labeled by 
\emph{map-code internal equations.} These deduction trees are 
ordered by tree nesting-depth, and -- second priority -- code length: 
$\dtree_k$ is the $k$\,th deduction tree in this order, it (internally) 
\emph{proves, deduces} $\piOR$-equation $u\,\checkeq_k\,v.$

For reaching our goal of \textbf{Termination-Conditioned Soundness}
for evaluation 
\begin{align*}
& \eps = \eps(u,x): \mrpiOR \times \X = \mrPRa \times \X
                        \iso \mrPRaX \times \X \parto \X, \text{with} \\
& \piOR \derives\ \Gamma(k,(u,v)) 
     \implies [\,u\,\checkeq^{\piOR}_k\,v 
         \implies \eps(u,a) \doteq \eps(v,a)\,],
\end{align*} 
below, $\Gamma$ ``the'' suitable Termination condition, we consider 
\emph{evaluation} of \emph{argumented deduction trees} $\dtree_k/a,$
top down ``argumented'' starting with \emph{given} argument, to wanted
equation $\eps(u,a) \doteq \eps(v,a).$


For fixing ideas, we \emph{redefine} -- with the above counting 
$\dtree_k$ of deduction trees -- internal \emph{proving} as
\begin{align*}
& \Pro_{\piOR} (k,\ u\,\checkeq\,v) 
    \defeq \Pro_{\piOR} (\dtree_k,\ u\,\checkeq\,v) \\
& \bydefeq [\,u\,\checkeq^{\piOR}_k\,v\,]: 
          \N \times \mrPRa^2 \iso \N \times (\mrPRaX)^2 \to 2.
\end{align*}

\smallskip
Each such deduction tree, deducing -- \emph{root} -- internal
equation $u\,\checkeq\,v$ can canonically be \emph{argumented}
\emph{top down} with suitable arguments for each of its
(node) equations, when given -- just \emph{one} -- argument to its
\emph{root} equation $u\,\checkeq\,v.$

\smallskip
\textbf{Example:} Internal version of equational ``simplification''
Theorem \\
$s\,a \dotminus s\,b = a \dotminus b,$  namely
  $\an{\code{s} \odot \code{\ell} \code{\dotminus} \code{s} \odot \code{r}}
          \,\checkeq_k\,\an{\code{\ell} \code{\dotminus} \code{r}},$
``still'' more formal -- we omit from now on Object subscripts 
                                      (for $\piXOR = \PRaX$-Objects):
  $$\code{\dotminus} \odot \an{\code{s} \odot \code{\ell};
                                       \code{s} \odot \code{r}}
     \,\checkeq_k\,\code{\dotminus} \odot \an{\code{\ell};\code{r}},$$
$k \in \N$ suitable.

\smallskip
Internal \emph{deduction tree} $\dtree_k$ in this case:

\bigskip
\begin{minipage} {\textwidth}
$\dtree_k\quad = \quad$ 

\bigskip
\cinferenceq{} 
{ $\an{\ncode{s\,\ell} \ncode{\dotminus} \ncode{s\,r}}
     \,\checkeq_k\,\an{\ncode{\ell} \ncode{\dotminus} r}$ 
}
{ \cinferenceq{}
  { \parbox{5.5 cm} { 
      $\an{\ncode{s\,\ell} \ncode{\dotminus} \ncode{s\,r}}$ \\
      \phantom{M} $\checkeq_i\,\an{\ncode{\pre\,s\,\ell} 
                                      \ncode{\dotminus} \ncode{r}}$ } 
  } 
  { \cinferenceq{} 
    { \parbox{5.5 cm} { 
        $\an{\ncode{s\,\ell} \ncode{\dotminus} \ncode{s\,r}}$ \\
                    \phantom{M}$\checkeq_{ii}\,\an{\ncode{s\,\ell} 
                                        \ncode{\dotminus} \ncode{s\,0}} 
                                       \ncode{\dotminus} \ncode{r}$ } 
    }
    { \cinferenceq{}
      { \parbox{5.5 cm} { 
          $\an{\ncode{s\,\ell} \ncode{\dotminus} \ncode{s\,s\,r}}$ \\ 
          \phantom{M} $\checkeq_{iii}\,\an{\ncode{\pre\,s\,\ell} 
                                       \ncode{\dotminus} \ncode{s\,r}}$ } 
      }
      { \parbox{5.5 cm} { 
          $\an{\ncode{\ell} \ncode{\dotminus} \ncode{s\,r}}$ \\ 
          \phantom{M} $\checkeq_{iiii}\,
            \an{\ncode{\pre\,\ell} \ncode{\dotminus} \ncode{r}}$ \\
          $(\mr{definition}\ \text{of}\ \dotminus\,)\,\bs{.}$ }
      }
    } 
  }
      \cinferenceq{}
      { \parbox{4 cm} { 
          $\an{\ncode{\pre\,s\,\ell} \ncode{\dotminus} \ncode{r}}$ \\ 
          \phantom{M} $\checkeq_j\,\an{\ncode{\ell} \ncode{\dotminus} \ncode{r}}$ } 
      }
      { \parbox{4 cm} { 
          $\ncode{\pre\,s\,\ell}\,\checkeq_{ij}\,\ncode{\ell}$ \\
          $(\mr{definition}\ \text{of}\ \pre)\,\bs{.}$ }
      }
}        
\end{minipage}

\bigskip
When argument -- here for example $\an{a;7} \in \an{\N^2} \subset \X:$
$a \in \N$ free, and $7 \bydefeq s\,s\,s\,s\,s\,s\,s\,0: \one \to \N$ 
a constant: \emph{defined} natural number, is given to this 
(deduction) \emph{root,} it spreads down ``canonically'' to this 
tree $\dtree_k$ to give \emph{argumented deduction tree}

\bigskip
\begin{minipage} {\textwidth}
$\dtree_k/\an{a;7} \quad = \quad$ 

\bigskip
\cinferenceq{} 
{ $\ncode{\dotminus} \odot \an{\ncode{s}/a;\ncode{s}/\,7}
                        \sim \ncode{\dotminus}/\an{a;7}$ }
{ \cinferenceq{}
  { \parbox{4.5 cm} { 
      $\ncode{\dotminus} \odot \an{\ncode{s}/a;\ncode{s}/\,7}$ \\
      \phantom{M} 
        $\sim \ncode{\dotminus} \odot \an{\ncode{\pre\,s}/a;7}$ }
  } 
  { \cinferenceq{} 
    { \parbox{6 cm} { 
        $\ncode{\dotminus} \odot \an{\ncode{s}/a;\ncode{s}/\,7}$ \\
         \phantom{M}
         $\sim \ncode{\dotminus} \odot 
                        \an{\ncode{s}/a \ncode{\dotminus} \ncode{s\,0};7}$ } 
    }
    { \cinferenceq{}
      { \parbox{6 cm} { 
          $\ncode{\dotminus} \odot 
             \an{\ncode{s}/a;\ncode{s\,s}/\,7}$ \\ 
             \phantom{M}
             $\sim \ncode{\dotminus} \odot
                     \an{\ncode{\pre\,s}/a;\ncode{s}/\,7}$ } 
      }
      { \parbox{5.5 cm} { 
          $\ncode{\dotminus} \odot \an{a;\ncode{s}/\,7}$ \\ 
          \phantom{M}   
            $\sim \ncode{\dotminus} \odot \an{\ncode{\pre}/a;7}$ \\
          $(\mr{definition}\ \text{of}\ \dotminus\,)\,\bs{.}$ } 
      }
    } 
  }
      \cinferenceq{}
      { \parbox{4 cm} { 
          $\ncode{\dotminus} \odot \an{\ncode{\pre\,s}/a;7}$ \\ 
          \phantom{M}
          $\sim \ncode{\dotminus}/\an{a;7}$ }
      }
      { \parbox{4 cm} { 
          $\ncode{\pre\,s}/a \sim a$ \\
          $(\mr{definition}\ \text{of}\ \pre)\,\bs{.}$ } 
      }
}        
\end{minipage}

\bigskip
When evaluated -- by \emph{deduction tree evaluation} $\eps_d$ --
on \emph{argument} $\an{a;7} \in \an{\N^2}$ above -- this 
deduction tree, say $\dtree_k,$ \emph{should} (and will) give the 
following \ul{inference tree} $\eps_d(\dtree_k/\an{a;7}$ 
in \ul{Ob}j\ul{ect Level Lan}g\ul{ua}g\ul{e}:

\bigskip
$\eps_d(\dtree_k/\an{a;7}\quad = \quad$

\bigskip
\begin{minipage} {\textwidth}
\cinferenceq{} 
{ $s\ a \dotminus s\ 7 = a \dotminus 7$ }
{ \cinferenceq{}
  { \parbox{5 cm} {$s\ a \dotminus s\ 7
                    = \pre(s\ a) \dotminus 7$}
  } 
  { \cinferenceq{ \,(\mr{U_3})\, }
    { \parbox{5 cm} {$s\ a \dotminus s\ 7
                      = (s\ a \dotminus s\ 0) \dotminus 7$} 
    }
    { \cinferenceq{}
      { \parbox{5 cm} {$s\ a \dotminus s\ s\ 7 
                        = \pre(s\ a \dotminus s\ 7)$} 
      }
      { \parbox{5 cm} {$a \dotminus s\ 7 = \pre(a \dotminus 7)$} 
      }
    } 
  }
  \cinferenceq{}
  { \parbox{4.5 cm} { $\pre(s\ a) \dotminus 7 = a \dotminus 7$ }
  }
  { \parbox{4.5 cm} {$\pre\ s\ a = a$} 
  }
}        
\end{minipage}

\bigskip
Deduction- and \ul{Inference} trees above contain some ``macros'', 
for example \NAME{Goodstein}'s uniqueness rule $\mr{(U_3),}$ which is a
\textbf{Theorem} of $\PR,\ \PRa,$ and hence of $\piOR.$
Without such macros, concrete \ul{inferences}/deductions would 
become very deep and long. But theoretically, we can describe
these trees and their evaluation rather effectively by (primitive)
Recursion on \textbf{axioms} and axiom \textbf{schemata} of our 
Theorie(s), $\piOR.$

\smallskip
\textbf{Deduction Trees for Theory $\piOR:$} We \textbf{introduce} now
the family $\dtree_k,\ k \in \N$ of $\piOR$'s 
(internal) -- ``fine grain'' -- \emph{deduction trees:}
``fine grain'' is to mean, that each (\NAME{Horn} type) \emph{implication}
in such a tree falls in one of the following cases:
\begin{itemize}

\item Node entry is an equation directly given by (internalised) 
\emph{axiom.}

\item A bar stands for an implication of -- at most -- two 
``down stairs'' (internal) \emph{premise-}equations implying -- 
``upwards'' -- a \emph{conclusion}-equation, \emph{directly} 
by a suitable (internal) instance of an \textbf{axiom} schema 
of the Theory considered, here Theorie(s) $\piOR.$
\end{itemize}

So we are lead to \textbf{define} the natural-numbers-indexed family   
$\dtree_k$ as follows:   
  $$\dtree_k = \dtree^{\piOR}_k: \N \to \Bintree_{\mrPRa} \subset \X$$
is PR \textbf{given} by
\begin{align*}
& \dtree_0 = t_0 = \an{\code{id} \checkeq_0 \code{id}}
  \bydefeq \an{\code{id};\code{id}} \in \Bintree_{\mrPRa}, \\
& \dtree_k = \an{\an{u_k;v_k};\an{\dtree_{i(k)};\dtree_{j(k)}}}: 
                                         \N \to \Bintree_{\mrPRa}^2,
\end{align*}
the latter written \textbf{symbolically}
\bigskip
\cinferenceq{ \dtree_k\ =\quad }
{ $u_k\,\checkeq_k\,v_k$ }
{ \cinferenceq{}
  { $u_i\,\checkeq_i\,v_i$ }
  { $t_{ii}$ \quad $t_{ji}$ }
      \cinferenceq{}
      { $u_j\,\checkeq_j\,v_k$ }
      { $t_{ij}$ \quad $t_{jj}$ }
}
\bigskip
with -- as always below -- left \resp right \emph{predecessors} abbreviated
$i :\,= i(k),\ j :\,= j(k): \N \to PR^2,$ and recursively:
$ii :\,= i(i) = i(i(k))$ \etc 

$\Bintree_{\mrPRa} \subset \X$ above denotes the (predicative) subset
of those (nested) lists of natural numbers which code binary
trees with nodes labeled by $\PRa$ code pairs, meant to code
internal $\PRa \bf{\iso} \PRaX$ equations.

\smallskip
\textbf{Argumented Deduction Trees as Similarity Trees:} 
Things become easier, in particular so \emph{evaluation}
of \emph{argumented,} \emph{instantiated} deduction trees,
if treated in the wider frame of \emph{Similarity trees} 
  $$\Stree \defeq \Bintree_{\an{\mrPR \times \Xb}^2} \subset \N.$$
By \textbf{definition,} $\Stree$ is the predicative set of (coded) 
\emph{binary trees} with nodes labeled by \emph{Similarity} pairs 
$u/x \sim v/y,$ of pairs of \emph{map-code/argument pairs,} 
called ``Similarity pairs'', since in the interesting, 
\emph{legitimate} cases, they are expected to be converted 
into \emph{equal} pairs, by (deduction-) tree evaluation $\eps_d.$ 

\smallskip
General form of $t \in \Stree:$

\bigskip
\cinferenceq{ t \quad = \quad }
{ $u/x \sim v/y $}
{ \cinferenceq{}
  { $u'/x' \sim v'/y'$ }
  { \quad $t' \hfill \tilde{t}'$\quad } 
      \cinferenceq{}
      { $u''/x'' \sim v''/y''$ }
      { \quad$t'' \hfill \tilde{t}''$\quad }
}

\bigskip

$t',\ldots,\tilde{t}'' \in \Stree$ have (strictly) lesser $\depth$ 
than $t.$ 

In the \emph{legitimate} cases these pairs are ``expected'' to 
become \emph{equal} under $\Stree$-evaluation $\eps_d$ below 
-- \emph{argumented deduction tree evaluation:} \emph{legitimate}
are just \emph{argumented deduction trees,} of form $\dtree_k/x.$ 


We will \textbf{define} $\Stree$-\emph{evaluation} 
$\eps_d: \Stree \parto \Stree$ iteratively as $\mrCCIO$ via 
a PR \emph{evaluation step} $e_d = e_d(t): \Stree \to \Stree$ 
and a \emph{complexity} $c_d = c_d(t): \Stree \to \N[\omega].$

\smallskip
$[\,$Ordinal $O$ is here always choosen to extend $\N[\omega].$
Notation $\eps_d,\ e_d,\ c_d$ is choosen because \emph{restriction} 
to argumented \emph{\ul{d}eduction} trees ``is meant''.$]$ 

\smallskip
This construction of $\eps_d$ will extend \emph{basic} evaluation  
$\eps: \mrPR \times \X \parto \mrPR \times \X \onto \X,$
by suitable extension of basic \emph{step} 
$e: \mrPR \times \X \to \mrPR \times \X,$
and basic descending \emph{complexity} 
$c_{\eps} (u,a) = c_{\mrPR} (u): \mrPR \times \X 
                       \onto \mrPR \to \N[\omega].$

\bigskip
We will see in next section that \textbf{definition} of tree
evaluation step $e_d = e_d(t)$ needs formal definition
of \emph{argumentation} of arbitrary (legitimate)
deduction trees, 
  $(\dtree_k,x) \mapsto \TreeArg(\dtree_k,x) = \dtree_k/x \in \Stree.$

This will be the first, formally long, task to accomplish.
For making things homogeneous, we identify pure, 
argument-free trees, node-labeled with map pairs $u \sim v,$
with \emph{dummy argumented} trees, in $\dumTree \subset \Stree,$
dummy arguments given to (left and right sides of) all
of its \emph{similarity pairs:}
 
$\an{u \sim v} \mapsto \an{u/\Box \sim v/\Box},$ in particular
$\dtree_k$ is identified with $\dtree/\Box \in \dumTree \subset \Stree$
obtained this way.

\bigskip
We now give \textbf{Tree-\emph{Argumentation}} -- by 
\textbf{case distinction} PR on \emph{nesting depth} of
(arbitrary) $t \in \dumTree,$ for suitable \emph{arguments}
to be \emph{spread down,} from \emph{root} of $t,$ arguments
out of $\X,$ in particular out f $\an{\X \times \N} \subset \X$
\etc

\medskip
\textbf{Cases} of Tree-Argumentation, by \textbf{equation} \resp
\NAME{Horn} clause \emph{meant to deduce} \emph{root}
(or \emph{branch}) equation $u \sim v$ from left and right 
antecedents, see figure above of $t$ with this (general) \emph{root,}

This type of display of up-to-two explicit (binary) levels, plus 
recursive mention of lower branches, will suffice all our needs: 
two levels are enough for dislay of \NAME{Horn} type implications, 
from (up to two) equations to one equation.

\bigskip \noindent
-- (unconditioned) \textbf{equational} case $Equ\Case \subset \Stree$ 
for $\TreeArg:$
\begin{align*}
& \an{u/\Box \sim v/\Box}/x \defeq \an{u/x \sim v/x} \\
& \bydefeq \an{\an{u;x};\an{v;x}}: (\mrPRaX^)2 \times \X \to \Stree:
\end{align*}
replace the ``waiting'' dummy arguments by two equal (!) ``real'' ones.

\smallskip
This case covers in particular reflexivity of equality, associativity 
of composition, bi-neutrality of identities, terminality of $!,$ 
Godements and Fourman's equations for the induced, as well as the 
\emph{equations} for iteration.

\smallskip
-- \textbf{symmetry of equality} case \emph{Sym\Case:} straight forward.

\smallskip
-- \textbf{transitivity-of-equality} case (basic \textbf{forking} case): 
for $t \in \dumTree$ of form

\bigskip
\cinferenceq{ t\ =\quad }
{ $u/\Box \sim w/\Box$ }
{ \cinferenceq{}
  { $u/\Box \sim v/\Box$ } 
  { $t'$ \quad $\tilde{t'}$ } 
      \cinferenceq{}
      { $v/\Box \sim w/\Box$ }
      { $t''$ \quad $\tilde{t''}$ }
}

\bigskip
(hence $t',\,\tilde{t'},\,t'',\,\tilde{t''}$ all in $\dumTree$),
we \textbf{define} recursively: 

\bigskip
\cinferenceq{ t/x \ \defeq\quad }    
{ $u/x \sim w/x$ }
{ \cinferenceq{}
{ $u/x \sim v/x$ } 
{ $t'/x$ \quad $\tilde{t'}/x$ }
    \cinferenceq{}
    { $v/x \sim w/x$ }
    { $t''/x$ \quad $\tilde{t''}/x$ }
}  

 \bigskip
-- \textbf{composition compatibility} case: $t \in \dumTree$
of form 

\bigskip
\cinferenceq{ t\ =\quad }
{ $v \odot u/\Box \sim v'\odot u'/\Box$ }
{ \cinferenceq{}
  { $v/\Box \sim v'/\Box$ }
  { $t'$ \hfill $\tilde{t'}$ }
      \cinferenceq{}
      { $u/\Box \sim u'/\Box$ }
      { $t''$ \hfill $\tilde{t''}$ }
}

\bigskip
with all \emph{branches} in $\dumTree$ (or empty). Here we \textbf{define} 

\bigskip
\cinferenceq{ t/x \  \defeq\quad }
{ $v \odot u/x \sim v'\odot u'/x$ }
{ \cinferenceq{}
  { $v/\Box \sim v'/\Box$ }
  { $t'$ \hfill $\tilde{t'}$ } 
      \cinferenceq{}
      { $u/x \sim u'/x$ }
      { $t''/x$ \hfill $\tilde{t''}/x$ }
}

$[\,$Actual argument is given to pair $u \sim u'$ of first factors, 
and -- recursively -- to its deduction tree.$]$


\bigskip
-- \textbf{compatibility-of-cylindrification} case: straight forward

\bigskip   
Remain the following two cases:

\bigskip
-- $\FR!\Case,$ of \textbf{Uniqueness} of \textbf{initialised iterated:} 

\bigskip
for $t\ =$

\bigskip
\cinferenceq{}
{ $w/\Box \sim \an{v^{\S} \odot \an{\code{id} \code{\times} u}}/\Box$ }    
{ \cinferenceq{}
  { $\an{w \odot \an{u;\code{0}}}/\Box \sim u/\Box$ }
  { $t'$ \hfill $\tilde{t'}$ }
      \cinferenceq{} 
      { $\an{w \odot \an{v \code{\times} \code{s}}}/\Box
                                      \sim \an{v \odot w}/\Box$ }
      { $t''$ \hfill $\tilde{t''}$ }
}

\bigskip
we \textbf{define}

\bigskip
\begin{minipage} {\textwidth}
$t/\an{x;n}\ \defeq$

\bigskip
\cinferenceq{}
{ $w/\an{x;n} 
       \sim v^{\S} \odot \an{\code{id} \code{\times} u}/\an{x;n}$
}
{ \cinferenceq{}
  { $w \odot \an{u;\code{0}}/x \sim u/x$ }
  { $t'/x$ \hfill $\tilde{t'}/x$ } 
        \cinferenceq{}    
        { $w \odot \an{v \code{\times} \code{s}}/\an{x;n}
                     \sim \an{v \odot w}/\an{x;n}$ } 
        { $t''/\an{x;n}$ \hfill $\tilde{t''}/\an{x;n}$ }
} 
\end{minipage}

\bigskip
``For \textbf{example}'', \emph{fixing iteration count} and taking
another variable name, $a,$ instead of $x,$ we get,
with $7 \bydefeq s^7 \circ 0 = s\,s\,s\,s\,s\,s\,s\,0: \one \to \N:$

\bigskip
$t/\an{a;7}\ \defeq$

\bigskip
\begin{minipage} {\textwidth}
\cinferenceq{}
{ $w/\an{a;7} 
       \sim v^{\S} \odot \an{\code{id} \code{\times} u}/\an{a;7}$
}
{ \cinferenceq{}
  { $w \odot \an{u;\code{0}}/a \sim u/a$ }
  { $t'/a$ \hfill $\tilde{t'}/a$ } 
        \cinferenceq{}    
        { $w \odot \an{v \code{\times} \code{s}}/\an{a;7}
                     \sim \an{v \odot w}/\an{a;7}$ } 
        { $t''/\an{a;7}$ \hfill $\tilde{t''}/\an{a;7}$ }
} 
\end{minipage}

\bigskip
-- \textbf{final, extra} case $\piO\Case,$ of 
\textbf{on-terminating (``finite'') descent,} 
\emph{extra} for \emph{axis} Theory $\piO\bfR$ -- corresponding to 
schema $(\piO)$ of \emph{on-termination} of \emph{descending chains} 
in \emph{Ordinal} $O \succeq \N[\omega].$ This case is hard -- and
logically not self-evident, because it is \emph{self-referential}
in a sense: 

The first thing to do is \emph{internalisation} of 
(\NAME{Horn}) clause $(\piOR).$ We begin with \emph{internalisation}
of \textbf{definitions} $\DeSta[\,c\,|\,p\,]\,(a): A \to 2,$
-- of \emph{\ul{De}scent} + \emph{\ul{Sta}tionarity} -- of
\emph{complexity} $c,$ with each application of (predecessor) step
$p,$ as well as \emph{\ul{Ter}mination} \ul{C}omparison
formula (predicate) into -- obvious --  

\textbf{Definitions} -- ``abbreviations'' -- defining 
$\PRa \bf{\iso} \PRaX$ maps 
$\desta = \desta(u,v): \mrPRa \times \cds{\X,O} \to \cds{\X,2}$ 
(internal descent + stationarity), and 
$\terc = \terc(u,v,w): \mrPRa \times \cds{\X,O} \times \cds{\X,2} 
\to \cds{\X,2}$ (internal \emph{termination comparison}), are
immediate, ``term by term.''

Free variable $w \in \cds{\X,2}$ stands for an internal \emph{comparison} 
predicate, and $\terc(u,v,w)$ says -- internally -- that reaching 
\ul{c}omplexity zero: \emph{\ul{ter}minating,} when iterating $u$ 
``sufficiently'' often, makes \emph{comparison} $w$ (internally) true:

All this when ``completely'' \emph{evaluated} on suitable 
\emph{argument} out of $\X.$

The internal conclusion (\emph{root}) equation for $w$ then is  
$w \,\checkeq\, \code{\true}.$

\bigskip
\textbf{Putting all this together} we \textbf{arrive} at the following 
\textbf{type} of dummy argumented tree $t$ in the actual $\piO\Case:$

\bigskip
\cinferenceq{ t\ =\quad }
{ $w/\Box \sim \code{\true}/\Box$ }    
{ \cinferenceq{}
  { $\desta(u,v)/\Box \sim \code{\true}/\Box$ }
  { $t'$ \hfill $\tilde{t}'$ }
      \cinferenceq{} 
      { $\terc(u,v,w)/\Box \sim \code{\true}/\Box$ }
      { $t''$ \hfill $\tilde{t}''$ }
}

\bigskip
with, as always above, \emph{branches} 
$t',\,\tilde{t}',\,t'',\,\tilde{t}'' \in \dumTree \subset \Stree$
all \emph{dummy argumented} Similarity trees.
  
In analogy to the cases above, we are led to \textbf{define} for $t$ 
of the actual form:

\bigskip
\begin{minipage} {\textwidth}
$t/x\ \defeq$

\bigskip
\cinferenceq{}
{ $w/x \sim \code{\true}/x$ }
{ \cinferenceq{}
  { $\desta(u,v)/x \sim \code{\true}/x$ }
  { $t'/x$ \hfill $\tilde{t}'/x$ } 
      \cinferenceq{}    
      { $\terc(u,v,w)/\an{x;n_+} \sim \code{\true}/\an{x;n_+}$ } 
      { $t''/\an{x;n_+}$ \hfill $\tilde{t}''/\an{x;n_+}$ } 
    
}
\end{minipage}

\smallskip
These are the \emph{regular cases.} Cases not covered up to here
are considered \emph{irregular,} and \emph{aborted} by 
deduction-tree evaluation step $e_d = e_d(t): \Stree \to \Stree$
to be \textbf{defined} below, into 
$\an{\id/\Box \sim \id/\Box} \in \dumTree \subset \Stree.$

\smallskip
\textbf{Dangerous Bound} in case $(\piO)$ above: If one wants 
to \emph{spread down} a given argument, down from the \emph{root} 
of a dummy argumented tree to (the nodes of) its \emph{branches,} 
one may think that it be necessary to give all arguments needed 
on the way top down already to the \emph{root equation.}

In our actual ``argumentation case'' above, we did \textbf{not} give 
right component of a pair 
$\an{x;n} \in \an{\X}^2$ to the \emph{root} equation, only its left 
component $x.$ Only right subtree gets 
``full'' argument -- of form $\an{x;n_+}$ -- substituted 
at \emph{actual argumentation step.}

\emph{Logically,} argument (part) $n_+ \in \N$ 
has the character of a \emph{bound} 
variable, \emph{hidden} to the equation on top, here \\
``$w/x \sim \code{\true}$'', and to all equations way up to the
``global'' \emph{root} of the deduction tree provided with 
\emph{arguments} so far. 

``Free'' variable $n_+$ is to mean here \emph{classically} a 
variable which is \emph{universally bound} within an implication, 
more specifically: a variable which is \emph{existentially bound} 
in the \emph{premise} of (present) implication, since this variable 
does not appear within the \emph{conclusion} of the implication.

In classical Free-Variables Calculus, we would have to make sure that
the \emph{fresh} Free Variable -- here ``over'' $\N$ -- given to the right
hand branch above, \ie to $\terc(u,v,w)$ and its deductive descendants, 
gets not the \emph{name} of any (free) variable already occurring as 
a component of ``$x$'' in the present context. This possible conflict 
would be resolved \emph{classically} by \ul{countin}g names of 
Free Variables -- here of \emph{type} $\N$ -- given during 
\emph{argumentation,} and by giving to such a variable to be 
introduced in \emph{fresh} -- as in present case -- an \emph{indexed} 
name with index not used so far: this motivates notation ``$n_+$''
for this ``fresh'' variable.

\smallskip
In our \emph{categorical} Free-Variables Calculus -- with Free Variables
\ul{inter}p\ul{reted} as (nested) \emph{projections,} we 
interprete this \emph{fresh} variable $n_+$ \emph{introduced}
in ``critical'' argumentation case above, as -- additional -- 
\emph{right projection}
  $$\an{n_+} :\,= \an{r_{\X,\N}}: 
         \X \supset \an{\X \times \N} \to \an{\N},$$
of extended Cartesian product $\an{\X \times \N},$ extending
argument domain $\X$ for \emph{root} 
$\an{w/\Box \sim \code{\true}/\Box}.$ This way, categorically, 
variable $\an{n_+}$ behaves in fact -- intuitively -- 
as a \emph{fresh} Free Variable in the actual context.



\section{Evaluation Step on Map-Code/Argument Trees} 

We attempt now to extend basic evaluation $\eps$ of map-code 
argument pairs which has been given above as iteration of step
  $$e = e(u,x) = (e_{map} (u,x),e_{arg} (u,x)): 
                     \mrPR \times \X \to \mrPR \times \X,$$
into a -- terminating (?) -- evaluation $\eps_d$ of 
\emph{Similarity trees} $t,$ of general form displayed earlier. 

This evaluation comes -- in the present framework -- 
as a ($\mrCCIO)$ iteration of a suitable (descent) \emph{step} 
  $$e_d = e_d(t): \Stree \to \Stree,$$
on the set $\Stree subset \N$ of \emph{Similarity trees.}
 
\smallskip
$[\,$ $\Stree$ will host -- see below -- in particular all the 
\emph{intermediate results} of (iteratively) applying 
\textbf{deduction-tree evaluation step} $e_d$ to trees of form 
$t = \dtree_k/x:$ pure \emph{decuction} trees, \emph{argumented} 
by (suitable) constants or variables, \emph{argumentation} see
foregoing section.$]$

\smallskip
\textbf{Definition} of \emph{argumented-deduction-tree evaluation step}
  $$e_d = e_d(t): \Stree \to \Stree$$
recursively (PR) on $\depth(t),$ \ie on the \emph{nesting depth}
of $t,$ as a (binary) tree. More precisely: by recursive case 
distinction on the form of the two upper layers of $t.$ 
 
\bigskip
* For $t$ \emph{near flat,} \ie of form

\bigskip
\begin{minipage} {\textwidth}
\cinferenceq{ t\ =\quad }
{ $u/x \sim v/y$ }
{ $\an{\code{id}/x' \sim \code{id}/y'}$
      \quad $\an{\code{id}/x' \sim \code{id}/y'}$
}
\end{minipage}

\bigskip
we \textbf{define}
  $e_d(t) \defeq \emph{root} (t) = \an{u/x \sim v/y} \in \Stree.$

\smallskip
$[\,$In real \emph{deduction-life} we expect here $x' \doteq y'.]$

\bigskip
``The'' \textbf{exception} is the following 
\textbf{argument shift simplification} case -- arising in 
\emph{deduction} context below from the (internalised) schema of 
composition \textbf{compatibility} with equality (\emph{between} maps): 

\bigskip
$\bullet$\ Exceptional tree $t \in \Stree$ is one of form

\bigskip
\begin{minipage} {\textwidth}
\cinferenceq{ t\ =\quad }
{ $v \odot \code{id}/x \sim v'\odot \code{id}/x $ }
{ \cinferenceq{}
  { $v/\Box \sim v'/\Box$ }
  { $t'$ \hfill $t''$ }
  \quad $\code{id}/x \sim \code{id}/x$ 
}
\end{minipage}

\bigskip
$t',t'' \in \dumTree,$ pure map code trees, \emph{dummy argumented} 
at each argument place. $t'$ and/or $t''$ may be empty.

\smallskip 
\textbf{Note} that in this -- at least at surface -- \emph{legitimate} 
case, left and right argument, $x,$ of \emph{root ``equation''} 
of $t$ is the \emph{same.} If not, $t$ would be considered 
\emph{illegitimate,} and aborted by $e_d$ into 
$t_0/\Box \defeq \an{\id/\Box \sim \id/\Box}.$

\smallskip
For $t$ of exceptional (but regular) form above, we now \textbf{define}
recursively:  

\bigskip
\cinferenceq{ e_d(t)\ \defeq\quad }
{ $\an{v/x \sim v'/x}$ }
{ $t'/x$ \quad $t''/x$ } 

\bigskip
This is \textbf{shift} and \textbf{simplification:} 
right branch with its pair of identities is obsolete, 
its (common) argument $x$ is shifted, \emph{formally substituted,} 
into $v$ and $v'$ as well as into the trees ``responsable for the 
proof'' of hitherto not (yet) argumented \emph{equation,} formally:
``Similarity'' $v/\Box \sim v'/\Box.$

\smallskip
\textbf{Comment:} Present \textbf{case} is the first and only 
``surface'' case, where \textbf{definition} for evaluation step 
$e_d$ on ``deduction trees'' coming nodewise with variables, needs 
\emph{substitution, instantiation} of a (general) variable 
-- here $x \in \X$ -- into a general (!) ``deduction tree''.
 
By that reason, we had to consider the whole bunch of (quasi) 
legitimate \textbf{cases} of ``deduction'' trees 
and their ``natural'' spread down \emph{argumentation} into 
Similarity trees: $\dtree_k/x \in \Stree.$
 
\bigskip
$\bs{*}$\ \emph{Standard Case} which applies ``en cours de route'' of
stepwise tree-evaluation $\eps_d,$ step $e_d,$ where step
$e_d: \Stree \to \Stree$ is to apply basic evaluation step
$e: \mrPR \times \X \to \mrPR \times \X$ to all map-code/argument
pairs labeling the nodes of tree $t \in \Stree$ in question:

This is the case when $t \in Etree$ is of form
\begin{minipage} {5 cm}
\cinferenceq{\qquad t\ =\quad }
{ $u/x \sim v/y$ }
{ $t'$ \hfill $t''$ }
\end{minipage}

\bigskip
\ul{and} \emph{not} \emph{exceptional.} Here we \textbf{define} 
-- PR on $\depth(t):$

\bigskip  
\cinferenceq{ \qquad e_d(t)\ \defeq\quad }
{ $e(u/x) \sim e(v/y)$ }
{ $e_d(t')$ \quad $e_d(t'')$ } 

\bigskip

\textbf{SubException:} For $t' \in \dumTree$ we \textbf{define} in
this \emph{standard superCase:}
\begin{center}

\bigskip
\cinferenceq{ e_d(t) \defeq\quad }
{ $e(u/x) \sim e(v/y)$ }
{ $t'$ \quad $e_d(t'')$ }
\end{center}

\bigskip
Dummy tree $t'$ waits for \emph{later argumentation,} to come from 
evaluated right branch; an empty tree $t'$ in this case remains
empty under $e_d.$


\bigskip
What we still need, to become (intuitively) sure on 
\textbf{termination} of iteration 
  $$e_d^m(t): \Stree \times \N \to \Stree,$$
\ie to become sure that this iteration (stationarily) results 
in a tree $t$ of form 
  $t = \an{\code{id}/\bar{x} \sim \code{id}/\bar{y}},$ 
this for $m$ ``big enough'', is a suitable tree \textbf{complexity} 
  $$c_d = c_d(t): \Stree \to O \N[\omega],$$ 
which \textbf{strictly descends} -- above complexity zero -- 
with each application of \emph{step} $e_d.$ 

This just in order to give within $\piOR = \PRa+(\piO)$, 
by its schema $(\pi) = (\pi_{\N[\omega]})$ ($O \succeq \N_{\omega}$),
\emph{on-terminating descent} of argumented (deduction) tree 
evaluation $\eps_d,$ which is \textbf{defined} -- analogeously
to basic evaluation $\eps$ -- as the formally \emph{partial} map
  $$\eps_d = \eps_d(t/x)
     \bydefeq e_d^\S(t/x,\mu\set{m\,|\,c_d\ e_d^m(t/x) \doteq 0}): 
                                                \Stree \parto \Stree.$$
\smallskip
\textbf{Definition} of \emph{(argumented-)deduction tree complexity}
  $$c_d = c_d(t): \Stree \to \N[\omega] \preceq O$$ 
as natural extension of \emph{basic map complexity} 
  $$c = c_{\eps} (u,x) = c_{\mrPR} (u): 
       \mrPR \times \X \onto \mrPR \to \N[\omega]$$
to argumented ``deduction'' trees, \textbf{definition} in words:

\smallskip
$c_d(t)$ is $t$'s number of \emph{inference bars} plus the 
\emph{sum} of all \emph{map code complexities} $c_{\mrPR} (u)$
for $u \in \mrPR$ appearing in $t$'s node labels
(including the dummy argumented ones). The \emph{sum} is the
sum of polynomials in $\N[\omega]$ -- just here we need the
polynomial structure of Ordinal $O :\,= \N[\omega].$

\smallskip
$[\,$Formally this \textbf{definition} is PR on depth of tree $t.$ 
As in case $c_{\eps}$ for \emph{basic} evaluation 
$\eps = \eps(u,x): \mrPR \times \X \to \mrPR \times \X,$
the \emph{arguments} of the trees do not enter in this complexity.$]$

\smallskip
An easy (recursive) calculation of the -- different structural 
cases for -- trees $t \in \Stree$ \textbf{proves} 

\bigskip 
\textbf{Deduction-Tree Evaluation Descent Lemma:} 
Extended PR evaluation step
$e_d = e_d (t): \Stree \to \Stree$ \textbf{strictly descends} 
with respect to (PR) extended map code complexity  
$c_d = c_d (t): \Stree \to \N[\omega]$ \emph{above} complexity zero, \ie
  $$c_d (t) > 0 \implies c_d\ e_d (t) < c_d (t): 
                      \Stree \to \N[\omega]^2 \to 2,$$ 
and is stationary at complexity zero:
  $$c_d (t) \doteq 0 \implies e_d (t) \doteq t: \Stree \to 2.$$

$[\,$We have choosen complexity $c_d$ just in a manner to make sure
this stepwise \emph{descent.}$]$

\smallskip
So \emph{intuitively} we expect -- and can \ul{derive} in 
\textbf{set theory} -- that \emph{argumented-deduction-tree} 
evaluation $\eps_d: \Stree \to \Stree$ for $\piOR,$ 
\textbf{defined} as \emph{Complexity Controlled Iteration} 
($\mrCCIO$) of step $e_d$ -- descending complexity 
$c_d: \Stree \to \N[\omega] \preceq O$ -- always \emph{terminates,} 
with a \emph{correct} result of form
$\an{\id/\bar{x} \sim \id/\bar{y}},$ with $\bar{x} \doteq \bar{y},$
the latter when applied to a given argumented deduction tree 
of form $t = \dtree_k/x.$

\smallskip
We will not \textbf{prove} this termination: Termination
will be only a \textbf{Condition} in \emph{Main Theorem} next
section.


\section{Termination-Conditioned Soundness}

\emph{Termination Condition} -- a $\PRa$-\emph{predicate} -- for 
$\mrCCIO$'s was introduced above, and reads for (basic, iterative) 
\emph{evaluation}
\begin{align*}
& \eps = \eps(u,x) = e^{\mu\set{n\,|\,c_{\mrPR}e^n \doteq 0}}: 
                                 \mrPR \times \X \parto \X: \\
& [\,m\ \deff\ \eps(u,x)\,] \defeq [\,c_{\eps}\ e^m(u,x) \doteq 0\,]: 
                                \N \times \mrPR \times \X \to 2, \\
& m \in \N,\ u \in \mrPR,\ x \in \X\ \text{all}\ \free.
\end{align*}
Analogously for \emph{Argumented Deduction Tree evaluation}
defined as $\mrCCI$ ``over'' step $e_d = e_d(t): \Stree \to \Stree,$
$t$ an ``argumented deduction tree'', frame $\Stree,$ complexity
$c_d: \Stree \to \N[\omega]$ measuring \emph{descent.}

Here \emph{domination, truncation, quantitative ``definedness''} of 
termination reads
  $$[\,m\ \deff\ \eps_d(t)\,] \bydefeq [\,c_d\ e_d^m(t) \doteq 0\,]: 
                                 \N \times \Stree \to 2,\ m,\ t\ \free.$$
By definition of $\eps$ and $\eps_d$ -- in particular by 
stationarity at complexity zero, we obtain with
this ``free'' \emph{truncation} ($m \in \N$ free):
\begin{align*}
& [\,m\ \deff\ \eps(u,x)\,] \implies [\,c_{\mrPR}\,e^m(u,x) \doteq 0\,]
           \ \land\ [\,\eps(u,x) \doteq r\ e^m(u,x)\,],\ \ul{and} \\
& [\,m\ \deff\ \eps_d(t)\,] 
      \implies [\,c_d\,e_d^m(t) \doteq 0\,] 
                          \ \land\ [\,\eps_d(t) \doteq e_d^m(t)\,].
\end{align*}


Using the above abbreviations, we state the 

\bigskip
\textbf{\emph{Main Theorem}}, on \textbf{Termination-Conditioned Soundness:}

For theories $\piOR = \PRa+(\piO),$ of Primitive Recursion 
with (predicate abstraction and) \emph{on-terminating descent} 
in Ordinal $O \succeq \N[\omega]$ extending $\N[\omega],$ we have
\begin{enumerate} [(i)]

\item \emph{Termination-Conditioned Inner Soundness:}
\begin{align*}
  \piOR \derives\
  & [\,u\,\checkeq_k\,v\,]\ \land\ [\,m\ \deff\ \eps_d(\dtree_k/a)\,] \\
  & \implies m\ \deff\ \eps(u,a),\ \eps(v,a)\ \land\ : \\
  & \phantom{ \implies }
             \eps(u,a) \doteq r\ e^m(u,a) 
                  \doteq r\ e^m(v,a) \doteq \eps(v,a),  & (\bullet) \\ 
  & u,v \in \mrPR,\ a \in \X,\ m \in \N\ \free.
\end{align*}        

\emph{In words, this Truncated Inner Soundness says:} 
Theory $\piOR$ \ul{derives:}
    
\smallskip
\emph{
\textbf{If} for an internal $\piOR$ \emph{equation} $u\,\checkeq_k\,v$ 
the (minimal) \emph{argumented deduction tree} $\dtree_k/a$ for 
$u\,\checkeq_k\,v,$ top down \emph{argumented} with $a \in \X$ admits 
complete \emph{argumented-tree evaluation} -- \ie \textbf{If}
tree-evaluation becomes \textbf{stationary} after a finite number $m$ of 
evaluation steps $e_d$ --\,,} 

\emph{
\textbf{Then} both sides of this \emph{internal (!) equation} are completely 
evaluated on $a,$ by (at most) $m$ steps $e$ 
of original, \emph{basic} evaluation $\eps,$ into \emph{equal values.}
}

Substituting in the above ``concrete'' codes into $u$ \resp $v,$ 
we get, by \emph{Objectivity} of evaluation $\eps:$

\item \emph{Termination-Conditioned Objective Soundness for Map Equality:}

For $\piOR$ maps (\ie $\PRa$ maps) 
        $f,g: \X \supseteq A \to B \subseteq \X:$
\begin{align*}
\piOR \derives\
& [\,\code{f} \,\checkeq_k\, \code{g}
                             \ \land\ m\ \deff\ \eps_d(\dtree_k/a)\,] \\  
& \implies f(a) \doteq_B r\ e^m(\code{f},a) 
                      \doteq_B r\ e^m(\code{g},a) \doteq_B g(a):
\end{align*}
\emph{\textbf{If} an internal deduction-tree for (internal) equality of 
$\code{f}$ and $\code{g}$ is available, and \textbf{If} on this
tree -- top down argumented with a given $a \in A$ -- 
tree-evaluation \textbf{terminates,} will say: iteration of evaluation 
step $e_d$ becomes \textbf{stationary} after a finite number $m$ of steps,
\textbf{Then} equality $f(a) \doteq_B g(a)$ of $f$ and $g$ at 
this argument is the consequence.} 


Specialising this to case 
  $f :\,= \chi: A \to 2,\ g :\,= \true_A: A \to 2,$ we eventually get 
   
\item 
\emph{Termination-Conditioned Objective Logical Soundness:}
  $$\piOR \derives\
      \Pro_{\piOR} (k,\code{\chi}) \ \land\ m\ \deff\ \eps_d(\dtree_k/a)
                                       \implies \chi(a): \N^2 \to 2:$$ 
\emph{\textbf{If} tree-evaluation of a deduction tree of a 
predicate $\chi: \X \to 2$ -- the tree top down \emph{argumented} 
with ``an'' $a \in \X$ -- \textbf{terminates} after a finite number 
$m$ of tree-evaluation steps, \textbf{Then} $\chi(a) \doteq \true$ 
is the consequence.} 
\end{enumerate}

$[\,$The latter statement  reminds at the 
\emph{Second Uniform Reflection Principle}
$\mr{RFN'} (\T)$ in \NAME{Smorynski} 1977.$]$

\textbf{Proof} of ``axis'' \emph{Termination-Conditioned}
\emph{\textbf{Inner} Soundness:}

Without reference to \emph{formally partial} maps
$\eps: \mrPR \times \X \parto \X$ \\
and $\eps_d: \Stree \parto \Stree$ -- alone in $\piOR$ terms
$e: \mrPR \times \X \to \mrPR \times \X,$ 
$c_{\mrPR}: \mrPR \to \N[\omega],$ as well as
$e_d: \Stree \to \Stree$ and $c_d: \Stree \to \N[\omega]$ --
this \textbf{Theorem} reads:
\begin{align*}
\piOR \derives\
& u\,\checkeq_k\,v\ \land\ c_d\ e_d^m(\dtree_k/a) \doteq 0 \\
& \implies c_{\mrPR}\ r\ e^m(u,a) \doteq 0 
                     \doteq c_{\mrPR}\ r\ e^m(v,a) \\ 
& \phantom{ \implies }
             \ \land\ r\ e^m(u,a) \doteq r\ e^m(v,a): 
                                \N^2 \times \mrPR^2 \to 2  
                                            & (\breve{\bullet}) 
\end{align*}

\smallskip
\textbf{Proof} of $(\breve{\bullet})$ is by (primitive) recursion on 
$\depth(\dtree_k)$ of $k$\,th (internal) \emph{deduction tree} 
$\piOR$-\emph{proving} its \emph{root} $u\,\checkeq_k\,v.$
\emph{Argumented tree} $\dtree_k/a$ then has same depth,
and strictly speaking, we argue PR on $\depth(\dtree_k/a),$
by \emph{recursive case distinction} on the form of 
$\dtree_k/a.$

\smallskip
\emph{Flat} \textbf{SuperCase} $\depth(dtree_k) = 0,$ \ie SuperCase 
of \emph{unconditioned,} axiomatic (internal) \emph{equations} 
$u\,\checkeq_k\,v:$

We demonstrate our Proof strategy on the first involved of these cases,
namely \emph{associativity} of (internal) \emph{composition:} 
  $$\emph{Ass}\Case 
    \defeq [\,\dtree_k \doteq \an{\an{w \odot v} \odot u}
                \,\checkeq_k\, \an{w \odot \an{v \odot u}}\,]:
                                         \N \times \mrPR^3 \to 2.$$  
Here we first evaluate left hand side of equation substituted, 
``instantiated'' with (Free-Variable) \emph{argument} $a \in A:$ 
\begin{align*}
\piOR \derives\ 
& \emph{Ass}\Case \implies : \\
& m\ \deff\ \eps_d(\dtree_k/a) \\
& \implies [\,m\ \deff\ \eps(\an{w \odot v} \odot u,a)\,] \\
& \implies [\,m\ \deff\ \eps(u,a)\,]
               \ \land\,[\,m\ \deff\ \eps(w \odot v,\eps(u,a))\,] \\
& \phantom{\implies}
             \ \land\,\eps(\an{w \odot v} \odot u,a) 
                        \doteq \eps(w \odot v,\eps(u,a)) \\
& [\implies \text{the above}\ ] \\
& \phantom{\implies}
           \ \land\,[\,m\ \deff\ \eps(v,\eps(u,a))\,]
           \ \land\,\eps(v \odot u,a) \doteq \eps(v,\eps(u,a)) \\
& \phantom{\implies}
           \ \land\,[\,m\ \deff\ \eps(w,\eps(v \odot u,a))\,] \\
& \phantom{\implies}
           \ \land\,\eps(w \odot v,\eps(u,a)) 
                  \doteq \eps(w,\eps(v \odot u,a))
\end{align*}
Same way -- evaluation on a composed works step $e$ by step $e$
successively, it does not care here on brackets 
$\an{\ldots}$ -- we get for the right hand side of the equation:
\begin{align*}
& \piOR \derives\ \emph{Ass}\Case 
                  \implies [\,m\ \deff\ \eps_d (\dtree_k/a)\ \implies: \\
& m\ \deff\ \eps(w \odot \an{v \odot u},a)
        \,\land\,\eps(w \odot \an{v \odot u},a)
                       \doteq \eps(w,\eps(v,\eps(u,a)))\,]. 
\end{align*}
Put together:
\begin{align*}
\piOR \derives\ 
& \an{\an{w \odot v} \odot u}
       \,\checkeq_k\, \an{w \odot \an{v \odot u}} \implies 
                [\,m\ \deff\ \eps_d (\dtree_k/a) \implies : \\
& [\,m\ \deff\ \eps(\an{w \odot v} \odot u,a)\,] 
             \,\land\,[\,m\ \deff\ \eps(w \odot \an{v \odot u},a)\,] \\
& \,\land\,\eps(\an{w \odot v} \odot u,a)
                      \doteq \eps(w,\eps(v,\eps(u,a)))     
                            \doteq \eps(w \odot \an{v \odot u},a).]
\end{align*}
This proves assertion $(\breve{\bullet})$ in this 
\emph{associativity-of-composition} case.

\smallskip
Analogeous \textbf{Proof} for the other \textbf{flat,} equational cases,
namely \emph{Reflexivity of Equality,} \emph{Left and Right Neutrality
of Identities,} \emph{Functor property of Cylindrification,}
\NAME{Godement} \emph{equations for induced into Cartesian (!) product,}
\NAME{Fourman}'s \emph{equation for uniqueness of the induced,}
and finally, the two \emph{equations (!) for the (internally) iterated.}

\smallskip
We give the \textbf{Proof} for the latter case explicitely, since
it is logically the most involved one for Theory $\PR,$ and 
``characteristic'' for treatment of (internal) 
\emph{potential infinity.}  
 
For commodity, we choose -- equivalent -- ``bottom up'' 
presentation of this iteration case, namely \emph{iteration step} 
equation $f^\S(a,s\ n) = f^\S(f(a),n)$ instead of earlier axiom
$f^\S(a,s\ n) = f\ f^\S(f(a),n),$ formally: 
  $$f^\S \circ (A \times s) = f^\S \circ (f \times \N): 
                                A \times \N \to A \times \N \to A.$$
The \textbf{anchor} case statement for the internal iterated
$u^{\code{\S}}$ is trivial: apply evaluation step $e$ once.

\smallskip
\emph{Bottom up iteration} \textbf{step,}
    \textbf{Case} of \emph{genuine iteration equation:} 
\begin{align*}
& \piOR \derives \emph{iteq}\Case(k,u) \\
& [\,\defeq [\,\dtree_k \doteq  
           \an{u^{\code{\S}} \odot \an{\code{\id} \code{\times} \code{s}}
             \,\checkeq_k
                \,u^{\code{\S}} \odot \an{u \code{\times} \code{\id}}}\,]\,] \\
& \implies: m\ \emph{defines}\ \emph{all instances of}
                       \ \eps\ \emph{below,}\ \emph{and:} \\
& \eps(u^{\code{\S}} \odot 
    \an{\code{id} \code{\times} \code{s}},\an{a;n}) & (1) \\
& \doteq \eps(u^{\code{\S}},\eps
           (\code{id} \code{\times} \code{s},\an{a;n})) \\
& \doteq \eps(u^{\code{\S}},\an{a;s\ n})
      \doteq \eps(u^{\code{\S}} \odot 
                  \an{u \code{\times} \code{\id}},\an{a;n}). & (2)
\end{align*}
This common (termination conditioned) \emph{evaluation result} 
for both sides -- (1) and (2) -- of \ $\checkeq_k\, \in \mrPR^2,$ 
is what we wanted to show in this general iteration equality case. 

\smallskip
$[\,$Freyd's uniqueness case, to be treated below, is not an equational 
case, it is a genuine \NAME{Horn} case.$]$

\bigskip
Let us turn to the -- remaining -- \emph{genuine} \NAME{Horn} 
cases for assertion $(\breve{\bullet}).$ 

\smallskip
\textbf{Comment:} All of our \emph{arguments} below are to be 
\emph{formally} just Free Variables -- ``undefined elements'' --
or map constants such as $0,s\,0: \one \to \N.$ But since the 
variables usually occur in \emph{premise} \textbf{and} 
\emph{conclusion} of the \NAME{Horn} clauses -- to be \ul{derived} --
of assertion $(\breve{\bullet}),$ they mean \emph{the same} throughout such a clause: 
In this sense their ``multiple'' occurences are \emph{bounded together,} 
with meaning: \emph{for all.} ``But'' if such a variable occurs 
-- within an \emph{implication} -- only in the \emph{premise,} it means 
intuitively an \emph{existence,} to \emph{imply} the \emph{conclusio,}
\cf discussion of \emph{tree-argumentation} in the $(\piO)$-case.

\smallskip
\textbf{Proof} of Termination-Conditioned Soundness for the
``deep'', genuine \NAME{Horn} \textbf{cases} of $\dtree_k,$
\NAME{Horn} type (at least) at \emph{deduction} of \emph{root:}

\smallskip
\emph{Symmmetry- and Transitivity-of-equality} cases are immdediate.

\smallskip
-- \textbf{Compatibility} Case of 
\textbf{composition} with \textbf{equality:} 

\bigskip
\cinferenceq{ \dtree_k/a \quad = \quad }
{ $\an{v \odot u}/a \sim \an{v'\odot u'}/a$ }
{ \cinferenceq{}
  { $v/\Box \sim v'/\Box$ }
  { $\dtree_{ii(k)}/\Box \quad \dtree_{ji(k)}/\Box$ }
        \quad $\dtree_j/a$
} 

\bigskip
with two \textbf{subcases:}

\smallskip
-- \emph{exceptional,} \textbf{shift} case $u = u' = \code{\id},$
  $\dtree_j = t_0 = \an{\code{\id} \sim \code{\id}}:$ \\ 
In this subcase, to be treated separately because of exceptional
definition of step $e_d$ in this case, namely -- recursively --
\begin{align*}
& e_d(\dtree_k/a) \bydefeq \dtree_i/a\ \emph{(shift to left branch),}
                                         \ \text{and hence ``then''} \\
& \piOR \derives\ 
    m\ \deff\ \eps_d (\dtree_k/a) \implies: \\
& \eps_d(\dtree_k/a) \doteq \eps_d(e_d(\dtree_k/a)) 
                                    \doteq \eps_d(\dtree_i/a) \\
& \qquad 
    \text{whence, by induction hypothesis}
                              \ (\breve{\bullet}_i)\ \text{also:} \\
& \land\ \eps(v,a) \doteq \eps(v',a),\ \text{and hence, trivially:} \\
& \land\ \eps(v \odot \code{id},a) \doteq \eps(v' \odot \code{id},a):
                 \quad \emph{Soundness} \qquad (\breve{\bullet}_k). 
\end{align*}

\textbf{Genuine Composition Compatibility} \textbf{Case:} \emph{not} 
both $u,u'$ code of identity: This case is similar to -- and 
combinatorially simpler than the above. It is easily \textbf{proved}
by \emph{recursion} on $\depth(\dtree_k):$ we have just to 
\emph{evaluate} -- \emph{truncated soundly} -- argumented 
tree $\dtree_j/a.$ This branch evaluation is given by hypothesis
because of $\depth(\dtree_j/a) < \depth(\dtree_k/a).$

\smallskip
-- \textbf{Case} of Freyd's (internal) \textbf{uniqueness} 
of the iterated, is \textbf{case} of tree $t = \dtree_k/\an{a;n}$ 
of form 

\bigskip
\begin{minipage} {\textwidth}
$t = \dtree_k/\an{a;n} \quad = \quad $

\bigskip
\cinferenceq{} 
{ $w/\an{a;n} \sim \an{v^{\code{\S}} \odot 
                        \an{u \code{\times} \code{\id}}/\an{a;n}}$ }
{ \cinferenceq{} 
  { $w \odot \an{\code{\id};\code{0}}/a \sim u/a$ }
  { $\dtree_{ii} \hfill \dtree_{ji}$ }
        \quad 
         \cinferenceq{}
         { $w \odot \an{\code{\id} \code{\times} \code{s}}/\an{a;n}
                                      \sim \an{v \odot w}/\an{a;n}$ }
         { $\dtree_{ij} \hfill \dtree_{jj}$ }
}
\end{minipage}

\bigskip
\textbf{Comment:} $w$ is here an internal \emph{comparison candidate} 
fullfilling the same internal PR equations as 
$\an{v^{\code{\S}} \odot \an{u \code{\times} \code{\id}}/\an{a;n}}.$
It should is -- \emph{Soundness} -- evaluated identically to the latter, 
under \emph{condition} that evaluation of the corresponding argumented 
deduction tree terminates after finitely many steps, say after 
$m$ steps $e_d.$

\smallskip
Soundness \textbf{assertion} $(\breve{\bullet}_k)$ for the 
present Freyd's \emph{uniqueness} \textbf{case} is \textbf{proved} 
PR on $\depth(\dtree_i),\,\depth(\dtree_j) < \depth(\dtree_k),$ 
by established ``double recursive'' equations -- this time for 
evaluation of the \emph{iterated} -- established above for our 
\emph{dominated, truncated} case. These equations give in fact:
\begin{align*}
\piOR \derives\ 
& fr!\Case \implies: m\ \text{\emph{defines} 
                              all the following $\eps$-terms, and} \\ 
& \eps(w,\an{a;0}) \doteq \eps(u,a) 
    \doteq \eps(v^{\code{\S}} 
             \odot \an{u \code{\times} \code{\id}},\an{a;0}),  
                                          \ \text{as well as} & (\bar{0}) \\ 
\piOR \derives\ 
& fr!\Case \implies: m\ \text{\emph{defines} 
                              all the following $\eps$-terms, and} \\ 
& \eps(w,\an{a;s\,n})
       = \eps(w \odot \an{\code{\id} \code{\times} \code{s}},\an{a;n}) 
                                        \doteq \eps(v \odot w,\an{a;n}) \\
& \doteq \eps(v,\eps(w,\an{a;n})) & (\bar{s}).
\end{align*}
But the same is true for 
  \quad $v^{\code{\S}} \odot \an{u \code{\times} \code{\id}}$ \quad
in place of $w,$ once more by (truncated) double recursive equations for 
$\eps,$ this time with respect to the \emph{initialised internal iterated.} 

\smallskip
$(\bar{0})$ and $(\bar{s})$ put together show, by \textbf{induction}
on \emph{iteration count} $n\in \N$ \\
-- all other free variables $k,m,u,v,w,a$ together form the 
\emph{passive parameter} for this induction --  
\emph{truncated Soundness} assertion $(\breve{\bullet})$ of the 
Theorem for this \emph{Freyd's uniqueness} case, namely:
\begin{align*}
& \piOR \derives\ \emph{fr!}\Case \implies: 
          m\ \text{\emph{defines} all the following $\eps$-terms, and} \\
& \eps(w,\an{a;n}) \doteq 
    \eps(v^{\code{\S}} \odot \an{u \code{\times} \code{id}},\an{a;n}). 
                                                    & (\breve{\bullet}_k)
\end{align*}

\textbf{Final Case,} not so ``direct'', is internal version of case 
$(\piO)$ of ``finite'' descent -- in Ordinal $O \succeq \N[\omega]$ --
of (``endo driven'') $\mrCCIO$'s: \emph{Complexity Controlled Iterations}
with \emph{complexity values} in $O.$ In a sense, treatment of this 
\textbf{axiom} has something of reflexive, since it \emph{constitutes} 
theory $\piOR = \PRa+(\piO),$ and since \emph{on-termination} of 
evaluations $\eps$ and -- ``derived'' -- $\eps_d$
is forced by ``just'' this axiom, for $O :\,= \N[\omega].$

\textbf{Proof} strategy for this case is ``construction''
of ``super'' \emph{predecessor} $p_{\pi} = p_{\piO},$
``super'' \emph{complexity} $c_{\pi},$ and \emph{test} predicate
$\chi_{\pi},$ such that $p_{\pi}$ descends as long as $c_{\pi} > 0,$
is stationary at $0$ and \textbf{proves} 
\emph{Termination Conditioned Soundness} in present case by
application of schema $(\piO)$ itself (!) to \emph{data} 
$p_{\pi},\ c_{\pi},\ \chi_{\pi}.$

\smallskip
For treatment of this final case, we rely on \emph{internalisation}
of \textbf{Abbreviations} $\DeSta\,[\,p,c\,]: A \to 2:$
\ul{De}scent + \ul{Sta}tionarity of $\mrCCIO$ (given for step
$p: A \to A$ and Complexity $c: A \to O$), as well as
$\TerC\,[\,p,c,\chi\,]: A \to 2:$ \ul{Ter}mination \ul{C}omparison. 

The internal version of ``the above'' is -- with \\ 
$u \in \mrPR = \cds{\X,\X}_{\PRa}$ internalising \emph{iteration step} 
     $p: A \to A,$ \\
$v \in \cds{\X,O}$ internalising \emph{complexity} $c: A \to O,$ and \\
$w \in \cds{\X,2}$ internalising \emph{test} $\chi: A \to 2$ -- 
present argumented deduction tree

\bigskip
\begin{minipage} {\textwidth}
$\dtree_k/a \quad = \quad$

\bigskip
\cinferenceq{}
{ $w/a \sim \code{\true}$ }
{ \cinferenceq{}
  { $\desta(u,v)/a \sim \code{true}$ }
  { $\dtree_{ii}/a \quad \dtree_{ji}$ }
    \quad \cinferenceq{} 
           { $\terc(u,v,w)/\an{a;n_+} \sim \code{\true}$ }
           { $\dtree_{ij}/\an{a;n_+} \quad \dtree_{jj}/\an{a;n_+}$ } 
}
\end{minipage}

\bigskip 
\noindent
Here $\desta(u,v) \defeq$ \\ 
$[\,v \code{>} \code{0} \code{\impl} v \odot u \code{<} v\,] 
       \ \code{\land}\ [\,v \code{\doteq} \code{0} 
          \code{\impl} u \code{\doteq} \code{\id}\,]$ \\ 
internalises $\DeSta\,[p,c\,];$ 
    internalisation of $TerC\,[p,c,\chi\,]$ is \\
$\terc(u,v,w) \defeq 
    \an{v \odot u^{\code{\S}} \code{\doteq} \code{0}}  
                     \code{\impl} w \odot \code{\ell}.$

\smallskip
\textbf{Comment:} In the present $\piO\Case,$ (Free-Variable) 
\emph{argument} argument $n_+ \in \N$ for logical (right) 
predecessor-branch $\dtree_j$ within present instance $\dtree_k/a$ 
above, is not part of \emph{argument} argument ``given'' to 
(\emph{root} of) $\dtree_k.$
 
It is thought to be \emph{universally quantified} within ``its''
(argumented) right branch $\dtree_j/\an{a;n_+},$ so in fact it is 
thought to be \emph{existentially quantified} since it appears there
just in the \emph{premise,} \cf \textbf{discussion}
-- \textbf{Dangerous Bound} -- in foregoing section, on 
\emph{deduction-tree argumentation:} $n_+$ is
here a \emph{fresh} NNO variable, categorically seen as
``fresh'' name of a right projection. 

\smallskip
In what follows, we name this \emph{fresh} NNO-variable $n_+$ 
``back'' into $n.$ As you will see, there will result from
this no confusion, since we work just on two \emph{actual} 
levels of our argumented deduction tree $\dtree_k/a,$ 
only the right (argumented) branch comes with a ``visible'' 
``extra'' NNO variable, now called $n,$ giving
substitution, \emph{instantiation} $dtree_j/\an{a;n}.$  

\smallskip
We now attempt to show the assertion proper, $(\breve\bullet),$ 
for present $\piO\Case,$ via the original, \emph{objective,} 
schema $(\piO)$ \emph{itself.} We use for this the following 
``super'' \textbf{instance} of this schema: 

\bigskip
-- First we choose the (common) \emph{complexity/step} \textbf{Domain} 
  $A_{\pi} \subset \N \times \mrPR^3 \times A$ -- short for
``$A_{\piO}$'' -- predicatively \textbf{defined} as 
\begin{align*}
& A_{\pi} = A_{\pi} (a_{\pi}) = A_{\pi} (m,(u,v,w),a) \\ 
& \defeq [\,m\ \deff\ \eps(u,a),\ \eps(v,a), 
           \ \eps(v \odot u,a), \eps(w,a)\,] \\
& \N \times \mrPR^3 \times A \supseteq 
  \N \times (\cds{A,O} \times \cds{A,A} \times \cds{A,2}) \times A) 
                                                              \to 2, \\
& \text{and \emph{composit} Free Variable} \\
& a_{\pi} \defeq (m,(u,v,w),a)\ [\ = \id_{A_{\pi}}\ ]: 
                                         A_{\pi} \to A_{\pi}: 
\end{align*}
All of $a_{\pi}$'s \emph{components} free -- (nested) 
\emph{projections} -- in particular so  \emph{``dominating'',} 
formally: \emph{truncating,} $m \in \N,$ as well as 
$u \in \cds{A,A},\ v \in \cds{A,O},\ w \in \cds{A,2},$ and $a \in A.$

\smallskip
$[\,$ $A \subseteq \X$ (as well as $O$) are considered as 
\ul{meta-variables}, ranging over the subobjects of $\X,$ 
``\ie'' over the Objects of $\PRa$ -- and the Ordinals (of $\PRa$) 
extending $\N[\omega]$ respectively.$]$

In present internal \emph{proof,} \emph{deduction tree,} we have, 
with respect to \emph{left predecessor} branch 
  $$\dtree_i = \dtree_{i(k)} \in \Stree,$$
of actual deduction tree $\dtree_k,$ in particular with regard 
to its \emph{root:}   
  $$\piO\Case(k,(u,v,w))/a \implies 
      \roott\,\dtree_i/a \doteq \an{\desta(u,v)/a \sim \code{true}/a}.$$

\bigskip
-- Next ingredient for present application of \textbf{descent} schema 
is \textbf{complexity} 
  $$c_{\pi} = c_{\pi} (a_{\pi}): A_{\pi} \to O:$$ 
Here we choose Objectivisation of \emph{internal} complexity $v$ 
by \textbf{dominated, truncated evaluation,} namely
  $$c_{\pi} = c_{\pi} (a_{\pi}) = c_{\pi} (m,(u,v,w),a)
     \defeq r\ e^m(v,a) = \eps(v,a): A_{\pi} \to O.$$
The latter equation -- termination with $m$ -- follows by 
\textbf{definition} of Domain $A_{\pi}$ of $c_{\pi}.$

\smallskip
$[\,$(Just) here we need Ordinal $O \succeq \N[\omega]$ to 
extend $\N[\omega]:$ In the present approach,
\emph{syntactical complexity} of PR map codes 
takes values in $\N[\omega].$ But it is not excluded a priori that
in another attempt \eg Ordinal $\N^2$ would do.$]$

\bigskip
-- As \textbf{predecessor step} $p_{\pi}$ for present application of 
\textbf{descent} schema $(\piO),$ again within Theory $\PRa,$ 
we choose 
  $p_{\pi} = p_{\pi} (a_{\pi}): A_{\pi} \to A_{\pi},$ 
\emph{dominated,} \emph{truncated} by Free Variable $m \in \N,$ as
\begin{align*}
& p_{\pi} (a_{\pi}) = p_{\pi} (m,(u,v,w),a) \\  
& \defeq (m,(u,v,w),r\ e^m(v,a)) = (m,(u,v,w),\eps(v,a)): 
                                            A_{\pi} \to A_{\pi}.
\end{align*} 
Here again, as for \emph{complexity} $c_{\pi}$ above, \textbf{definition} 
of Domain $A_{\pi}$ provides \emph{termination}
  $m\ \deff\ \eps(v,a) \doteq_A r\ e^m(v,a)$ 
of (iterative) evaluation $\eps.$  

\bigskip
-- In choice of \emph{comparison predicate} 
$\chi_{\pi} = \chi_{\pi} (a): A_{\pi} \to 2$
we are free: a \emph{suitable} choice -- suitable for the needs of
\textbf{proof} in the actual case -- leads, analogeously to the
other ``$(\piO)$-data'', to externalisation via \textbf{evaluation}
of an \emph{arbitrary} internal predicate (free variable)
$w \in \cds{A,2} \subset \mrPR,$ as follows -- same receipt:
  $$\chi_{\pi} (a_{\pi}) = \chi_{\pi} (m,(u,v,w),a) 
      \defeq r\ e^m(w,a) = \eps(w,a): A_{\pi} \to 2.$$
\emph{Termination} $m\ \deff\ \eps(w,a) \doteq r\ e^m(w,a)$ of 
$\eps(w,a): A_{\pi} \to 2$ is as for complexity $c_{\pi}$
and predecessor $p_{\pi}$ above.

\bigskip
For due application of this -- now completely defined -- 
\textbf{instance} of schema $(\piO)$ -- which constitutes Theory 
$\piOR$ --  we check the two \textbf{antecedents,} as follows:
\begin{align*}
& \piOR \derives\ \DeSta_{\pi} (a_{\pi}): A_{\pi} \to 2: 
                            \ \text{\emph{left antecedent,} \textbf{and}} \\
& \piOR \derives\ \TerC_{\pi} (a_{\pi},n): A_{\pi} \times \N \to 2
                            \ \emph{right antecedent:}
\end{align*} 
By \textbf{definition} -- with \emph{composit} Free Variable 
  $a_{\pi} = (m,(u,v,w),a) \in A_{\pi}\ \text{above,}$ \\
actual \textbf{Left antecedent} reads: 
\begin{align*} 
& \DeSta_{\pi} (a_{\pi})
  = [\,c_{\pi} (a_{\pi}) > 0
         \implies c_{\pi}\,p_{\pi} (a_{\pi}) < c_{\pi} (a_{\pi})\,]  \\
& \phantom{\DeSta_{\pi} (a_{\pi})}
    \qquad \land\ [\,c_{\pi} (a_{\pi}) \doteq 0_O 
              \implies p_{\pi} (a_{\pi}) \doteq_{A_{\pi}} a_{\pi}\,]: 
                                                        A_{\pi} \to 2, \\
& \text{explicitely:} \\
& \DeSta_{\pi} (m,(u,v,w),a) 
    = [\,m\ \emph{defines}\  
         \text{all of the following instances of}\ \eps\,]\ \text{and} \\
& [\,\eps(v,a) > 0 \implies \eps(v,\eps(u,a)) < \eps(v,a)\,]
     \ \land\ [\,\eps(v,a) \doteq 0 
                    \implies \eps(u,a) \doteq_A a\,]: \\
& A_{\pi} \to 2,
\end{align*}
the latter $m$-terminations again by choice of Domain $A_{\pi}.$ 

\bigskip
-- \textbf{Right Antecedent}
  $$\TerC_{\pi} (a_{\pi},n) = \TerC((m,(u,v,w),a),n):
                                              A_{\pi} \times \N \to 2$$
then is -- for present $(\piO)$-\textbf{proof} instance 
``necessarily'' -- \textbf{defined} as
\begin{align*}
\TerC_{\pi} (a_{\pi},n) 
& \defeq [\,c_{\pi}\ p_{\pi}^\S(a_{\pi},n) \doteq 0 
                               \implies \chi_{\pi} (a_{\pi})\,] \\
& = [\,c_{\pi}\ p_{\pi}^n (a_{\pi}) \doteq 0 
                                 \implies \chi(a_{\pi})\,]: A_{\pi} \to 2. 
\end{align*}  

$[\,$(Free) \emph{iteration count} $n \in \N$ -- formally: $n_+ \in \N,$
see above -- comes in (only) here. $n$ is to count the number of iterated 
``applications'' of $e$ -- formally: \emph{evaluation steps} -- applied to 
\emph{internal endo} $u,$ on a given \emph{argument} $a \in A,$ 
for \emph{Comparison} with (evaluation of) internal \emph{test predicate} 
$w,$ again evaluated on $a.]$

\smallskip
We spell out \textbf{premise} equation 
$c_{\pi}\ p_{\pi}^n(a_{\pi}) \doteq 0:$
\begin{align*}
& [\,c_{\pi}\ p_{\pi}^n (a_{\pi}) \doteq 0\,] 
    \quad\ [\ = [\,c_{\pi}\ p_{\pi}^n (m,(u,v,w),a) \doteq 0\,]\ ] \\
& = [\,m\ \deff\ \eps(v,\bar{a}) \doteq 0\,] 
      \quad \text{with}\ \bar{a} = r\ e^n(u,a): A_{\pi} \to A_{\pi} \to A; \\
& \qquad
    \text{with auxiliary, \emph{dependent} variable}
                      \ \bar{a}\ \text{\emph{eliminated:}} \\
& = [\,m\ \deff\ \eps(v \odot u^{[n]},a) 
             \doteq \eps(v,\eps(u^{[n]},a)) \doteq 0\,].
\end{align*}
$[\,$ $u^{[n]} = u \odot \ldots \odot u$ is -- PR \emph{defined} -- 
$n$-fold \emph{code expansion,} see intermediate map-argument in 
iterative (basic) evaluation $\eps$ above.$]$

\smallskip
The above \textbf{defines} -- formally PR -- \textbf{premise equation}
$c_{\pi}\ p_{\pi}^n(a_{\pi}) \doteq 0.$

\smallskip
\textbf{Test predicate} $\chi_{\pi}: A_{\pi} \to 2$ in 
\text{right antecedent} $TerC(a_{\pi}): A_{\pi} \to 2$ 
is -- by \emph{choice} above -- 
  $$\chi_{\pi} (a_{\pi}) = \chi_{\pi} (m,(u,v,w),a)
      \bydefeq [\,m\ \deff\ \eps(w,a) \doteq r\ e^m(w,a)\,]: 
                                                A_{\pi} \to 2.$$
Putting things together into the actual \textbf{right antecedent} gives
\begin{align*}
\TerC(a_{\pi},n) 
& = [\,c_{\pi}\ p_{\pi}^n (a_{\pi}) \doteq 0 
                  \implies \chi_{\pi} (a_{\pi})\,] \\ 
& = [\,c_{\pi}\ p_{\pi}^n (m,(u,v,w),a) \doteq 0
                  \implies \chi_{\pi} (m,(u,v,w),a)\,] \\
& = [\,m\ \deff\ \eps(v,\eps(u^{[n]},a))\ \land\ m\ \deff\ \eps(w,a) \\
& \qquad 
    \land\ [\,\eps(v,\eps(u^{[n]},a)) \doteq 0 
                   \implies \eps(w,a)\,]\,]: A_{\pi} \times \N \to 2.
\end{align*} 
``Regular'' \emph{Termination} of all instances of 
$\eps: \mrPR \times \X \parto \X$ is here given
again by choice of $A_{\pi}: \N \times (\mrPR^3 \times A) \to 2.$ 

\smallskip
\textbf{Comment:} Free Variable $m \in \N$ -- ocurring in our 
\emph{premises} only -- means here intuitively assumption of 
\emph{``existence''} of a sufficiently large number -- $m$ -- 
such that $m$ iterations of evaluation step 
$e: \mrPR \times A \to \mrPR \times A$ suffice for \emph{regular} 
-- not \emph{genuinely truncated} -- $m$ fold iteration of step $e$ 
to give the wanted result $\eps(u,a) :\,= e^m(a).$

Intuitively such $m$  ``disappears'' -- better: 
is \emph{hidden} into the \emph{potentially infinite} -- in all 
of our (complexity controlled) iterations considered; and axiom 
schema $(\piO)$ which constitutes Theory $\piOR$ -- has just the 
sense to approximate -- without enriching the language 
(of Theory $\PRa$) -- this intuition of finite termination of 
$\PRa$ based, formally \emph{partial} evaluation. 

\smallskip 
So far the \emph{data.}

\bigskip
We now verify the needed \textbf{properties} of the two \emph{Antecedents} 
of schema $(\piO)$ for the actual instance  
  $$A_{\pi}\,,\ \DeSta_{\pi} (a_{\pi}): A_{\pi} \to 2,\ \text{and}\ 
         \TerC_{\pi} (a_{\pi},n): A_{\pi} \times \N \to 2:$$

\begin{itemize}
\item 
\textbf{Strict Descent} above complexity $0,$ and 
    \textbf{Stationarity} at $0:$ 
\begin{align*}
\piOR \derives\ 
& \piO\Case(k,(u,v,w))/a \implies: \\
& m\ \deff\ \eps_d(\dtree_i,a)\ \land\ \text{(``and gives further'')} \\
& m\ \deff\ \eps(\desta(u,v),a)\ \land\ 
       \doteq \eps(\code{\true},a) \doteq \true.
\end{align*}
This gives in particular 
  $\piOR \derives\ \DeSta_{\pi} (m,(u,v,w),a): A_{\pi} \to 2,$ 

the latter in particular by $\eps$-\emph{Objectivity} applied to
\textbf{definition} $(*)$ of $\desta(u,v)$ above, and by 
$m$-dominated (formally: $m$-truncated) \textbf{Double Recursive equations} 
for (iterative) evaluation $\eps: \mrPR \times \X \parto \X.$
  
\item \textbf{Termination Comparison} for \emph{comparison predicate}
$\chi_{\pi}: A_{\pi} \to 2:$ 
\begin{align*}
\piOR \derives\ 
& \piO\Case(k,(u,v,w))/\an{a;n} \implies: \\
& m\ \deff\ \eps_d(\dtree_j,\an{a;n})\ \land\ 
                               \text{(``gives further'')} \\
& m\ \deff\ \eps(\terc(u,v,w),\an{a;n}) \doteq \true, 
                                                \ \text{whence} \\
\piOR \derives\ 
& \TerC_{\pi} ((m,(u,v,w),a),n): A_{\pi} \to 2.
\end{align*}
The latter again by -- dominated, formally: truncated -- 
``characteristic'' (Double Recursive) equations for 
$\eps: \mrPR \times \X \parto \X.$
\end{itemize}

So we have verified \textbf{both Antecedents} for (objective) 
schema $(\piO),$ in its here needed \textbf{instance} 
$A_{\piO},\ \DeSta_{\piO},\ \TerC_{\piO}.$

\smallskip
\textbf{Postcedent} of this \emph{on-terminating descent} schema 
for theory $\piOR$ then gives
\begin{align*}
\piOR \derives\ 
& \chi_{\pi} (m,(u,v,w),a): A_{\pi} \to 2,\ \text{namely} \\
\piOR \derives\ 
& \piO\Case(k,(u,v,w))/a \implies \chi_{\pi}\,,
                              \ \text{and hence in particular} \\ 
\piOR \derives\
& \piO\Case(k,(u,v,w))/a \implies: \\
& m\ \deff\ \eps_d(\dtree_k/a) \implies 
        \eps(w,a) \doteq \true \doteq \eps(\code{\true_A},a): 
                                           & (\breve\bullet_k). 
\end{align*}
So in this \textbf{final case} too, (internal) \emph{root} equation 
  $$\roott\ \dtree_k \bydefeq \an{w\,\checkeq_k \code{true_A}}$$
is evaluated -- formally: \emph{termination-conditioned} evaluated --
into expected \textbf{objective} predicative equation:
  $$\piOR \derives\ [\,m\ \deff\ \eps_d(\dtree_k/a)\,] 
               \implies \eps(w,a) \doteq_A \eps(\code{\true_A},a).$$
This means that \emph{dominated, formally: truncated} evaluation
$\eps_d$ of \emph{argumented deduction trees} evaluates -- in case of
\emph{Termination} -- not only the \emph{map code}/\emph{argument} 
pairs in $\dtree_i/a = \dtree_{i(k)}/a$ as well as in 
$\dtree_j(k)/\an{a;n}$ into equal \emph{values,} but 
-- recursion -- by this also those of $\dtree_k/a,$ 
$a \in A \subseteq \X,$ all this in the present, last regular 
case of $(k,a) \in \N \times A \subseteq \N \times \X,$ 
and its associated \emph{deduction tree} $\dtree_k/a,$ 
$a$ (recursively) substituted, \emph{instantiated}  into 
\emph{pure, variable-free} internal (equational) \emph{deduction tree} 
$\dtree_k$ for any internal equation, general form $u\,\checkeq_k\,v.$ 

\smallskip
This -- exhaustive -- \emph{recursive case distinction} shows 
\emph{Dominated,} formally: \emph{truncated,} and more intuitive: 
\textbf{Termination-Conditioned,} \textbf{Soundness} 
for Theory $\piOR,$ relative to itself, and hence also
the other assertions of \textbf{Main Theorem,}
on \emph{Termination-Conditioned Soundness} \quad\textbf{\qed}

\smallskip
\textbf{Remark:} Universal set $\X \subset \N$ seems to give a 
good service: without it, we would have be forced (?) to define 
evaluation $\eps$ as a \ul{famil}y
  $$\eps = [\,\eps_{A,B}: 
                \cds{A,B} \times A \parto B\,]_{A,B \in \Obj_{\PRa}}$$
\ul{meta-indexed} over pairs of Objects of Theory $\PRa,$
as is usual in Category Theory for \emph{axiomatically} given evaluation
  $$\epsilon 
    = [\,\epsilon_{A,B}: B^A \times A \to B\,]_{A,B \in \Obj_{\mathbf{C}}},$$
$\mathbf{C}$ a (Cartesian) Closed Category in the sense of 
\NAME{Eilenberg}\,\&\,\NAME{Kelly} 1966 and 
\NAME{Lambek}\,\&\,\NAME{Scott} 1986.
(Observe our typographic distinction between the two ``evaluations'').

At least formally, a \emph{constructive} \textbf{definition} 
of evaluation as one single -- formally partial --
$\PRa$ map $\eps = \eps(u,x): \cds{\X,\X} \times \X \parto \X$ is 
``necessary'' or at least makes things simpler.    

So both, the typified approach -- traditional in Categorical main stream, as
well as the \NAME{Ehresmann} type one starting with just one \emph{class}
of maps -- and partially defined composition -- are usefull in our context:
\emph{Universal set} $\X$ -- of \emph{(codes of) strings} of natural
numbers here makes the join.

\smallskip
From this \emph{Main Theorem,} we get, as shown in detail in 
\textbf{Summary} above -- use of schema $(\tilde{\pi}_O),$
on absurdity of infinitely descending $\mrCCIO$'s ``in'' Ordinal
$O,$ \emph{contraposition of} and therefore equivalent to
schema $(\piO)$ -- the following 

\smallskip
\emph{Self-Consistency} \textbf{Corollary} for Theories $\piOR:$
  $$\piOR \derives\ \neg\,\Pro_{\piOR} (k,\code{\false}): \N \to 2:$$
Theory $\piOR,$ $O \succeq \N[\omega],$ \ul{derives} its own
-- Free-Variable -- (internal) \emph{non-Provability} of 
$\code{\false},$ \ie it \ul{derives} its own (Free-Variable)
\emph{Consistency Formula.}


\section{An Implicational, Local Variant 
                                                  of  Axiom of \emph{Descent} }

We consider an \ul{alternative} \emph{Descent} \textbf{axiom}
over $\PRa,$ namely the following \emph{implicational,} by that
\emph{equational} schema, to replace \emph{Descent} axiom
$(\piO),$ namely
\inference{ (\pidO) }
{ $c = c(a): A \to O$ (\emph{complexity}), \\
& $p = p(a): A \to A$ (\emph{``predecessor'' step}) \\
& $\chi = \chi(a): A \to 2$ \\
& \qquad 
    \text{(arbitrary) \emph{``test'' predicate} for circumscription 
                                               of ``$\exists\,n$'',} \\
& \text{logically: $\chi$ a $\ul{\free}$ meta-\ul{Variable} over 
                                     $\PRa$-predicates on $A$} 
}
{ $[\,[\,[\,\DeStad[\,c\,|\,p\,] (a,n) 
                \implies c\,p^n(a) \doteq 0_O\,]$ \\
& \qquad\qquad\qquad\qquad
    $\implies \chi(a)\,] \implies \chi(a)\,] = \true: 
                                      A \times \N \to 2:$ 
}
\emph{For ``each'' $a$ ``exists'' $n \in \N$ terminating $p^n$ 
into $c\,p^n(a) \doteq 0$, existence expressed ``locally'' via
2 implications, local at ``given'' $a \in A,$ and 
concerning ``test'' predicate} ($\ul{\free}$ predicate Variable) 
$\chi = \chi(a): A \to 2.$

\smallskip
\textbf{Definition} of \emph{individualised} \emph{Descent condition,} 
above, descent condition concerning ``only'' a ``given'', 
(finite) \emph{sequence} of length $n,$ starting at given $a:$
\begin{align*}
& \DeStad[\,c\,|\,p\,] (a,n) \defeq 
      \underset{n' \leq n} {\,\land\,} \DeSta\,[\,c\,|\,p\,] (p^{n'} (a)):
                                                A \times \N \to 2, \\
& \quad
    \text{where, \textbf{recall:}} \\
& \DeSta = \DeSta\,[\,c\,|\,p\,] (a) \bydefeq \\
& [\,c(a) > 0 \implies c\,p(a) < c(a)\,]
                                  \ \text{\ul{De}scent (\emph{main})} \\
& \land\ [\,c(a) \doteq 0 \implies p(a) \doteq_A a\,]
                       \ \text{\ul{Sta}tionarity (\emph{auxiliary})}
\end{align*} 
\textbf{Strengthening Remark:} This (equational) \textbf{axiom} 
infers ``original'' schema $(\piO)$ by inferential 
\ul{modus} p\ul{onens}: \ul{Antecedent} of $(\piO)$ makes true 
(first) \emph{premise} $\DeStad[\,c\,\,|\,p\,] (a,n)$ of $(\pidO)$'s 
\ul{Postcedent}, for $a \in A\ \free\ (!),$ and then gives -- by 
\emph{boolean Free Variables tautology} -- \ul{Postcedent}

\smallskip 
  $\pidOR\ \derives\ \chi(a) = \true_A: A \to 2,\ a \in A\ \free,
    \ \text{of schema}\ (\piO)\ \text{for theory}\ \pidOR.$
 
\medskip
We turn to (equivalent) Free-Variables \textbf{Contraposition} to 
\emph{local,} implicational schema $(\pidO).$ It reads:
\inference{ (\tildepidO) }
{ $c = c(a): A \to O,$ $p = p(a): A \to A$ in $\PRa$ ``given'', \\
& $\psi = \psi(a): A \to 2$
     (meta \ul{free}) \emph{``absurdity test'' predicate}  
}
{ $\pidOR\ \derives\ 
     [\,[\,\psi(a) \impl
       \DeStad[\,c\,|\,p\,] (a,n)\,\land\,c\,p^n(a) > 0\,]
                                       \impl \neg\,\psi(a)\,]:$ \\
& \qquad\qquad 
    $A \times \N \to 2.$
}
  %
  %
\textbf{Interpretation of $(\pidO)$ and $(\tildepidO):$}
\begin{enumerate} [(i)]
\item
Implicational schema $(\pidO)$ says intuitively: 
for any $a \in A$ ``given'', there \emph{``exists''} $n \in \N$ 
such that \emph{descent} $c\,p^0(a) > \ldots > c\,p^n(a)$ during 
$n$ steps, \emph{implies (stationary) termination} 
$c\,p^n(a) \doteq 0_O$ after $n$ steps.

\item
In particular: If chain $[\,c\,|\,p\,]$ satisfies earlier
descent condition $\DeSta\,[\,c\,|\,p\,] (a),$ mainly:
$c(a') > 0 \implies c\,p(a') < c(a')$
for all (consecutive) arguments of form $a' = p^{n'} (a),\ n' \leq n,$ 
``any'' $n$ given, then this chain must become \emph{stationary} after 
finitely many steps $n' \mapsto n'+1.$ All this 
\emph{individually, ``locally''} for $a \in A$ given. 

\item
If $[\,c\,|\,p\,]$ satisfies $\DeSta$ \emph{globally:}
for $a \in A$ free, then chain above must be stationary after 
finitely many steps for all $a$ (with termination index still
individual for each $a.$) This case is just (Interpretation of) 
\textbf{Strengthening Remark} above: $(\pidO)$ \ul{infers} $(\piO).$

\item  
(Equivalent) Free-Variables Contraposition $(\tildepidO)$ of 
$(\pidO):$ 
  $$[\,\psi(a) \impl [\,\DeStad(a,n)\,\land\,c\,p^n(a) > 0\,]\,]
                        \impl \neg\,\psi(a)\ \emph{interprets:}$$
\emph{$\DeSta\,[\,c\,|\,p\,] (p^n(a))$ for (individual) 
$a \in A$ and for all $n \in \N,$ but nevertheless
infinite descent at ``this'' $a,$ is absurd: 
any condition $\psi = \psi(a)$ on $A$ which implies that 
absurdity for the given $a,$ must be false on 
that $a.$}
\end{enumerate}

Theorie(s) $\pidOR = \PRa+(\pidO)$ now \textbf{inherit} directly 
all of the assertions on formally partial, $\hatPRa$ evaluation 
$\eps = \eps(u,a): \mrPRa \times \X \parto \X$ as well as
of \emph{argumented-deduction-tree evaluation} 
$\eps_d: \Stree \parto \Stree,$ with the following exceptions, 
where schema $(\piOR)$ enters explicitely:

\smallskip 
\textbf{Tree Argumentation, extra Case:} For this we need 
``abbreviation''
\begin{align*}
& \DeStad[\,c\,|\,p\,] (a,n): A \times \N \to 2, \\
& \quad
    \text{this predicate reads more formally:} \\
& \bydefeq \mrpr\,[\,\true: A \to 2,
                  \ b\,\land\,\DeSta\,[\,c\,|\,p\,] (p^{n'} (a))\,]: 
                                                  A \times \N \to 2.
\end{align*}
Here $b :\,= r_{A \times \N,2}: (A \times \N) \times 2 \to 2$ is right
projection, and \\
  $$\mrpr\,[\,g: A \to B,\ h: (A \times \N) \times B \to B\,]: 
                                                A \times \N \to B$$
is (unique) \textbf{definition} of a $\mrPRa$ map, out of \emph{anchor} 
$g$ and \emph{step} $h,$ by the \emph{full schema} $(\mrpr)$ of Primitive
Recursion.

Still more formally, without use of Free Variables, we have
\begin{align*}
& \DeStad[\,p\,|\,c\,] = 
    \mrpr[\,\true_A\,,\ r_{A \times \N,2}\,\land\,
                        [\,\DeSta\,[\,c\,|\,p\,] \circ p^{\S} 
                                 \circ \ell_{A \times \N,2}\,]\,]: \\ 
& A \times \N \to 2. 
\end{align*}
We \emph{internalise} this \emph{sequential descent,} $\DeStad,$ into
\begin{align*}
\destad(u,v) \defeq   
& \code{\mrpr}\,[\code{\true_A}; \code{r}\,\code{\land}\,
                            [\,\desta(u,v) \odot v^{\code{\S}} 
                                          \odot \code{\ell}\,]\,]: \\ 
& \cds{A,O} \times \cds{A,A} \to \cds{A \times \N,2}, 
\end{align*}
where $\desta = \desta(u,v)$ is internal version of 
$\DeSta\,[\,c\,|\,p\,]$ \textbf{defined} and used frequently
above: no change here.

This gives the following \textbf{type} of dummy argumented tree $t$ 
in the actual $\pidO\Case,$ with just one explicit level:

\bigskip
\cinferenceq{t\ =\quad }
{ $\an{\an{\an{\destad\,(u,v)\,\code{\impl}\,\an{\,u \odot v^{\code{\S}} 
                                        \code{\doteq 0}}}$ \\
& $\,\code{\impl}\,w}\,\code{\impl}\,w}/\Box \sim \code{\true}$ 
}    
{ $t'$ \hfill $\tilde{t}'$ }  

\bigskip
with
\emph{branches} 
$t',\,\tilde{t}' \in \dumTree \subset \Stree$
\emph{dummy argumented} Similarity trees.
  
In analogy to the other \emph{equational} cases (for theorie(s) $\piOR,$
we are led to \textbf{define} for $t$ the actual, \emph{argumented} form:

\bigskip
\cinferenceq{ t/\an{a;n} \defeq\quad }
{ $\an{\an{\an{\destad(u,v)\,\code{\impl}\,
            \an{u \odot v^{\code{\S}}/\an{a;n}\,\code{\doteq 0}}}$ \\
& $\,\code{\impl}\,w/a}\,\code{\impl}\,w/a}\,\sim\,\code{\true}$ 
}    
{ $t'/\an{a;n}$ \hfill $\tilde{t}'/\an{a;n}$ }

\bigskip
This completes \emph{tree argumentation,} by consideration of
the \textbf{final,} extra case, final case here treating 
schema $(\pidO)$ for theorie(s) $\pidOR,$ replacing original 
one(s) $(\piO),$ for theorie(s) $\piOR.$

\bigskip
\textbf{Definition} of \emph{map-code/argument} trees, $\Stree,$
of (PR) \emph{tree-complexity} $c_d: \Stree \to O$ as well as  
(PR) \emph{tree-evaluation step} $e_d: \Stree \to \Stree$ 
carry over -- suitably modified -- from theorie(s) $\piOR$
to present theorie(s) $\pidO.$ The same then is true for the
``finite'' \textbf{Descent} of \emph{map-code/argument tree}
evaluation $\eps_d: \Stree \parto \Stree.$ This $\eps_d$ is 
the $\mrCCIO$ \textbf{defined} by these (modified) complexity 
$c_d$ and iteration of step $e_d:$ iteration 
\emph{as long as complexity $0_O$ is not ``yet'' reached.}

\smallskip
From this we get, in analogy to that for theorie(s) $\piOR,$ the
(modified) 

\medskip
\emph{\textbf{Main Theorem}} for theorie(s) $\pidOR,$ again on 
\textbf{Termination-Conditioned Soundness:}

It is conceptually unchanged: replace \emph{Descent} Theory $\piOR$ 
by \emph{``even'' local} \emph{Descent} Theory $\pidOR,$ and read 
internal equality (enumeration) $\checkeq_k: \N \to \mrPRa^2$ 
as internal equality of $\pidOR$ (just this makes the difference.) 

\emph{Termination-Conditioned Inner Soundness} reads, 
for theories \\ $\pidOR = \PRa+(\pidO):$ 
  %
\begin{align*}
\pidOR \derives\
& [\,u\,\checkeq_k\,v\,]   
      \,\land\, [\,m\ \deff\ \eps(u,a),\ \eps(v,a)\,] \implies: \\
& \eps(u,a) \doteq r\ e^m(u,a) 
              \doteq r\ e^m(v,a) \doteq \eps(v,a), & (\bullet) \\  
& u,v \in \mrPRa,\ a \in \X,\ m \in \N\ \free.
\end{align*}        
\emph{Interpretation:} Unchanged, see \emph{Main Theorem} for 
theorie(s) $\piOR$ above.

Same for the \textbf{consequences:} 
\begin{itemize}
\item
\emph{Termination-Conditioned Objective Soundness for Map-Equality,} \\
which gives in particular

\item 
\emph{Termination-Conditioned Objective Logical Soundness:}
  $$\pidOR \derives\
      \Pro_{\pidOR} (k,\code{\chi})
         \,\land\,[\,m\ \deff\ \eps_d(\dtree_k/a)\,]
                           \implies \chi(a): \N^2 \times A \to 2.$$

\end{itemize}
\textbf{(Modified) Proof} of Termination-Conditioned 
\emph{Inner Soundness:}

There is no change necessary in all \textbf{Cases} except the 
\textbf{extra,} final case characterising theory $\piOR$ 
\resp $\pidOR:$
The standard, non-\textbf{extra} cases can be 
\textbf{proved} already within $\PRa,$ with
$u\,\checkeq_k\,v$ designating $\PRa$'s internal-equality enumeration,
as well when designating the \emph{stronger} ones of $\piOR$ \resp the 
still stronger ones of present theorie(s) $\pidOR.$

Remains to \textbf{prove} \emph{Termination-Conditioned Inner Soundness} 
for

\smallskip
\textbf{Extra Case} for theory $(\pidO),$ corresponding to 
its characteristic, \emph{extra} \textbf{axiom} $(\pidO).$
 
For this, \textbf{recall:} 
\begin{align*}
& \desta = \desta(u,v) \bydefeq \\ 
& \quad
    \an{u \code{>0}\,\code{\impl}\,u \odot v \code{<} u}
      \,\land\,\an{v \code{\doteq 0}
             \,\code{\impl}\,u \code{\doteq} \code{\id}}: \\
& \quad
    \cds{\X,O} \times \cds{\X,\X} \to \cds{X,2} = \cds{\X,2}_{\PRa}.
\end{align*}  
Free variable $w \in \cds{\X,2}$ is to internalise \emph{test} 
predicate $\chi: A \to 2.$ 

\smallskip
Finally \textbf{recall} from above completely formal internalisation 
\begin{align*}
& \destad (u,v): \cds{\X,O} \times \cds{\X,\X} \to \cds{\X \times \N,2}
                                                 \ \text{given by} \\
& \destad(u,v) \defeq  
      \code{\mrpr}\,[\code{\true}; \code{r}\,\code{\land}\,
                               [\,\desta(u,v) \odot v^{\code{\S}} 
                                          \odot \code{\ell}\,]\,]: \\ 
& \cds{\X,O} \times \cds{\X,\X} \to \cds{\X \times \N,2}. 
\end{align*}
What we have to \textbf{prove} in this case -- taking into account
just the only explicit equation in the corresponding deduction tree -- is
\begin{align*}
\pidOR\ \derives\ 
& m\ \deff\ \text{all $\eps$ terms below} \implies: \\
& [\,[\,\eps(\destad(u,v),\an{a;n}) 
           \implies [\,\eps(u \odot v^{\code{\S}},\an{a;n}) \doteq 0\,] \\ 
& \implies \eps(w,a)\,] \implies \eps(w,a)\,] \doteq \true: 
                                             & (\bullet^{{}_{\bullet}}) \\
& \N \times (\cds{\X,O} \times \cds{\X,X} \times \cds{\X,2}) 
                                         \times \an{\X \times \N} \to 2.
\end{align*}
For reduction of this case ``to itself'', we \textbf{define} here
-- in (simpler) parallel to the $\piOR$ setting --
a special \textbf{instance} for schema $(\pidO),$ ``consisting'' out of 
a ``super Domain'' $A_{\pi},$ a ``super complexity'' 
$c_{\pi}: A_{\pi} \to O,$ a ``super step'' $p_{\pi}: A_{\pi} \to A_{\pi},$
as well as a ``super test predicate'' $\chi_{\pi}: A_{\pi} \to 2,$
such that in fact ``finite descent'' is given -- and such that this
instance of $(\pidO)$ is able to \ul{derive} our assertion 
$(\bullet^{{}_{\bullet}})$ in present case. Here are the data 
for this instance:
\begin{align*}
A_{\pi} \defeq  
& \set{(m,(u,v,w),a) 
    \in \N \times (\cds{\X,O} \times \cds{\X,X} \times \cds{\X,2})
                                                        \times \X\,|\, \\
& \qquad\qquad\qquad\qquad\qquad
    m\ \deff\ \eps(u,a),\eps(v,a),\eps(\destad(u,v),a),\eps(w,a)} \\ 
& \subset \N \times \mrPRa^3 \times \X.   
\end{align*} 
Introduce Free Variable 
 $a_{\pi} \defeq (m,(u,v,w),a) \in A_{\pi} 
                      \subset \N \times \mrPRa^3 \times \X,$ \\
and \textbf{define}
\begin{align*}
& c_{\pi} = c_{\pi} (a_{\pi}) \defeq r\,e^m(u,a): A_{\pi} \to O,
   \ c_{\pi} (a_{\pi}) = \eps(u,a): A_{\pi} \to O\ \text{for short,} \\
& \quad
    (\text{termination property of $m$ ``fixed'' within}\ 
                                              a_{\pi} \in A_{\pi}.) \\
& p_{\pi} (a_{\pi}) = p_{\pi} (m,(u,v,w),a) 
    \defeq (m,(u,v,w),\eps(v,a)): A_{\pi} \to A_{\pi}.
\end{align*}
Finally, externalised ``super test predicate'' is 
taken, suitable for actual \textbf{proof,}
  $$\chi_{\pi} = \chi_{\pi} (a_{\pi}) = \chi(m,(u,v,w),a)
          = \eps(w,a) \bydefeq r\,e^m(w,a): A_{\pi} \to 2.$$
These fixed, next step is calculation of $\DeSta$ for above ``super'' data:
\begin{align*}
& \DeSta[\,c_{\pi}\,|\,p_{\pi}\,]\,(a_{\pi}) \\
& = [\,c_{\pi} (a_{\pi}) > 0_O \implies  
                 c_{\pi}\,p_{\pi} (a_{\pi}) < c_{\pi} (a_{\pi})\,]
                                                     & (\emph{Descent}) \\
& \quad
    \,\land\,[\,c_{\pi} (a_{\pi}) \doteq 0 \implies 
               c_{\pi} (a_{\pi}) \doteq a_{\pi}\,].  & (\emph{Stationarity})
\end{align*}
By \textbf{definition} of these data, this calculation gives:
\begin{align*}
& \DeSta[\,c_{\pi}\,|\,p_{\pi}\,]\,(a_{\pi}) \\
& = [\,m\ \deff\ \text{all instances of $\eps$ below}\,]\,\,\land\,: \\ 
& \quad 
    [\,\eps(u,a) > 0_O \implies \eps(u,\eps(v,a)) < \eps(u,a)\,] \\
& \quad\ 
    \,\land\,[\,\eps(u,a) \doteq 0 \implies \eps(v,a) \doteq_A a\,]:  
                 \N \times \mrPRa^3 \times \N \supset A_{\pi} \to 2.
\end{align*}
But this is equality between ($m$-dominated) iteration predicates
\begin{align*}
& \DeSta\,[\,c_{\pi}\,|\,p_{\pi}\,]\,(m,(u,v,w),a) \implies: \\  
& [\,m\ \deff\ \eps(\destad(u,v),a)\,] \\
& \qquad
    \land\,\DeSta\,[\,c_{\pi}\,|\,p_{\pi}\,]\,(m,(u,v,w),a)
                                    \doteq \eps(\destad(u,v),a): \\ 
& \N \times (\cds{\X,O} \times \cds{\X,\X} \times \cds{\X,2}) 
                                                     \times \X \to 2,
\end{align*}
We \emph{Objectivise} internal continous descent $\desta(u,v),$ 
via evaluation $\eps$ on $\an{a;n} \in \an{\X;\N}:$ we expect to get 
just instance $\DeStad[\,c_{\pi}\,|\,p_{\pi}]\,\an{a;n}$ of 
\emph{Objective sequential Descent:}
\begin{align*}
& m\ \deff\ \text{all $\eps$ terms in $(\bullet^{{}_{\bullet}})$ \emph{implies:}} \\
& m\ \deff\ \text{all $\eps$ terms below}\ \,\land\,: \\
& \eps(\destad(u,v),\an{a;n}) \\
& \doteq \eps(\code{\mrpr}\,[\,\code{\true_\X};\code{r}\,\code{\land}\,
           [\,\desta(u,v) \odot v^{\code{\S}} \odot \code{\ell}\,]\,],
                                                              \an{a;n}) \\
& \doteq \eps(\code{\underset{n' \leq n} {\,\land\,}}
                \desta(u,v) \odot v^{\code{\S}},\an{a;n'}) \\ 
& \doteq \underset{n' \leq n} {\,\land\,} 
           \eps(\desta(u,v),\eps(v^{\code{\S}},\an{a;n'})) \\ 
& \doteq \underset{n' \leq n} {\,\land\,} 
             \eps(\desta(u,v),p_{\pi}^{n'} (m,(u,v,w),a) \\
& \qquad
    \text{with}\ a_{\pi} :\,= (m,(u,v,w),a),\ p_{\pi}^{n'} (a_{\pi}) 
       \in A_{\pi} \subset \N \times \mrPRa^3 \times \X,
                                       \ \text{for}\ n' \leq n \\
& \bydefeq \underset{n' \leq n} {\,\land\,} 
             \DeSta\,[\,c_{\pi}\,|\,p\,]\,(p_{\pi}^{n'} (a_{\pi})) \\
& \bydefeq \DeStad[\,c_{\pi}\,|\,p_{\pi}\,]\,(a_{\pi},n) \\
& = \DeStad[\,c_{\pi}\,|\,p_{\pi}\,]\,((m,(u,v,w),a),n): \\ 
& \N \times (\cds{\X,O} \times \cds{\X,\X} \times \cds{\X,2})
                                          \times \an{A \times \N} \to 2.
\end{align*} 
This is wanted externalisation
\begin{align*}
& m\ \deff\ \text{all $\eps$ terms in $(\bullet^{{}_{\bullet}})$ \emph{implies:}} \\
& \eps(\destad(u,v),\an{a;n})
    \doteq \DeStad[\,c_{\pi}\,|\,p_{\pi}\,]\,((m,(u,v,w),a),n): 
                                                      & (\eps\ \desta) \\
& \N \times (\cds{\X,O} \times \cds{\X,\X} \times \cds{\X,2}) \to 2. 
\end{align*}
This given, we attempt, again by Objectivisation via $\eps$ of 
$(\bullet^{{}_{\bullet}}),$ to show the \emph{``finite'' descent} property for our 
\textbf{instance} $A_{\pi}$ \etc, \ie essentially for $\DeStad,$ 
as follows:
\begin{align*}
& m\ \deff\ \text{all $\eps$ terms in $(\bullet^{{}_{\bullet}})$ \emph{implies:}} \\
& [\,[\,\DeStad[\,c_{\pi}\,|\,p_{\pi}\,]\,(a_{\pi},n) 
        \implies \chi_{\pi} (a_{\pi})\,] \implies \chi_{\pi} (a_{\pi})\,] \\
& = [\,[\,\DeStad[\,c_{\pi}\,|\,p_{\pi}\,]\,((m,(u,v,w),a),n) 
        \implies \eps(w,a)\,] \implies \eps(w,a)\,] \\
& \doteq [\,[\,\eps(\destad(u,v),\an{a;n}) 
                 \implies \eps(w,a)\,] \implies \eps(w,a)\,]: 
                            \qquad (\text{just $(\bullet^{{}_{\bullet}})$}) \\
& \N \times (\cds{\X,O} \times \cds{\X,\X} \times \cds{\X,2})
                                          \times \an{A \times \N} \to 2. 
\end{align*}
This shows that our hypothesis $(\bullet^{{}_{\bullet}})$ is equivalent to
\emph{``finite'' sequential descent} of \textbf{instance}
$\bfan{\bfan{A_{\pi}, c_{\pi}, p_{\pi}},\chi_{\pi}}.$ 

But this is an instance ``for'' \textbf{axiom} $(\pidOR)$ of our 
Theory $\pidOR = \PRa+(\pidO).$ So that axiom shows remaining
assertion $(\bullet^{{}_{\bullet}}),$ \emph{Inner Soundness} for the final,
\emph{``self-referential''} case. This \textbf{proves} the
\textbf{Main Theorem} for theorie(s) $\pidOR.$ 

\smallskip
By use of (contrapositive) characteristic schema $(\tildepid_O)$
of theory $\pidOR = \PRa+(\pidO)$ (absurdity of infinitely descending 
iterative $O$-chains), we get -- in complete analogy to the \textbf{proof} 
for theorie(s) $\piOR$ in \textbf{Summary} above:  

\medskip
\emph{Self-Consistency} \textbf{Corollary} for Theories $\pidOR:$
  $$\pidOR \derives\ \neg\,\Pro_{\pidOR} (k,\code{\false}): \N \to 2,
                                                  \ k \in \N\ \free:$$
Theory $\pidOR,$ $O \succeq \N[\omega],$ \ul{derives} its own
-- Free Variable -- (internal) \emph{non-Provability} of 
$\code{\false},$ \ie it \ul{derives} its own (Free Variable)
\emph{Consistency Formula.}


\section{Unconditioned Objective Soundness}

As is well known, Consistency Provability and Soundness are strongly tied
together. Above we have shown that already \emph{Termination-Conditioned}
Soundness entails Consistency Provability. Here we ``easily''
\ul{derive} Full, Unconditioned Objective (!) Soundness from Consistency 
Provability, for all of our \emph{Descent Theories} $\bfPi,$
strengthenings of $\PRa,$ $\bfPi$ standing from now on for one arbitrary 
such theory, namely $\piOR$ of \emph{on}-terminating 
\emph{Complexity Controlled Iterations,} or $\pidOR$ of 
``\ond-terminating'' $\mrCCIO$'s, with complexity values in Ordinal $O,$ 
$O$ one of the (Order) extensions of Ordinal $\N[\omega]$ 
introduced above, \ie one of
$\N[\omega],\ \N[\xi_1,\ldots\,\xi_{\ul{m}}],\ \X,\ \text{and}\ \E.$

\smallskip
We start with the observation that \emph{Consistency}(-formula) 
\ul{Derivabilit}y \\
$\bfPi \derives\ \neg\,[\,0\,\checkeq\,1\,]: \N \to 2$
is equivalent to \ul{derivabilit}y
  $$\bfPi \derives\ [\,\nu_2(a)\,\checkeq_k\,\nu_2(b)\,] 
             \implies a \doteq b: \N \times (2 \times 2) \to 2: 
                                                       \ (\bs{*})$$
Test with $(a,b) \in \set{(0,0),(0,1),(1,0),(1,1)}.$ Cases
$(0,1)$ and $(1,0)$ are (each) just \emph{Consistency} 
\ul{derivabilit}y, the remaining two are trivial.

Formally this test is based on the fact, that 
  $$(0,0),\ (0,s\,0),\ (s\,0,0),\ (s\,0,s\,0): \one \to 2 \times 2$$
are the 4 coproduct injections of coproduct (sum) 
$2 \oplus 2 \defeq 2 \times 2.$  

\smallskip
Now $(\bs{*})$ is -- by \textbf{definition} -- just \emph{injectivity} 
of \emph{internal numeralisation} 
  $$\nu_2 = \nu_2(a): 2 \to [\one,2]_{\bfPi} 
                = \cds{\one,2}_{\PRa}/\checkeq^{\bfPi}.$$
This \emph{numeralisation} is defined within general Arithmetical
theories by
\begin{align*} 
& \nu_\N = \nu(n): \N \to [\one,\N] = \cds{\one,\N}/\checkeq 
                                      \ \text{PR as follows:} \\
& \nu(0) \defeq \code{0}: \one \to [\one,\N], \\ 
& \nu(s\,n) \defeq \code{s} \odot \nu(n): \N \to [\one,\N],
                                   \ \text{whence in particular:} \\
& \nu(\num(\ul{n})) = \code{\num(\ul{n})} = \code{s \ldots s \circ 0} \\
& \text{for external numeralisation}
     \ \num: \ul{\N} \bs{\lto} \bfS(\one,\N).
\end{align*}  
Further -- externally \ul{PR}: 
\begin{align*}
& \nu_{A \times B} = \nu_{A \times B} (a,b) \defeq \an{\nu_A(a);\nu_B(b)}: \\
& A \times B \to [\one,A] \times [\one,B] \ovs{\iso} [\one,A \times B]. 
\end{align*}
For an abstraction Object $\set{A\,|\,\chi},$ 
           as in particular $2 = \set{\N\,|\ < s\,0},$ \\
$\nu_{\set{A\,|\,\chi}}$ is defined by (double) restriction, 
of $\nu_A: A \to [\one,A].$

\smallskip
\textbf{Naturality Lemma for Internal Numeralisation:} For each $\bfPi$
map ($\PRa$ map) $f: A \to B$ the following \textsc{diagram}
commutes -- in category $\bfPi{Q} = \bfPi+\mr{Quot} \bs{\sqsup} \bfPi:$ 
Theory $\bfPi$ enriched by (virtual) Quotients by equivalence 
\emph{Relations,} such as in particular 
$\,\checkeq\ =\ \checkeq_k: \N \to \cds{\X,\X}^2:$

\bigskip
\begin{minipage} {\textwidth}
\xymatrix @+2em{ 
& A
  \ar[rr]^f
  \ar[d]^{\nu_A}
  \ar @{} [drr]|=
  & & B
      \ar[d]^{\nu_B}
\\
\cds{\one,A}/\checkeq
\ar @{=} [r]
& [\one,A]
  \ar[rr]^{ [\one,f] }
  & & [\one,B]
      \ar @{=} [r]
      & \cds{\one,B}/\checkeq
}
\end{minipage}

\bigskip
\textbf{Proof:} We have to show equality in the following 
Free-Variable setting which displays the assertion, by 
\textbf{definition} of functor $[\one,f]: [\one,A] \to [\one,B]:$

\bigskip
\begin{minipage} {\textwidth}
\xymatrix{ 
{A \owns a} 
\ar @{|->} [rrr]^f
\ar @<2.0 ex> @{|->} [d]^{\nu_A} 
& & & {f(a) \in B}
      \ar @<-2.0 ex> @{|->} [d]^{\nu_B}
\\
{ [\one,A] \owns \nu_A(a) }
\ar @{|->} [rr]^{ [\one,f] } 
& & {\ncode{f} \odot \nu_A(a)}
    \ar @{} [r]|\checkeq
    & \nu_B(f(a)) \in [\one,B]
}              
\end{minipage}

\bigskip
This internal equality $\ncode{f} \odot \nu_A(a)\ \checkeq\ \nu_B(f(a))$
\ \ is \textbf{proved} straightforward by external structural 
\ul{recursion} on the structure of $f: A \to B$ in $\PRa,$ beginning
with the maps constants $0,\ s,\ \ell,$ using internal associativity
of ``$\odot$'', and (objective) PR on the iteration count for the case 
of an iterated.
 
\smallskip
\textbf{Injectivity Lemma for Internal Numeralisation:} Injectivity
of $\nu_2: 2 \to [\one,2]_{\bfPi},$ given by Consistency 
\ul{derivabilit}y, extends to injectivity
of all $\nu_A = \nu_A(a): A \to [\one,2],$ first to 
$\nu_\N = \nu(n): \N \to [\one,\N]$ essentially by considering 
truncated subtracction, and then immediately to the other
Objects of $\PR$ and $\PRa.$

\smallskip
This leads to our final result here, namely

\smallskip   
\textbf{(Unconditioned) Objective Soundness Theorem} for $\bfPi:$
\begin{itemize}
\item For each pair $f,\ g: A \to B$ of $\PRa$-maps:
  $$\bfPi \derives\ [\,\code{f}\,\checkeq_k\,\code{g}\,] 
             \implies [\,f(a) \doteq_B g(a)\,]: \N \times A \to 2,$$
whence by specialision:

\item For each $\PRa$ predicate $\chi = \chi(a): A \to 2:$
  $$\bfPi \derives\ \Pro_{\bfPi} (k,\code{\chi}) 
                      \implies \chi(a): A \to 2:$$
Availability of an (Internal) \emph{Proof} of (code of) a predicate 
implies \emph{truth} of this predicate at each argument. 
\end{itemize}

\textbf{Proof} of first \textbf{assertion:} Consider the following
commutative \textsc{diagram} -- in Theory $\bfPi{Q} \bs{\sqsup} \bfPi:$

\bigskip
\begin{minipage} {\textwidth}
\xymatrix @+2em{ 
A
\ar[rr]^f_g
\ar[d]^{\nu_A}
\ar @{} [drr]|=
& & B
    \ar[d]^{\nu_B}
\\
[\one,A]
\ar[rr]^{ [\one,f] }_{ [\one,g] }
& & [\one,B]
}
\end{minipage}

\bigskip
This gives
\begin{align*}
\bfPi \derives\ 
& [\one,f]\,(\nu_A(a))\ [\ \bydefeq \code{f} \odot \nu_A(a)\ ] \\
& \quad \,\checkeq_{j(k,a)}\,\code{g} \odot \nu_A(a) 
            \quad(\text{by hypothesis $\code{f}\,\checkeq_k\,\code{g}$}), \\
& \implies (\nu_B \circ f) (a) = (\nu_B \circ g) (a) 
                       \quad \text{by commutativity above} \\
& \implies f(a) \doteq_B g(a): \N \times A \to B^2 \to 2 \\
              \ \text{by injectivity of}\ \nu_B. 
\end{align*}

This taken together gives first -- and then second -- 
assertion of the Theorem\quad\textbf{\qed}

\smallskip
Analysis of \textbf{Proof} above shows that we can take (internal)
\emph{Consistency} as an additional condition for a an arithmetical
theory $\bfS$ instead using it as \ul{derived} property of our
(self-consistent) theories $\bfPi.$ This then gives, for such
general theory $\bfS,$ with $\bfS^+ \defeq \bfS+\Con_{\bfS}:$

\smallskip
\textbf{Consistency Conditioned Injectivity of Internal Numeralisation:}
  $$\bfS^+ \derives\ \nu_A(a)\,\checkeq^{\bfS}_k\,\nu_A(a') 
              \implies a \doteq_A a': \N \times A^2 \to 2.$$
$[\,$Note the difference between frame $\bfS^+$ and internal equality
taken within weaker theory $\bfS$ itself.$]$

\smallskip
\textbf{Consistency Conditioned Soundness:}
\begin{itemize}
\item 
for $\PRa$-maps $f,\ g: A \to B:$
  $$\bfS^+ \derives\ [\,\code{f}\,\checkeq^{\bfS}_k\,\code{g}\,]
                      \implies f(a) \doteq_B f(b): \N \times A \to 2.$$
\item
in particular for a predicate $\chi = \chi(a): A \to 2:$
  $$\bfS^+ \derives\ \Pro_{\bfS} (k,\code{\chi}) \implies \chi(a): 
                                                 N \times A \to 2.$$
Again: Here (internal) $\bfS$-\emph{Provability} is the premise.
It coincides with \emph{Provability} of frame $\bfS^+$ only for
self-consistent $\bfS,$ as for example for theorie(s) 
$\bfPi = \bfPi^+$  considered above.
\end{itemize}
  
\smallskip
(Conditioned) injectivity of internal numeralisation, and 
naturality invite to consider an \ul{inferential} form 
of (conditioned) $\omega$-Completeness:


\smallskip
\textbf{$\omega$-Completeness Theorem, Inference Form:}
\begin{itemize}
\item
Strengthenings $\bfS$ of $\PRa$ are \emph{Consistency-conditioned}
$\omega$-\emph{inference-complete,} \ie

\inference{ (\mr{\ul{Com}p}_{\omega}^{\bfS/\bfS^+}) }
{ $\chi = \chi(a): A \to 2$ in $\PRa,$ \\
& $k = k(a): \N \to \emph{Proof}_{\bfS}$ in $\PRa,$ \\
& $\bfS^+ \derives\ \Pro_{\bfS} (k(a),\code{\chi} \odot \nu_A(a)): 
                                                            A \to 2$ 
}
{ $\bfS^+ \derives\ \chi: A \to 2.$ }

\item 
Axis case: Self-consistent theories $\bfPi$ are (``unconditioned'') 
\emph{inferential $\omega$-self-complete,} they admit the special 
schema derived from the above:

\inference{ (\mr{\ul{Com}p}_{\omega}^{\bfPi}) }
{ $\chi = \chi(a): A \to 2$ in $\PRa,$ \\
& $k = k(a): \N \to \emph{Proof}_{\bfPi}$ in $\PRa,$ \\
& $\bfPi \derives\ \Pro_{\bfPi} (k(a),\code{\chi} \odot \nu_A(a)): 
                                                            A \to 2$
}
{ $\bfPi \derives\ \chi: A \to 2,$ and hence, by internalisation: \\
& $\bfPi \derives\ \Pro_{\bfPi} (k[\chi],\code{\chi}): \one \to 2,$ \\
& $k[\chi]: \one \to \emph{Proof}_{\bfPi}$ the code of 
                                     $\bfPi$ \ul{Proof} of $\chi.$
}
\end{itemize}

$[\,$The latter \emph{internalisation} of $\bfPi-\ul{derivation}$ 
of $\chi$ into an (internal) \emph{Proof} of $\bfPi$ itself for 
$\code{\chi}$ is decisive: it works because of self-consistency 
$\bfPi = \bfPi^+.$ Schema $(\mr{\ul{Com}p}_{\omega}^{\bfPi}),$
with last poscedent, almost says that $\one$ is a \emph{separator}
Object for internalised theory $\bfPi:$ test with all internal
points, even: with all internal \emph{numerals,} establishes 
internal equality, at least for ``concrete'' code pairs 
$\code{f}, \code{g} \in \cds{A,B},$ coming coded from objective
map pairs $f,\ g: A \to B$ of $\bfPi.]$

\smallskip
\textbf{Proof:} Look at $\nu$-naturality \textsc{diagram} in foregoing
section, and take special case $\chi: A \to 2$ for $f: A \to B.$ 
Then consider Free-Variable \textsc{diagram} chase for this $f,$
subsequent \textsc{diagram.} By commutativity of that rectangle 
we have
  $$\code{\chi} \odot \nu_A(a)\,\checkeq^{\bfS}_{j(a)}\,\nu_2(\chi(a)),$$
suitable $j = j(a): A \to \emph{Proof}_{\bfS} \subset \N.$ But by 
antecedent, we have also 
\begin{align*}
& \code{\chi} \odot \nu_A(a)\,\checkeq^{\bfS}_{k(a)}\,\code{\true}, 
                                                     \ \text{whence} \\
& \nu_2(\chi(a))\,\checkeq^{\bfS}_{k'(a)}\,\code{\true} = \nu_2(\true).
\end{align*}
(Consistency conditioned) \emph{injectivity} of internal numeralisation 
$\nu$ then gives $\chi(a) \doteq true,$ $a \in A$ free. Taken 
together: Given the antecedent $\bfS^+$ \ul{derivation}, we get
$\bfS^+ \derives\ \chi(a): A \to 2,\ a \in A\ \free.$ This is
what we wanted to show.

The ``axis'' case of a self-consistent theory, such as $\bfPi,$ then is
trivial, and gives \emph{(Unconditioned) inferential 
$\omega$-Completeness.}



\section*{Coda: Termination Conditioned Soundness 
                                                                            for Theory $\PRa$}

Termination-conditioned (!) (Objective) Soundness holds ``already'' 
for \emph{basic} PR Theory $\PRa,$ and hence also for its embedded 
Free-Variables \emph{fundamental} (categorical) Theory 
$\PR \bs{\sqsub} \PRa.$
The argument is use of following \textbf{Reduction} schema $(\rho_O)$ 
of predicate-truth, \emph{Reduction} ``along'' a given $\mrCCIO.$

Eventually we will \textbf{prove} by this schema of $\PRa\,(!)$
\ul{Consistenc}y of \emph{Descent} Theories $\bfPi$ \ul{relative}
to $\PRa.$

\bigskip
\textbf{Theorem:} Theory $\PRa$ admits the following \textbf{Schema} of  

\emph{Reduction} along $\mrCCIO$'s for \emph{Ordinal} $O$:
\inference{ (\rho_O) }
{ $[\,c: A \to O\,|\,p: A \to A\,]$ is a $\mrCCIO$ in $\PRa,$ \\
& $\chi = \chi(a): A \to 2$\ $\PRa$-predicate to be investigated, \\
& $\PRa \derives\ c(a) \doteq 0_O \implies \chi(a): A \to 2$
                                        \emph{predicate anchor,} \\
& $\PRa \derives\ \chi(p(a)) \implies \chi(a): A \to 2$ 
                            \ \emph{reduction step}
}
{ $\PRa \derives\ [\,m\ \deff\ \wh_O[\,c\,|\,p\,]\,(a) \implies \chi(a)\,]: 
                                                     A \times \N \to 2.$
}
\ul{Postcedent} meaning: 
\emph{Termination-of-$\while$-\emph{loop} conditioned} truth of $\chi(a),$
``individual'' $a.$

\smallskip
\textbf{Proof} by (Free-Variables) Peano induction on free variable
$m \in \N:$

\emph{Anchor} $m \doteq 0:$ \text{obvious by Antecedent}\ $(\anchor).$

\smallskip
Induction ``hypothesis'' on $m:$\quad 
  $m\ \deff\ \mu_O[\,c\,|\,p\,]\,(a) \implies \chi(a).$

\smallskip 
\emph{Peano Induction Step:}
\begin{align*}
  \PRa \derives\ 
  & m+1\ \deff\ \mu_O[\,c\,|\,p\,]\,(a') \\ 
  & \implies m\ \deff\ \mu_O[\,c\,|\,p\,]\,(p(a')) \doteq m \\
  & \phantom{\implies}\quad
      \text{by iterative definition of}\ \mu_O[\,c\,|\,p\,] \\
  & \implies \chi(p(a')) \ \text{by induction hypothesis} \\
  & \implies \chi(a'): A \times \N \to 2,
\end{align*}
the latter by \textbf{Antecedent} \emph{Reduction step}
\ \textbf{\qed}

\smallskip
For \textbf{Proof} of \emph{Termination-Conditioned Objective Soundness} 
of $\PRa$ by itself, we now consider the following instance of this 
Reduction schema $(\rho_{\rO})$ of $\PRa:$

\begin{itemize} 
\item
\emph{Domain} $\rA \defeq \N \times \Stree = \N \times \Stree_{\PRa},$
$\Stree$ above without the additional data coming in by schema
$(\piO)$ with its ``added'' (internal) deduction structure.
  
\item
\emph{Ordinal} $\rO \defeq \N \times \N[\omega]$ with 
hierarchical order: first priority to left component.

\item
\emph{``Predecessor'' step}
$p :\,= \re = \re(m,t) \defeq (m \dotminus 1,e_d(t)): \rA \to \rA,$ \\ 
(deduction) tree evaluation $e_d$ above, again ``truncated''
to the (internal) deduction data of $\PRa.$

\item
\emph{Tree complexity} 
$\rc = \rc(m,t) \defeq (m,c_d(t)): \rA \to \rO,$ 
$\PRa$ truncation as for $\re$ above.

\item
Finally the predicate to be \emph{reduced} 
                        with respect to its \emph{truth:} 

\smallskip
$\rph = \rph(m,t) 
   \defeq [\,m\ \deff\ \eps(\emph{root}_{\ell} (t))
                       \doteq \eps(\emph{root}_r (t))\,]:$ \\
$\N \times \Stree \to 2 \times \X^2 \xto{2 \times \doteq} 
                                       2 \times 2 \xto{\land} 2.$
\end{itemize} 
Here $\emph{root}_{\ell} (t)$ and $\emph{root}_r(t)$ are the left and
right entries, of form $u/x$ \resp $v/y,$ of 
$\emph{root} (t) = \an{u/x \sim v/x}$ say.

\smallskip
\textbf{Verification} of this instance of reduction schema 
$(\rho_{\rO})$ is now as follows:

\smallskip
\emph{Anchoring:} 
\begin{align*}
\PRa \derives\ 
& \rc(m,c_d(t)) \doteq (0,0) \implies : \\ 
& \rph(m,t) \doteq [\,0\ \deff\ \eps(\code{\id}/x \doteq \eps(\code{\id}/y)
      \doteq [\,x \doteq y\,] \doteq \true, \\
& \text{the latter necessarily for (flat) legimate $t$ of this form.} 
\end{align*}

\smallskip
\emph{Reduction Step} for $\rph:$
\begin{align*}
\PRa \derives\ 
& \rph\ \re(m,t) \bydefeq 
   [\,m \dotminus 1\ \deff\ \eps(\emph{root}_{\ell}\ e_d(t)) 
                              \doteq \eps(\emph{root}_r\ e_d(t))\,] \\
& \implies [\,m \ \deff\ \eps(\emph{root}_{\ell} (t)) 
                           \doteq \eps(\emph{root}_r(t))\,].
\end{align*}
This implication is \textbf{proved} -- logically -- by recursive 
case distinction on the two surface levels of $t,$ cases given in 
the main text above, the $(\piO)$ case truncated. Formally, this 
recursion is PR on (minimal) number $m$ of steps $e_d$ for complete 
tree evaluation of $t.$

\smallskip
Out of this \textbf{Antecedent,} schema $(\rho_{\rO})$ gives as its

\textbf{Postcedent} 
\begin{align*}
\PRa \derives\ 
& [\,m\ \deff\ \wh_{\rO} [\,\rc\,|\,\re]\,(m',dtree^{\mrPRa}_k/x)\,] 
                                                          \implies: \\ 
& [\,m'\ \deff\ \eps(\emph{root}_{\ell} (\dtree^{\mrPRa}_k/x)) 
                         \doteq \eps(\emph{root}_r(\dtree^{\mrPRa}_k/x))\,]: \\
& \N^2 \times (\N \times \X) \to 2,\ m,m',k \in \N,\ x \in \X\ \free,
\end{align*}
in particular, with $m :\,= m':$
\begin{align*}
\PRa \derives\ 
& [\,m\ \deff\ \wh_{\N[\omega]} 
    [\,c_d\,|\,e_d\,]\,(dtree_k/x)\,] \implies: \\ 
& [\,m\ \deff\ \eps(\emph{root}_{\ell} (\dtree_k/x)) 
                         \doteq \eps(\emph{root}_r(\dtree_k/x))\,]: \\
& \N \times (\N \times \X) \to 2,\ m,k \in \N,\ x \in \X\ \free.
\end{align*}
This is in fact 

\smallskip
\textbf{Termination-Conditioned Soundness Theorem} for \emph{basic}
PR Theory $\PRa,$ which holds by consequence also for 
\emph{fundamental} PR Theory $\PR \bs{\sqsub} \PRa.$

\smallskip
Can we reach from this \emph{Self-Consistency} for $\PRa$ as well,
in the manner we have got it for theorie(s) 
$\piOR = \PRa+(\piO) = \PRa+(\tilde{\pi}_O)?$ 

If you look at this derivation in the \textbf{Summary} above, you 
find as the final, decisive step, \ul{inference} from
\begin{align*}
\piOR \derives\ 
& \code{\false}\,\checkeq_k\,\code{\true} 
        \implies c_d\ e_d^m(\dtree_k/0) > 0: \N^2 \to 2,\ \text{\ul{to}} \\
\piOR \derives\ 
& \neg\,[\,\code{\false}\,\checkeq_k\,\code{\true}\,]: \N \to 2,
                                       \ k \in \N\ \text{free (!)}.
\end{align*}
This comclusion gets its \emph{legitimacy} by application of 
schema $(\tilde{\pi})$ to its suitable Antecedent with in particular
\emph{absurdity condition} $\psi$ -- for \emph{infinite} descent --
choosen as  
  $$\psi = \psi(k) :\,= [\,\code{\false}\,\checkeq_k\,\code{\true}\,]: 
                                                          \N \to 2.$$
Same for a general one out of theories $\bfPi,$ namely $\bfPi$ one of
$\piOR,\ \pidOR.$

\smallskip
If such -- formal, axiomatic -- absurdity of infinite descent is
\emph{not} available in the theory, infinite descent of in particular
$c_d\ e_d^m(dtree_k/0) > 0\ \text{(``for all'' $m$)}$ could not
be excluded: internal provability $\code{\false}\,\checkeq\,\code{\true}$
could ``happen'' formally by just ``the fact'' that (internal) 
\emph{deduction tree} for (internal) \emph{Theorem} 
$\code{\false}\,\checkeq\,\code{\true}$ cannot be externalised, 
by (iterative) deduction tree evaluation $\eps_d,$ in a finite 
number of its steps $e_d.$

So, in this sense, addition of highly plausible schema $(\tildepi)$
\resp $(\tildepid)$ is ``necessary'' -- at least it is sufficient -- 
for \ul{derivation} of (internal) \emph{Consistency,} this already 
for derivation of internal Consistency of Theory $\PRa.$

This latter result is not that astonishing, since Theory
$\piR = \pi_{\N[\omega]}\bfR$ is stronger than $\PRa,$ 
at least \ul{formall}y. Not to expect -- the G\"odel Theorems -- 
was finding of any \emph{Self-Consistent} 
(necessarily \emph{arithmetical}) theory, here theorie(s) $\bfPi,$ 
$\bfPi$ one of $\piOR,\ \pidOR,$ $O \succeq \N[\omega]:$ 

The most involved cases in the \textbf{proofs} leading to this 
Self-Consistency for theorie(s) $\bfPi$ -- in particular in (the two) 
Main Theorem(s) on \emph{Termination-Conditioned Inner Soundness,} 
and in the constructions leading to the notions used -- all come 
from ``this'' additional schema $(\Pi),$ schema $(\Pi)$ one of 
the schemata $(\piO)$ and $(\pidO)$ which constitute theorie(s) 
$\bfPi$ as (``pure'') strengthenings of $\PRa \bs{\sqsup} \PR.$

``Same'' discussion for (Unconditioned) \emph{Objective Soundness}
for $\bfPi,$ \ul{derived} in the above from \emph{Self-Consistency.} 
Conversely, this \emph{Objective Soundness} contains Self-Consistency 
as a particular case.

\smallskip
\textbf{Problem:} Is Theory $\piR,$ more general: are theories
$\bfPi$ (Objectively) Consistent \ul{relative} to \emph{basic} 
Theory $\PRa,$ and -- by that -- relative to \emph{fundamental} 
Theory $\PR \bs{\sqsub} \PRa$ of Primitive Recursion ``itself''?

In other words (case $\piR$): do \emph{Descent} data 
$c: A \to O :\,= \N[\omega],$ $p: A \to A,$ and availability of 
a $\PRa$ \emph{point} $a_0: \one \to A$ such that 
\begin{align*}
\PRa\ \derives\ 
& c\,p^\S(a_0,n) > 0_O: \\ 
& \one \times \N \xto{a_0 \times \N} A \times \N \xto{p^\S} A
               \xto{c} O \xto{\,> 0_O} 2,
\end{align*}
($n \in \N\ \free,$ intuitively: \emph{for all} $n \in \N:$ \ul{derived} 
\emph{non-termination} at $a_0$),
lead to a \ul{contradiction} within Theory $\PRa\,?$

\smallskip
We will take up this (relative) \textbf{Consistency Problem} again 
in terms of (recursive) \emph{Decision,} RCF\,5.

                   $$=====$$
 

\section*{References}



  
\smallskip
  
  \quad\NAME{J.\ Barwise} ed. 1977: \emph{Handbook of Mathematical Logic.}
  North Holland.
  
  \NAME{H.-B.\ Brinkmann, D.\ Puppe} 1969: 
  \emph{Abelsche und exakte Kategorien, Korrespondenzen.} 
  L.N. in Math. \textbf{96.} Springer.



 

  \NAME{S.\ Eilenberg, C.\ C.\ Elgot} 1970: \emph{Recursiveness.}
  Academic Press.
  
  \NAME{S.\ Eilenberg, G.\ M.\ Kelly} 1966: Closed Categories. 
  \emph{Proc.\ Conf.\ on Categorical Algebra}, La Jolla 1965, pp. 421-562. 
  Springer. 
  
  
  \NAME{G.\ Frege} 1879: \emph{Begriffsschrift.} Reprint in ``Begriffsschrift
  und andere Aufs\"atze'', Zweite Auflage 1971, I.\ Angelelli editor.
  Georg Olms Verlag Hildesheim, New York.
  
  \NAME{P.\ J.\ Freyd} 1972: Aspects of Topoi. 
  \emph{Bull.\ Australian Math.\ Soc.} \textbf{7,} 1-76.
  
  \NAME{K.\ G\"odel} 1931: \"Uber formal unentscheidbare S\"atze der
  Principia Mathematica und verwandter Systeme I. 
  \emph{Monatsh.\ der Mathematik und Physik} 38, 173-198.
  

  \NAME{R.\ L.\ Goodstein} 1971: \emph{Development of Mathematical
  Logic,} ch. 7: Free-Variable Arithmetics.  Logos Press.
  
  \NAME{F.\ Hausdorff} 1908: 
  Grundz\"uge einer Theorie der geordneten Mengen. 
  \emph{Math.\ Ann.} \textbf{65}, 435-505. 



  \NAME{D.\ Hilbert}: Mathematische Probleme. Vortrag Paris 1900. 
  \emph{Gesammelte Abhandlungen.} 
  Springer 
  1970. 
  
  
  \NAME{P.\ T.\ Johnstone} 1977: \emph{Topos Theory.} Academic Press

  \NAME{A.\ Joyal} 1973: Arithmetical Universes. Talk at Oberwolfach.
  
  
  \NAME{J.\ Lambek, P.\ J.\ Scott} 1986: \emph{Introduction to higher order 
  categorical logic.} Cambridge University Press.
  
  
  \NAME{F.\ W.\ Lawvere} 1964: An Elementary Theory of the Category of
  Sets. \emph{Proc.\ Nat.\ Acad.\ Sc.\ USA} \textbf{51,} 1506-1510.



  \NAME{S.\ Mac Lane} 1972: \emph{Categories for the working mathematician.} 
  Springer.
  
  
  \NAME{B.\ Pareigis} 1969: \emph{Kategorien und Funktoren}. Teubner.
  
  \NAME{R.\ P\'eter} 1967: \emph{Recursive Functions}. Academic Press.

  \NAME{M.\ Pfender} 1974: Universal Algebra in S-Monoidal Categories.
  Algebra-Berichte Nr. 20, Mathematisches Institut der Universit\"at
  M\"unchen. Verlag Uni-Druck M\"unchen.






  \NAME{M.\ Pfender} 2008: 
  Theories of PR Maps and Partial PR Maps. pdf file. 
  Condensed version as RCF\,1: Theories of PR Maps and Partial PR Maps. 
  arXiv: 0809.3676v1 [math.CT] 22 Sep 2008.
  
  \NAME{M.\ Pfender:} Evaluation and Consistency, Summary and 
  section 1 of version 1 of present work: arXiv 0809.3881v1 [math.CT] 
  23 Sep 2008. 

  \NAME{M.\ Pfender, M.\ Kr\"oplin, D.\ Pape} 1994: Primitive
  Recursion, Equality, and a Universal Set. 
  \emph{Math.\ Struct.\ in Comp.\ Sc.\ } \textbf{4,} 295-313.
  
  


  \NAME{W.\ Rautenberg} 1995/2006: \emph{A Concise Introduction to 
  Mathematical Logic.} Universitext Springer 2006.

  \NAME{R.\ Reiter} 1980: Mengentheoretische Konstruktionen in arithmetischen
  Universen. Diploma Thesis. Techn.\ Univ.\ Berlin.

  
  \NAME{L.\  Rom\`an} 1989: Cartesian categories with natural numbers object.
  \emph{J.\  Pure and Appl.\  Alg.} \textbf{58,} 267-278.
  


  \NAME{C.\ Smorynski} 1977: The Incompleteness Theorems. Part D.1
  in \NAME{Barwise} ed. 1977.

  \NAME{W.\ W.\ Tait} 1996: Frege versus Cantor and Dedekind: on the concept
  of number. Frege, Russell, Wittgenstein: \emph{Essays in Early Analytic
  Philosophy (in honor of Leonhard Linsky)} (ed. W.\ W.\ Tait). Lasalle:
  Open Court Press (1996): 213-248. Reprinted in \emph{Frege: Importance 
  and Legacy} (ed. M.\ Schirn). Berlin: Walter de Gruyter (1996): 70-113.


  \NAME{A.\ Tarski, S.\ Givant} 1987: \emph{A formalization of set theory 
  without variables}. AMS Coll.\ Publ.\ vol.\ 41.

\bigskip

\bigskip


  \noindent Address of the author: \\
  \NAME{M. Pfender}                       \hfill D-10623 Berlin \\
  Institut f\"ur Mathematik                              \\
  Technische Universit\"at Berlin         \hfill pfender@math.TU-Berlin.DE\\
  

\vfill

  \end{document}